\documentclass[a4paper,onecolumn, superscriptaddress,10pt,accepted=2019-07-24,issue=2, volume=1, shorttitle=papers/compositionality-1-2]{compositionalityarticle}
\pdfoutput=1
\usepackage[utf8]{inputenc}
\usepackage[T1]{fontenc}
\usepackage{amsmath,amsthm}
\usepackage{hyperref}

\usepackage{tikz}
\usepackage{lipsum}

\usepackage[numbers,sort&compress]{natbib}
\usepackage{amssymb,mathtools}
\usepackage{etex} 
\usepackage{stmaryrd}
\usepackage{centernot}
\usepackage{fancyhdr}
\usepackage[normalem]{ulem}
\usepackage{longtable}
\usepackage{calligra} 
\usepackage{dsfont}
\usepackage{enumerate}
\usepackage[all]{xy} 
\usepackage{tikz}
\usetikzlibrary{decorations.markings}
\usetikzlibrary{circuits.ee.IEC}
\usepackage{transparent}
\input diagxy

\usepackage{url}

\makeatletter
\renewcommand*\env@matrix[1][\arraystretch]{%
  \edef\arraystretch{#1}%
  \hskip -\arraycolsep
  \let\@ifnextchar\new@ifnextchar
  \array{*\c@MaxMatrixCols c}}
\makeatother

\numberwithin{equation}{section}

\theoremstyle{theorem}
\newtheorem{theorem}[equation]{Theorem}
\newtheorem{postulate}[equation]{Postulate}
\newtheorem{corollary}[equation]{Corollary}
\newtheorem{proposition}[equation]{Proposition}
\newtheorem{lemma}[equation]{Lemma}
\theoremstyle{definition}
\newtheorem{definition}[equation]{Definition}
\newtheorem{construction}[equation]{Construction}
\theoremstyle{definition}
\newtheorem{examplex}[equation]{Example}
\newenvironment{example}
  {\pushQED{\qed}\renewcommand{\qedsymbol}{$\blacklozenge$}\examplex}
  {\popQED\endexamplex}
\newtheorem{remarkx}[equation]{Remark}
\newenvironment{remark}
  {\pushQED{\qed}\renewcommand{\qedsymbol}{$\blacktriangle$}\remarkx}
  {\popQED\endremarkx}

\newtheorem{notation}[equation]{Notation}
\renewcommand\qedsymbol{$\blacksquare$}

\newcommand\define[1]{\emph{\textbf{#1}}}

\newcommand\llangle{\langle\!\langle}
\newcommand\rrangle{\rangle\!\rangle}

\newcommand{\stoch}{\;\xy0;/r.25pc/:(-3,0)*{}="1";(3,0)*{}="2";{\ar@{~>}"1";"2"|(1.09){\hole}};\endxy\!}
\newcommand{\xstoch}[1]{\;\xy0;/r.25pc/:(-5,0)*{}="1";(5,0)*{}="2";{\ar@{~>}"1";"2"|(1.09){\hole}^{#1}};\endxy\!}

\let\a=\alpha \let\b=\beta \let\g=\gamma \let\de=\delta \let\e=\epsilon
\let\z=\zeta \let\h=\eta \let\q=\theta  
\let\l=\lambda \let\r=\rho
\let\s=\sigma \let\t=\tau   \let\c=\chi 
\let\w=\omega         
  \let\S=\Sigma \let\U=\Upsilon  
\let\C=\Chi \let\W=\Omega
\def\vf{\varphi}

\newcommand{\be}{\begin{equation}}
\newcommand{\ee}{\end{equation}}
\def\ba{\begin{align}} 
\def\ea{\end{align}}
\newcommand{\bea}{\begin{eqnarray}}
\newcommand{\eea}{\end{eqnarray}}
\newcommand{\bx}{\begin{example}}
\newcommand{\ex}{\end{example}}
\newcommand{\bex}{\begin{exercise}}
\newcommand{\eex}{\end{exercise}}
\newcommand{\ban}{\begin{answer}}
\newcommand{\ean}{\end{answer}}
\newcommand{\bt}{\begin{theorem}}
\newcommand{\et}{\end{theorem}}
\newcommand{\bc}{\begin{corollary}}
\newcommand{\ec}{\end{corollary}}
\newcommand{\blem}{\begin{lemma}}
\newcommand{\elem}{\end{lemma}}
\newcommand{\bp}{\begin{problem}}
\newcommand{\ep}{\end{problem}}
\newcommand{\bn}{\begin{proposition}}
\newcommand{\en}{\end{proposition}}
\newcommand{\bd}{\begin{definition}}
\newcommand{\ed}{\end{definition}}
\newcommand{\bcon}{\begin{construction}}
\newcommand{\econ}{\end{construction}}
\newcommand{\bq}{\begin{question}}
\newcommand{\eq}{\end{question}}
\newcommand{\bprf}{\begin{proof}}
\newcommand{\eprf}{\end{proof}}
\newcommand{\br}{\begin{remark}}
\newcommand{\er}{\end{remark}}
\newcommand{\bs}{\begin{solution}}
\newcommand{\es}{\end{solution}}
\newcommand{\beqs}{\begin{eqnarray}}
\newcommand{\eeqs}{\end{eqnarray}}

\newcommand{\ov}[1]{\overline{#1}} 

\newcommand{\<}{\langle}
\renewcommand{\>}{\rangle}

\newcommand{\id}{\mathrm{id}}

\newcommand{\End}{\mathrm{End}}

\newcommand{\mC}{\mathcal{C}}
\newcommand{\mD}{\mathcal{D}}

\newcommand{\aand}{\qquad \& \qquad}
\newcommand{\tr}{{\rm tr} }

\def\C{{{\mathbb C}}}

\def\N{{{\mathbb N}}}

\def\Hi{{{\mathcal{H}}}}
\def\mI{{{\mathcal{I}}}}
\def\mJ{{{\mathcal{J}}}}
\def\mK{{{\mathcal{K}}}}
\def\mL{{{\mathcal{L}}}}
\def\mM{{{\mathcal{M}}}}
\def\mN{{{\mathcal{N}}}}
\newcommand{\CAlg}{\mathbf{C^*\text{-}Alg}}
\def\mA{{{\mathcal{A}}}}

\def\mB{{{\mathcal{B}}}}

\newcommand{\tCAT}{\mathcal{CAT}}

\newcommand{\op}{\mathrm{op}}
\newcommand{\St}{\mathbf{States}}

\newcommand{\rest}{\mathbf{rest}}
\newcommand{\pRep}{\mathbf{Rep}^{\bullet}}

\newcommand{\pGNS}{\mathbf{GNS}^{\bullet}}

\DeclareMathAlphabet{\mathcalligra}{T1}{calligra}{q}{n}
\DeclareFontShape{T1}{calligra}{q}{n}{<->s*[2.2]callig15}{}

\newcommand{\bigboxplus}{
  \mathop{
    \vphantom{\bigoplus} 
    \mathchoice
      {\vcenter{\hbox{\resizebox{\widthof{$\displaystyle\bigoplus$}}{!}{$\boxplus$}}}}
      {\vcenter{\hbox{\resizebox{\widthof{$\bigoplus$}}{!}{$\boxplus$}}}}
      {\vcenter{\hbox{\resizebox{\widthof{$\scriptstyle\oplus$}}{!}{$\boxplus$}}}}
      {\vcenter{\hbox{\resizebox{\widthof{$\scriptscriptstyle\oplus$}}{!}{$\boxplus$}}}}
  }\displaylimits 
}

\title[Stinespring's construction as an adjunction]{Stinespring's construction as an adjunction}
\author{Arthur J. Parzygnat}
\affiliation{Institut des Hautes \'Etudes Scientifiques,
35 Routes de Chartres, 91440, Bures-sur-Yvette, France}
\email{parzygnat@ihes.fr}

\begin{document}
\maketitle


\begin{abstract}
Given a representation of a unital $C^*$-algebra $\mA$ on a Hilbert space $\Hi$, together
with a bounded linear map $V:\mK\to\Hi$ from some other Hilbert space, 
one obtains a completely positive map on $\mA$ via restriction using the adjoint action associated to $V$. 
We show this restriction forms a natural transformation from
a functor of $C^*$-algebra representations to a functor of completely positive maps.
We exhibit Stinespring's construction as a left adjoint of this restriction.
Our Stinespring adjunction provides a universal property associated to minimal Stinespring dilations and morphisms of Stinespring dilations. We use these results to prove the purification postulate for all finite-dimensional $C^*$-algebras.
\end{abstract} 


\section{Introduction and outline}
\label{sec:intro}

Given a unital $C^*$-algebra $\mA$ and a state (a positive
and unital linear functional) on $\mA$, the Gelfand--Naimark--Segal (GNS) construction produces a
cyclic representation of $\mA$ in a natural way
compatible with the operation that produces a state from a
representation together with a unit vector. 
This naturality has been expressed as an adjunction
\be
\label{eq:GNSrest}
\xy0;/r.25pc/:
(-25,0)*+{\CAlg^{\op}}="1";
(25,0)*+{\tCAT}="2";
(0,0)*{\dashv};
{\ar@/^1.75pc/"1";"2"^{\St}};
{\ar@/_1.75pc/"1";"2"_{\pRep}};
{\ar@{=>}(2.5,-5.5);(2.5,5.5)_{\rest}};
{\ar@{=>}(-2.5,5.5);(-2.5,-5.5)_{\pGNS}};
\endxy
\ee
in a certain 2-category of functors from
the (opposite of the) category of $C^*$-algebras to the 
2-category of locally small categories%
\footnote{Technically, the collection of objects of the 2-category of locally small categories does not form a class. See Sections 3.49, 3.50, and 3.51 in \cite{AdHeSt06} and Sections 8 and 16 in \cite{Sh08} for the set-theoretic foundations of categories of categories and functor categories. In the present work, our constructions will be explicitly defined on objects and morphisms, and these set-theoretic issues will not affect our results.}
 \cite{PaGNS}. 

Stinespring's construction can be viewed as a generalization
of the GNS construction
by replacing states with operator-valued completely positive (OCP) maps.
In the present work, we extend our GNS adjunction
in Theorem~\ref{thm:stinespring} to include such OCP maps,
showing that Stinespring's construction can also be viewed
as an adjunction 
\be
\label{eq:Stinerest}
\xy0;/r.25pc/:
(-25,0)*+{\CAlg^{\op}}="1";
(25,0)*+{\tCAT}="2";
(0,0)*{\dashv};
{\ar@/^1.75pc/"1";"2"^{\mathbf{OCP}}};
{\ar@/_1.75pc/"1";"2"_{\mathbf{AnRep}}};
{\ar@{=>}(2.5,-5.5);(2.5,5.5)_{\rest}};
{\ar@{=>}(-2.5,5.5);(-2.5,-5.5)_{\mathbf{Stine}}};
\endxy
\ee
in the same 2-category of functors.

There are several subtle differences between the two adjunctions
(\ref{eq:GNSrest}) and (\ref{eq:Stinerest}).
The most notable one is that the category of OCP maps 
for a given $C^*$-algebra
is no longer discrete as it is for the GNS construction. 
This is a consequence of relaxing the unitality assumption on the positive maps and allowing the maps to be operator-valued instead of complex-valued. 
Furthermore, for similar reasons, the category of anchored representations has far more morphisms than the category of pointed representations introduced in \cite{PaGNS}. 
Our result provides a 
categorical sense of the minimality of Stinespring's construction
including many of its functorial properties for all $C^*$-algebras. As a consequence, we apply our results to prove a version of the purification postulate, which has been used by Chiribella, D'Ariano, and Perinotti as a crucial postulate in isolating quantum theory among other operational probabilistic theories \cite{CDP10}. 

Our results are closely related to, but distinct from, the universal property of 
minimal Stinespring dilations 
described by Westerbaan and Westerbaan
in \cite{WeWe16} (see also Section~2.3 in~\cite{We19}). 
In \cite{WeWe16}, the authors describe a
universal property for (minimal) \emph{Paschke dilations} for \emph{normal}
completely positive maps between \emph{von Neumann algebras}.
In the present work, 
we provide an alternative universal property for Stinespring dilations without restricting ourselves to \emph{normal} completely positive maps, since we allow the domains of our maps to be arbitrary $C^*$-algebras. 
Therefore, we characterize \emph{all} minimal Stinespring dilations by our universal property. 
However, our codomain for a completely positive map is assumed to be bounded operators on some Hilbert space, while Westerbaan and Westerbaan assume an abstract von Neumann algebra. Therefore, our results both provide a similar universal property for minimal Stinespring dilations of normal completely positive maps, but neither result subsumes the other. It would be interesting to see if there is a common generalization of both our results that includes Paschke dilations of arbitrary completely positive maps between arbitrary $C^*$-algebras. 

The outline of our paper is as follows. In Section~\ref{sec:operatorstates}, relevant background for 
completely positive maps is provided. An appropriate category of OCP maps for a $C^*$-algebra is defined and briefly contrasted with the
category of states introduced in \cite{PaGNS}. 
In Section~\ref{sec:anchoredrepns}, the category of anchored representations
for a $C^*$-algebra is introduced. Anchored representations are
generalizations of pointed representations described in \cite{PaGNS}.
Restriction from anchored representations to OCP maps
is defined in Section~\ref{sec:rest} and is shown to form a natural
transformation. Section~\ref{sec:Stinespring} contains our main result, Theorem~\ref{thm:stinespring}, which states that 
the restriction natural transformation has a left adjoint whose ingredients are determined by Stinespring's construction. 
Section~\ref{sec:generalizations} provides some immediate consequences of our adjunction, including its universal property. 
We provide a comparison to 
minimal Stinespring dilations in  
Corollary~\ref{cor:minimalStinespring}. 
The rest of Section~\ref{sec:generalizations} introduces partial isometries and their role in our Stinespring adjunction. For example, Theorem~\ref{thm:connectingStinespringrepns} shows that all morphisms between Stinespring dilations extend a minimal intertwining partial isometry. 
The first part of Section~\ref{sec:comparetoGNS}
compares our main theorem to the GNS adjunction from \cite{PaGNS}.
We then give a reformulation of the standard purification postulate suitable for our purposes in Postulate~\ref{thm:purificationpostulateprocess}. 
What follows is a detailed analysis of how this postulate is a consequence of our results. This leads to Theorem~\ref{thm:ChiribellaPurification} and Corollary~\ref{cor:fdcalgpurification}, which establish the essential uniqueness of purifications for finite-dimensional $C^*$-algebras. 
An index of notation is included  in Appendix~\ref{index} for the reader's convenience. 
Appendix~\ref{appendix} briefly reviews oplax-natural transformations, modifications, and 2-categorical adjunctions. Proposition~\ref{prop:adjunctionsinfunctorcategories} provides an equivalence in data regarding adjunctions in 2-categories of functors and their evaluations on objects. 

In this article, all $C^*$-algebras and $^*$-homomorphisms 
will be taken to be unital unless otherwise specified. 
A $^*$-homomorphism of $C^*$-algebras 
will be denoted diagrammatically as a straight arrow $\to$,
while a linear map (often a completely positive map) of $C^*$-algebras 
will be denoted as a curvy arrow $\stoch.$ 
This is largely motivated by the equivalence between the 
(opposite of the) category
of compact Hausdorff spaces with stochastic maps 
and the category of commutative $C^*$-algebras with positive maps \cite{Pa17}, \cite{FuJa13}.
This equivalence restricts to the usual commutative Gelfand--Naimark theorem
that describes the equivalence between the 
(opposite of the) category
of compact Hausdorff spaces with continuous functions 
and the category of commutative $C^*$-algebras with $^*$-homomorphisms \cite{GN43}. To the best of the author's knowledge, 
the usage of the $\stoch$ notation originated in the work
of Baez and Fritz on relative entropy \cite{BaFr14}. 
In what follows, ``iff'' stands for ``if and only if'' and is used solely in definitions.
The symbols $\text{\qedsymbol}$, $\blacklozenge$, and $\blacktriangle$ signal the end of a proof, example, and remark, respectively. 
$\mC_{i}$ denotes the $i$-morphisms of the category (or 2-category) $\mC$ 
(when $i=0,$ these refer to the objects of the category),
$\CAlg$ denotes the category of $C^*$-algebras and $^*$-homomorphisms, 
$\tCAT$ denotes the 2-category of categories, functors, and natural transformations, and  
$\mathbf{Fun}(\CAlg^{\op},\tCAT)$ denotes the 2-category of 
functors $\CAlg^{\op}\to\tCAT,$ oplax-natural transformations of such functors, 
and modifications of oplax-natural transformations. 
Although Appendix~\ref{appendix} reviews what is sufficient for our purposes, the reader is referred to 
our previous article on the GNS construction for further details
\cite{PaGNS}. 
Any claims made without proof follow easily from the definitions---this is particularly the case in Sections~\ref{sec:operatorstates},~\ref{sec:anchoredrepns},~and~\ref{sec:rest}.

\section{Operator-valued completely positive maps on $C^*$-algebras}
\label{sec:operatorstates}

The notion of a completely positive (CP) map will be used throughout
and will therefore be briefly reviewed.
Most of the concepts used here are introduced in the first few chapters
of Paulsen~\cite{Pa03}. 
General facts used without mention regarding Hilbert spaces can be found in Chapter~4 of Rudin~\cite{Ru87} and Chapter~12 of Rudin~\cite{Ru91}. 

\begin{notation}
\label{not:intro}
The set of natural numbers $1,2,3,\dots$ is denoted by $\N$. 
If $\mK$ is a Hilbert space, $\mB(\mK)$ denotes the $C^*$-algebra of bounded operators
on $\mK.$ If $a\in\mB(\mK),$ then $a^*$ denotes the adjoint of $a.$ 
If  $\mA$ is a $C^*$-algebra, $1_{\mA}$ denotes the unit in $\mA.$ 
The norm on $C^*$-algebras 
is always written using $\lVert\;\cdot\;\rVert$ without any subscripts,
while norms and inner products on Hilbert spaces frequently have subscripts. For example, $\<\;\cdot\;,\;\cdot\;\>_{\mK}$ denotes the inner product
on $\mK$ with linearity in the right variable and conjugate
linearity in the left variable. 
For $n\in\N,$ $\mM_{n}(\mA)$ denotes the $C^*$-algebra of $n\times n$ matrices
with entries in $\mA$. Addition and multiplication in $\mathcal{M}_{n}(\mA)$ are defined as for ordinary matrices. If $a_{ij}$ is the $ij$-th entry of $A\in\mathcal{M}_{n}(\mA),$ then the $ij$-th entry of $A^*$ is $a_{ji}^{*}$.
The norm on $\mM_{n}(\mA)$ is used only once in this article
(cf.~(\ref{eq:Stinespringbounded}))
and the reader is referred to the first several pages in 
Chapter~1 of Paulsen~\cite{Pa03} for details.
In particular, $\mM_{n}(\C)$ denotes the $C^*$-algebra of 
$n\times n$ complex matrices. In this case, the unit of 
$\mM_{n}(\C)$ will be denoted by $\mathds{1}_{n}.$ Any $C^*$-algebra of the form $\mathcal{M}_{n}(\C)$ will be referred to as a matrix algebra. 
\end{notation}

\bd
\label{defn:cpmaps}
Let $\mA$ and $\mB$ be $C^*$-algebras, 
and let $\vf:\mA\stoch\mB$ be a linear map.
The map $\vf$ is said to be \define{unital} iff $\vf(1_{\mA})=1_{\mB}.$
An element of a $C^*$-algebra $\mA$ is \define{positive}
iff it equals $a^*a$ for some $a\in\mA$. A linear map
between $C^*$-algebras is \define{positive} iff it is linear and 
it sends positive elements to positive elements. 
A positive map into $\C$ is referred to as a \define{positive linear functional}. A positive unital linear functional is called a \define{state}. 
If $n\in\N$, the \define{$n$-ampliation} of $\vf:\mA\stoch\mB$ is the linear map 
$\vf_{n}:\mM_{n}(\mA)\stoch\mM_{n}(\mB)$ defined by 
the assignment
\be
\label{eq:ampliation}
\mM_{n}(\mA)\ni\begin{bmatrix}a_{11}&\cdots&a_{1n}\\\vdots&&\vdots\\a_{n1}&\cdots&a_{nn}\end{bmatrix}
\xmapsto{\vf_{n}}
\begin{bmatrix}\vf(a_{11})&\cdots&\vf(a_{1n})\\\vdots&&\vdots\\\vf(a_{n1})&\cdots&\vf(a_{nn})\end{bmatrix}.
\ee
The map $\vf:\mA\stoch\mB$ 
is said to be \define{$n$-positive}
iff its $n$-ampliation is positive. 
The map $\vf:\mA\stoch\mB$ 
is said to be \define{completely positive}
iff it is $n$-positive for all $n\in\N.$ 
The shorthand PU, CP, and CPU may be used to denote a 
positive unital, completely positive, or completely positive unital map, 
respectively. 
\ed

When $\mB=\mB(\mK),$ bounded operators on a Hilbert space $\mK,$
the latter matrix in (\ref{eq:ampliation}) 
acts on an $n$-tuple of vectors in $\mK,$ i.e.~elements of $\mK\oplus\cdots\oplus\mK$ ($n$ times), to provide
another $n$-tuple of vectors in $\mK$ 
(one chooses an ordering for the direct sum here). 

\bx
\label{ex:cpmaps}
Some examples of CP maps follow.
\begin{enumerate}[(a)]
\itemsep0pt
\item
All $^*$-homomorphisms between $C^*$-algebras are CPU. 
\item
Let $T:\mK\to\mL$ be a bounded linear map between Hilbert spaces. 
Then, the map
\be
\label{eq:adjointactionmap}
\begin{split}
\mB(\mK)&\xstoch{\mathrm{Ad}_{T}}\mB(\mL)\\
A&\xmapsto{\;\;\qquad} TAT^*
\end{split}
\ee
is CP. 
This map is sometimes referred to as the \emph{adjoint action map}.
It is CPU whenever $T^*$ is an isometry. This notation is chosen so that if $S:\mathcal{L}\to\mathcal{M}$ is another bounded linear map between Hilbert spaces, then $\mathrm{Ad}_{S}\circ\mathrm{Ad}_{T}=\mathrm{Ad}_{ST}.$ 
\item 
Non-negative linear combinations of maps as in (\ref{eq:adjointactionmap}) are also CP. 
\end{enumerate}
Several more examples will appear later in this work. 
\ex

\blem
\label{lem:compCPisCP}
Each of the three classes of maps between $C^*$-algebras from Definition \ref{defn:cpmaps} (positive, unital, and CP) is closed under composition. 
All positive maps $\vf:\mA\stoch\mB$ between $C^*$-algebras
are self-adjoint in the sense that $\vf(a^*)=\vf(a)^*$ for all $a\in\mA.$ 
\elem

\bprf
See Exercise~2.1 in Paulsen~\cite{Pa03} for the last statement.
\eprf

\bd
\label{defn:opstate}
Let $\mA$ be a $C^*$-algebra. An \define{operator-valued CP map} (OCP map) on $\mA$ 
is a pair $(\mK,\vf)$ consisting of a Hilbert space $\mK$ 
and a CP map
$\vf:\mA\stoch\mB(\mK).$ 
When $\vf$ is unital, it is called an \define{operator state} on $\mA$ \cite{Ef78}. 
Let $(\mL,\psi)$ be another OCP map on $\mA.$ 
A \define{morphism of OCP maps}  $T:(\mK,\vf)\to(\mL,\psi)$ 
is a bounded linear map 
$T:\mK\to\mL$ such that 
\be
\label{eq:opstatesdiagram}
\xy0;/r.25pc/:
(-12.5,7.5)*+{\mK}="K1";
(12.5,7.5)*+{\mK}="K2";
(-12.5,-7.5)*+{\mL}="L1";
(12.5,-7.5)*+{\mL}="L2";
{\ar"K1";"K2"^{\vf(a)}};
{\ar"L1";"L2"_{\psi(a)}};
{\ar"K1";"L1"_{T}};
{\ar"K2";"L2"^{T}};
\endxy
\ee
commutes for all $a\in\mA.$ If $(\mK,\vf)$ and $(\mL,\psi)$ are operator states, a \define{morphism of operator states} is a morphism of OCP maps with $T$ an isometry. 
\ed

\br
\label{rmk:morphismopstates}
The choice of morphisms of OCP maps is subtle. 
One might have tried any number of reasonable variants.
For instance, one might request that the diagram
\be
\label{eq:opstatemorphismother}
\xy0;/r.25pc/:
(0,7.5)*+{\mA}="A";
(-12.5,-7.5)*+{\mB(\mK)}="K";
(12.5,-7.5)*+{\mB(\mL)}="L";
{\ar@{~>}"A";"K"_{\vf}|(0.84){\hole}};
{\ar@{~>}"A";"L"^{\psi}|(0.84){\hole}};
{\ar@{~>}"K";"L"_{\mathrm{Ad}_{T}}|(0.80){\hole}};
\endxy
\ee
commutes. Another possibility is that the diagram 
\be
\label{eq:opstatemorphismother2}
\xy0;/r.25pc/:
(0,7.5)*+{\mA}="A";
(-12.5,-7.5)*+{\mB(\mK)}="K";
(12.5,-7.5)*+{\mB(\mL)}="L";
{\ar@{~>}"A";"K"_{\vf}|(0.84){\hole}};
{\ar@{~>}"A";"L"^{\psi}|(0.84){\hole}};
{\ar@{~>}"L";"K"^{\mathrm{Ad}_{T^*}}|(0.80){\hole}};
\endxy
\ee
commutes. 
When $T$ is an isometry, 
commutativity of (\ref{eq:opstatemorphismother}) 
implies commutativity of (\ref{eq:opstatesdiagram}) and commutativity of 
(\ref{eq:opstatesdiagram}) 
implies commutativity of
(\ref{eq:opstatemorphismother2}). When $T^*$ is an isometry, these implications are reversed. In particular, when $T$ is unitary, all three are equivalent.
To see the first implication when $T$ is an isometry, diagram~(\ref{eq:opstatemorphismother}) says 
$T\vf(a)T^*=\psi(a)$ for all $a\in\mA.$ 
Applying $T$ on the right gives 
$T\vf(a)=\psi(a)T$ because $T^*T=\id_{\mK}.$ 
The second implication follows by applying $T^*$ on the left of this result. 
In summary (\ref{eq:opstatemorphismother})~$\implies$~(\ref{eq:opstatesdiagram})~$\implies$~(\ref{eq:opstatemorphismother2}) for operator states. 
Explicit counter-examples showing the reverse implications
fail in general are provided in Examples~\ref{ex:implication1irreversible}~and~\ref{ex:implication2irreversible}. 
We will see that 
commutativity of (\ref{eq:opstatemorphismother})
is too strong a requirement and 
commutativity of (\ref{eq:opstatemorphismother2})
is too weak a requirement for the purposes sought out in this work
(these points will be explained in footnotes).
No such simple comparison can be made for OCP maps.%
\footnote{This is easy to see for trivial reasons. For example, if $(\mK,\vf)\xrightarrow{T}(\mL,\psi)$ is a morphism of OCP maps and $\l\in[0,1)\cup(1,\infty)$, then $(\mK,\l\vf)\xrightarrow{\l T}(\mL,\l\psi)$ is a morphism of OCP maps. However, neither 
(\ref{eq:opstatemorphismother}) nor (\ref{eq:opstatemorphismother2}) commute in general. Along similar lines, (\ref{eq:opstatemorphismother}) and (\ref{eq:opstatemorphismother2}) do not imply each other (see also Example~\ref{ex:tracial}). 
}
These subtle points and the universal properties that will be discussed in this work would have been missed if we demanded $T$ to be unitary. 
\er

\bx
\label{ex:states}
As a special case, set $\mK=\C.$
An OCP map ${\vf:\mA\stoch\mB(\C)\cong\C}$ is (naturally isomorphic to) a positive linear functional.
If $\psi:\mA\stoch\C$ is another such OCP map, then a morphism
$T:(\C,\vf)\to(\C,\psi)$ of OCP maps consists of a linear map $T:\C\to\C.$ 
Such a linear map must be of the form 
$T(\l)=z\l$ for all $\l\in\C$ for some unique $z\in\C.$ 
If $z\ne0$, since $\vf$ and $\psi$ are linear, the diagram 
\be
\xy0;/r.25pc/:
(-12.5,7.5)*+{\C}="K1";
(12.5,7.5)*+{\C}="K2";
(-12.5,-7.5)*+{\C}="L1";
(12.5,-7.5)*+{\C}="L2";
{\ar"K1";"K2"^{\vf(a)}};
{\ar"L1";"L2"_{\psi(a)}};
{\ar"K1";"L1"_{T}};
{\ar"K2";"L2"^{T}};
\endxy
\ee
commutes for all $a\in\mA$ if and only if $\vf=\psi.$ 
Notice that when $T$ is an isometry, it is of the form 
$T(\l)=e^{i\q}\l$ for some $\q\in[0,2\pi).$ In this case, (\ref{eq:opstatemorphismother2})
and (\ref{eq:opstatemorphismother}) 
(and hence (\ref{eq:opstatesdiagram}) as well) 
are equivalent because $T$ is also unitary. 
\ex

\bx
\label{ex:tracial}
Fix $m,p\in\N$ and define the \define{tracial map}
\be
\label{eq:tracialmap}
\begin{split}
\mM_{m}(\C)&\xstoch{\tau}\mM_{p}(\C)\\
A&\xmapsto{\;\;\qquad}\frac{1}{m}\tr(A)\mathds{1}_{p},
\end{split}
\ee
where $\tr$ denotes the standard (un-normalized) trace. 
A quick computation shows $\t$ is unital. To see that it is CP, first note that since the trace is positive, 
the trace is CP because every positive map 
into $\C$ is automatically CP
(cf.~Theorem~3 in Stinespring~\cite{St55}). 
Second, note that $\t$ equals the composite of CPU maps
\be
\mathcal{M}_{m}(\C)\xstoch{\frac{1}{m}\tr}\C\xrightarrow{\;\;!\;\;}\mM_{p}(\C),
\ee 
where the second
map is the unique unital linear map sending $1$ to $\mathds{1}_{p}.$ Therefore, $\t$ is CPU by Example~\ref{ex:cpmaps} and Lemma~\ref{lem:compCPisCP}. Given $q\in\N,$ let $\s:\mathcal{M}_{m}(\C)\stoch\mathcal{M}_{q}(\C)$ denote the tracial map for these dimensions and let $T:\C^{p}\to\C^{q}$ be \emph{any} linear transformation. Then $(\C^{p},\t)\xrightarrow{T}(\C^{q},\s)$ is a morphism of OCP maps.
\ex

\bx
\label{ex:implication1irreversible}
The present example will show that 
(\ref{eq:opstatesdiagram})~$\centernot\implies$~(\ref{eq:opstatemorphismother})  for morphisms of operator states. 
Let $\mA:=\mathcal{M}_{2}(\C)$ and set
\begin{align}
\begin{aligned}
\mA&\xstoch{\vf}\C&\qquad\mA&\xstoch{\psi}\mathcal{M}_{2}(\C)&\qquad\C&\xrightarrow{\;\;\;\;T\;\;\;\;}\C^{2}\\
A&\xmapsto{\;\;\qquad}\frac{1}{2}\tr(A)&A&\xmapsto{\;\;\qquad}\frac{1}{2}\tr(A)\mathds{1}_{2}&\l&\xmapsto{\;\;\qquad}\frac{\l}{\sqrt{2}}\begin{bmatrix}1\\1\end{bmatrix}.
\end{aligned}
\end{align}
Then $\vf$ and $\psi$ are operator states and $T:(\C,\vf)\to(\C^{2},\psi)$ defines a morphism of 
operator states  by Example~\ref{ex:tracial}. Note that 
$T^*:\C^{2}\to\C$ is given by 
\be
T^*\left(\begin{bmatrix}x\\y\end{bmatrix}\right)
=\frac{x+y}{\sqrt{2}}.
\ee
Therefore, 
\be
\Big(\mathrm{Ad}_{T}\big(\vf(A)\big)\Big)\left(\begin{bmatrix}x\\y\end{bmatrix}\right)
=\frac{1}{4}\tr(A)(x+y)\begin{bmatrix}1\\1\end{bmatrix}
\ee
while
\be
\big(\psi(A)\big)\left(\begin{bmatrix}x\\y\end{bmatrix}\right)
=\frac{1}{2}\tr(A)\begin{bmatrix}x\\y\end{bmatrix},
\ee
which shows (\ref{eq:opstatemorphismother}) does not commute.
\ex

\bx
\label{ex:implication2irreversible}
The present example will show that 
(\ref{eq:opstatemorphismother2}) $\centernot\implies$ (\ref{eq:opstatesdiagram}).
Let $\mA:=\mathcal{M}_{2}(\C)$ and set
\begin{align}
\begin{aligned}
\mA&\xstoch{\vf}\C&\qquad\mA&\xstoch{\psi}\mathcal{M}_{2}(\C)&\qquad\C&\xrightarrow{\;\;\;\;T\;\;\;\;}\C^{2}\\
A&\xmapsto{\;\;\qquad}\frac{1}{2}\tr(A)&A&\xmapsto{\;\;\qquad}\frac{1}{2}A+\frac{1}{2}\begin{bmatrix}0&1\\1&0\end{bmatrix}A\begin{bmatrix}0&1\\1&0\end{bmatrix}&\l&\xmapsto{\;\qquad}\begin{bmatrix}\l\\0\end{bmatrix}.
\end{aligned}
\end{align}
Then $\vf$ and $\psi$ are operator states (cf.~Examples~\ref{ex:tracial}~and~\ref{ex:cpmaps}). Note that $T^{*}:\C^{2}\to\C$ is given by 
\be
T^*\left(\begin{bmatrix}x\\y\end{bmatrix}\right)=x
\ee
and
\be
\psi\left(\begin{bmatrix}a&b\\c&d\end{bmatrix}\right)=\frac{1}{2}\begin{bmatrix}a+d&b+c\\b+c&a+d\end{bmatrix}.
\ee
Therefore, although (\ref{eq:opstatemorphismother2}) holds, 
\be
\left(T \vf\left(\begin{bmatrix}a&b\\c&d\end{bmatrix}\right)\right)(\l)=\frac{\l}{2}\begin{bmatrix}a+d\\0\end{bmatrix}
\;\text{ and }\;
\left(\psi\left(\begin{bmatrix}a&b\\c&d\end{bmatrix}\right) T\right)(\l)=\frac{\l}{2}\begin{bmatrix}a+d\\b+c\end{bmatrix}
\ee
show that condition (\ref{eq:opstatesdiagram}) fails.
\ex

The following proposition provides 
the general structure
of morphisms of operator states. 

\bn
\label{prop:morphismofstatesdirectsum}
Let $\vf:\mA\stoch\mB(\mK)$ and $\psi:\mA\stoch\mB(\mathcal{L})$ be operator states on a $C^*$-algebra $\mA$ with $\mathcal{K}$ and $\mathcal{L}$ Hilbert spaces.
A morphism $(\mK,\vf)\xrightarrow{T}(\mL,\psi)$ of operator states exists if and only if
there exist closed subspaces $\mathcal{L}_{1},\mathcal{L}_{2}\subseteq\mathcal{L},$   
operator states $\psi_{j}:\mA\stoch\mB(\mathcal{L}_{j})$ for $j\in\{1,2\},$ and a unitary 
map $U:\mathcal{K}\to\mathcal{L}_{1}$ such that 
$\mathcal{L}=\mathcal{L}_{1}\oplus\mathcal{L}_{2},$
$U\vf(a) U^{*}=\psi_{1}(a),$ and 
$\psi(a)=\psi_{1}(a)\oplus\psi_{2}(a)$ for all $a\in\mA.$ 
\en

\bprf
{\color{white}{you found me!}}

\noindent
($\Rightarrow$)
For the forward direction, set $\mathcal{L}_{1}:=T(\mathcal{K}).$
Since $T$ is an isometry, $\mathcal{L}_{1}$ is a closed subspace of $\mathcal{L}$. 
Set $\mathcal{L}_{2}:=\mathcal{L}_{1}^{\perp},$ the orthogonal complement of $\mathcal{L}_{1}$ inside $\mathcal{L},$
$U:=\pi_{1} T,$ and $\psi_{j}(a):=\pi_{j}\psi(a) i_{j}$ for all $a\in\mA,$ where $\pi_{j}:\mathcal{L}\to\mathcal{L}_{j}$ denotes the projection onto the $j$-th factor
and $i_{j}:\mathcal{L}_{j}\to\mathcal{L}$ denotes the inclusion of the $j$-th factor.
Since $T$ is an isometry, $U$ is unitary. 
Since $T$ is a morphism of operator states, $\psi(a)\mathcal{L}_{1}\subseteq\mathcal{L}_{1}$ for all
$a\in\mA.$ Furthermore, since $\psi$ is positive, $\psi(a)^{*}=\psi(a^*)$ so that $T^*\psi(a^*)=\vf(a^*)T^*$ for all $a\in\mA$ upon taking the adjoint of (\ref{eq:opstatesdiagram}). Since $^*$ is an involution on $\mA,$ this is equivalent to $T^*\psi(a)=\vf(a)T^*$ for all $a\in\mA.$ Thus, $\psi(a)\mathcal{L}_{2}\subseteq\mathcal{L}_{2}.$ These two facts imply $\psi(a)=\psi_{1}(a)\oplus\psi_{2}(a)$ for all $a\in\mA.$ Finally, 
\be
\begin{split}
U\vf(a) U^{*}&=\pi_{1} T\vf(a) T^* i_{1}\quad\text{ by definition of $U$}\\
&=\pi_{1}TT^*\psi(a)TT^* i_{1}\quad\text{ by (\ref{eq:opstatesdiagram}) and Remark~\ref{rmk:morphismopstates}}\\
&=\pi_{1}\psi(a)i_{1}\quad\text{ since $T$ is an isometry onto $\mathcal{L}_{1}$}\\
&=\psi_{1}(a)\quad\text{ by definition of $\psi_{j}$}
\end{split}
\ee
for all $a\in\mA.$ 

\noindent
($\Leftarrow$)
To see the reverse direction, set $T:=i_{1}U$. 
Then $T$ is a morphism of operator states since the diagram
\be
\xy0;/r.25pc/:
(-12.5,15)*+{\mK}="K1";
(12.5,15)*+{\mK}="K2";
(-12.5,0)*+{\mL_{1}}="L1";
(12.5,0)*+{\mL_{1}}="L2";
(-12.5,-15)*+{\mL}="Lb1";
(12.5,-15)*+{\mL}="Lb2";
{\ar"K1";"K2"^{\vf(a)}};
{\ar"L1";"L2"|-{\psi_{1}(a)}};
{\ar"K1";"L1"^{U}};
{\ar"K2";"L2"_{U}};
{\ar"L1";"Lb1"^{i_{1}}};
{\ar"L2";"Lb2"_{i_{1}}};
{\ar"Lb1";"Lb2"_{\psi(a)}};
{\ar@/_2.0pc/"K1";"Lb1"_{T}};
{\ar@/^2.0pc/"K2";"Lb2"^{T}};
\endxy
\ee
commutes for all $a\in\mA.$ 
\eprf

OCP maps and operator states together with their morphisms form categories. 

\blem
\label{lem:opstA}
Let 
$(\mK,\vf)\xrightarrow{T}(\mL,\psi)\xrightarrow{S}(\mM,\chi)$ be a pair of
composable morphisms of OCP maps (or operator states) on a $C^*$-algebra $\mA$. The composite 
of said morphisms is defined to be 
the composite $ST$ of linear transformations and is 
a morphism of OCP maps (resp.~operator states). The identity 
$(\mK,\vf)\xrightarrow{\id_{(\mK,\vf)}}(\mK,\vf)$
is the identity linear transformation $\id_{\mK}.$
The collection of all OCP maps (or operator states) on $\mA$ with their morphisms forms a category, which is
denoted by $\mathbf{OCP}(\mA)$ (resp.~$\mathbf{OpSt}(\mA)$).  
Furthermore, $\mathbf{OpSt}(\mA)$ is a subcategory of $\mathbf{OCP}(\mA).$
\elem

\br
The result $T^*\psi(a)=\vf(a)T^*$ obtained in the proof of Proposition~\ref{prop:morphismofstatesdirectsum} shows $\mathbf{OCP}(\mA)$ is a $^*$-category. This will be explained in more detail in Proposition~\ref{prop:OCPandAnReparedagger}.
\er

\blem
\label{lem:opstf}
Let $f:\mA'\to\mA$ be a $^*$-homomorphism of $C^*$-algebras. 
The assignment
\be
\begin{split}
\mathbf{OCP}(\mA)&\xrightarrow{\mathbf{OCP}_{f}}\mathbf{OCP}(\mA')\\
(\mK,\vf)&\xmapsto{\qquad\;}(\mK,\vf\circ f)\\
\Big((\mK,\vf)\xrightarrow{T}(\mL,\psi)\Big)&\xmapsto{\qquad\;}\Big((\mK,\vf\circ f)\xrightarrow{T}(\mL,\psi\circ f)\Big)
\end{split}
\ee
defines a functor. Furthermore, it restricts to a functor $\mathbf{OpSt}(\mA)\xrightarrow{\mathbf{OpSt}_{f}}\mathbf{OpSt}(\mA')$.
\elem

\blem
\label{lem:opstfunctor}
The assignment 
\be
\begin{split}
\CAlg^{\op}&\xrightarrow{\mathbf{OCP}}\tCAT\\
\mA&\xmapsto{\;\;\quad\;\;}\mathbf{OCP}(\mA)\\
\left(\mA'\xrightarrow{f}\mA\right)&\xmapsto{\;\;\quad\;\;}\left(\mathbf{OCP}(\mA)\xrightarrow{\mathbf{OCP}_{f}}\mathbf{OCP}(\mA')\right)
\end{split}
\ee
defines a functor. The same is true for $\CAlg^{\op}\xrightarrow{\mathbf{OpSt}}\tCAT$. 
\elem

\section{Anchored representations of $C^*$-algebras}
\label{sec:anchoredrepns}

\bd
\label{defn:anchoredrepn}
Let $\mA$ be a $C^*$-algebra. An \define{anchored representation}
of $\mA$ is a quadruple $(\mK,\Hi,\pi,V)$ consisting of two Hilbert spaces $\Hi$ and $\mK,$ 
a $^*$-homomorphism $\pi:\mA\to\mB(\Hi),$ and 
a bounded linear map $V:\mK\to\Hi.$ 
When $V$ is an isometry, $(\mK,\Hi,\pi,V)$ is called a \define{preserving anchored representation}. 
Let $(\mL,\mI,\rho,W)$ be another anchored representation of $\mA.$ 
A \define{morphism of anchored representations}
$(T,L):(\mK,\Hi,\pi,V)\to(\mL,\mI,\rho,W)$ consists of bounded linear maps
$T:\mK\to\mL$ and $L:\Hi\to\mI$ satisfying the following two conditions. 
First, the diagram 
\be
\label{eq:intertwiner}
\xy0;/r.25pc/:
(-12.5,7.5)*+{\Hi}="H1";
(12.5,7.5)*+{\Hi}="H2";
(-12.5,-7.5)*+{\mI}="I1";
(12.5,-7.5)*+{\mI}="I2";
{\ar"H1";"H2"^{\pi(a)}};
{\ar"I1";"I2"_{\rho(a)}};
{\ar"H1";"I1"_{L}};
{\ar"H2";"I2"^{L}};
\endxy
\ee
commutes for all $a\in\mA$. Second, the diagrams
\be
\label{eq:morphismanchoredrepns}
\xy0;/r.25pc/:
(-12.5,7.5)*+{\mK}="K";
(12.5,7.5)*+{\Hi}="H";
(-12.5,-7.5)*+{\mL}="L";
(12.5,-7.5)*+{\mI}="I";
{\ar"K";"H"^{V}};
{\ar"L";"I"_{W}};
{\ar"K";"L"_{T}};
{\ar"H";"I"^{L}};
\endxy
\qquad\text{and}\qquad
\xy0;/r.25pc/:
(12.5,7.5)*+{\mK}="K";
(-12.5,7.5)*+{\Hi}="H";
(12.5,-7.5)*+{\mL}="L";
(-12.5,-7.5)*+{\mI}="I";
{\ar"H";"K"^{V^*}};
{\ar"I";"L"_{W^*}};
{\ar"K";"L"^{T}};
{\ar"H";"I"_{L}};
\endxy
\ee
both commute. When $V,W,T,$ and $L$ are isometries, then $(T,L)$ is said to be a \define{morphism of preserving anchored representations}.
\ed

\br
\label{rmk:unitarymorphismofpanrep}
For preserving anchored representations, 
commutativity of the right diagram in 
(\ref{eq:morphismanchoredrepns}) holds whenever $T$ is unitary
and the left diagram in (\ref{eq:morphismanchoredrepns}) commutes. 
To see this, the left diagram says $WT=LV.$ 
Taking the adjoint of this condition gives $T^*W^*=V^*L^*.$
Applying $T$ on the left and $L$ on the right gives
$W^*L=TV^*,$ which is the diagram on the right in  
(\ref{eq:morphismanchoredrepns}). 
\er

When the bounded linear maps in these definitions are isometries, 
the following lemma shows the isometry $T:\mK\to\mL$ 
is redundant and can be constructed from $L:\Hi\to\mI$. 

\blem
\label{lem:alternatedefnofmorphismofanchrep}
Let $\mK,\Hi,\mL,$ and $\mI$ be Hilbert spaces and let $V:\mK\to\Hi,$ 
$W:\mL\to\mI,$ and $L:\Hi\to\mI$ be isometries. 
Then $L$ satisfies 
\be
\label{eq:anchoredrepmorphismreduce}
L\big(V(\mK)\big)\subseteq W(\mL)
\quad\text{ and }\quad
L\big(V(\mK)^{\perp}\big)\subseteq W(\mL)^{\perp}
\ee
if and only if
there exists a (necessarily unique) isometry $T:\mK\to\mL$ such that  
the diagrams in (\ref{eq:morphismanchoredrepns}) commute. 
\elem

In (\ref{eq:anchoredrepmorphismreduce}), $W(\mL)^{\perp}$ stands for the orthogonal complement of $W(\mL)\subseteq\mathcal{I}$ and likewise for $V(\mK)^{\perp}\subseteq\Hi.$ 

\bprf[Proof of Lemma~\ref{lem:alternatedefnofmorphismofanchrep}]
First note that since $V$ and $W$ are isometries, their images are closed. Hence,  
\be
\label{eq:orthogonaldecomp}
\Hi=V(\mK)\oplus V(\mK)^{\perp}\quad\text{ and }\quad
\mI=W(\mL)\oplus W(\mL)^{\perp}.
\ee
\noindent
($\Rightarrow$)
Assume $L$ satisfies (\ref{eq:anchoredrepmorphismreduce}). 
Then the diagram
\be
\label{eq:WLV}
\xy0;/r.25pc/:
(45,0)*+{\mathcal{K}}="K";
(20,10)*+{V(\mK)\oplus V(\mK)^{\perp}}="H";
(20,-10)*+{V(\mK)}="VK";
(-20,10)*+{W(\mL)\oplus W(\mL)^{\perp}}="I";
(-20,-10)*+{W(\mL)}="WL";
(-45,0)*+{\mathcal{L}}="L";
{\ar"K";"H"_(0.35){V}};
{\ar"K";"VK"^(0.35){\pi_{V(\mK)} V}};
{\ar"H";"I"_{L}};
{\ar"VK";"WL"^{\pi_{W(\mL)} L i_{V(\mK)}}};
{\ar"I";"L"_(0.60){W^*}};
{\ar"WL";"L"^(0.60){(\pi_{W(\mL)} W)^*\;}};
{\ar"VK";"H"|-{i_{V(\mK)}}};
{\ar"I";"WL"|-{\pi_{W(\mL)}}};
\endxy
\ee
commutes by construction. In this diagram, $\pi_{W(\mL)}$ and $i_{V(\mK)}$ denote projection and inclusion maps, respectively. 
Note that $\pi_{W(\mL)} L i_{V(\mK)}$ is 
an isometry because $L\big(V(\mK)\big)\subseteq W(\mL)$ by (\ref{eq:orthogonaldecomp}). 
Setting $T:=W^* L V$, 
it follows from (\ref{eq:WLV}) that $T$ is the composite of the two unitary maps 
$\pi_{V(\mK)} V$ and ${(\pi_{W(\mL)} W)^*}$ and the isometry 
$\pi_{W(\mL)} L i_{V(\mK)}$. Therefore, $T$ is an isometry. 
Finally, the diagrams in (\ref{eq:morphismanchoredrepns}) commute because they are given by 
\be
\xy0;/r.25pc/:
(-15,30)*+{\mK}="K";
(15,30)*+{V(\mK)\oplus V(\mK)^{\perp}}="H";
(-15,10)*+{V(\mK)\oplus V(\mK)^{\perp}}="Hleft";
(-15,-30)*+{\mL}="L";
(15,-30)*+{W(\mL)\oplus W(\mL)^{\perp}}="I";
(-15,-10)*+{W(\mL)\oplus W(\mL)^{\perp}}="Ileft";
{\ar"K";"H"^(0.3){V}};
{\ar"L";"I"_(0.3){W}};
{\ar"Ileft";"L"_{W^*}};
{\ar"K";"Hleft"_{V}};
{\ar"Hleft";"Ileft"_{L}};
{\ar"H";"I"^{L}};
{\ar"Ileft";"I"^{P_{W(\mL)}}};
{\ar@{=}(-1,-4);(12,17)};
(3,7.5)*{\rotatebox{58}{\tiny $L(V(\mK))\subseteq W(\mL)$}};
{\ar@{=}(-2,-20);(-12,-26)};
(-6,-25)*{\rotatebox{33}{\tiny $W^*W=\id_{\mL}$}};
\endxy
\quad\text{and}\quad
\xy0;/r.25pc/:
(15,30)*+{\mK}="K";
(-15,30)*+{V(\mK)\oplus V(\mK)^{\perp}}="H";
(15,10)*+{V(\mK)\oplus V(\mK)^{\perp}}="Hright";
(15,-30)*+{\mL}="L";
(-15,-30)*+{W(\mL)\oplus W(\mL)^{\perp}}="I";
(15,-10)*+{W(\mL)\oplus W(\mL)^{\perp}}="Iright";
{\ar"H";"K"^(0.7){V^*}};
{\ar"I";"L"_(0.7){W^*}};
{\ar"H";"I"_{L}};
{\ar"Iright";"L"^{W^*}};
{\ar"K";"Hright"^{V}};
{\ar"Hright";"Iright"^{L}};
{\ar"H";"Hright"_{P_{V(\mK)}}};
{\ar@{=}(-12,-17);(1,4)};
(-3,-7.5)*{\rotatebox{58}{\tiny $L(V(\mK)^{\perp})\subseteq W(\mL)^{\perp}$}};
{\ar@{=}(2,20);(12,26)};
(6,25)*{\rotatebox{33}{\tiny $V^*V=\id_{\mK}$}};
\endxy
,
\ee
respectively, where $P_{W(\mL)}$ is the projection onto $W(\mL)$ inside $\mI$, and similarly for $P_{V(\mK)}.$ 

\noindent
($\Leftarrow$)
Assume an isometry $T$ satisfying (\ref{eq:morphismanchoredrepns})
exists. Commutativity of the diagram on the left requires $L\big(V(\mK)\big)\subseteq W(\mL)$, while commutativity of the diagram on the right requires $L\big(V(\mK)^{\perp}\big)\subseteq W(\mL)^{\perp}.$
\eprf

\begin{notation}
\label{not:morphismanchoredrepns}
A morphism of preserving anchored representations will be denoted
by the pair $(T,L):(\mK,\Hi,\pi,V)\to(\mL,\mI,\rho,W)$, 
even though the isometry $T$ is uniquely determined by $L$ as illustrated
in Lemma~\ref{lem:alternatedefnofmorphismofanchrep}. This is because $T$ need not be uniquely determined by $L$ in the general case of anchored representations and because we will use $T$ to relate anchored representations to OCP maps in this work. 
\end{notation}

\bx
\label{ex:pRep}
In the special case $\mK=\C$, 
an anchored representation of $\mA$ consists of 
a $^*$-representation $\pi:\mA\to\mB(\Hi)$ and 
a linear map
$V:\C\to\Hi.$
Such a map is uniquely characterized by the vector 
$\W:=V(1)\in\Hi,$ which is a unit vector if and only if $V$ is an isometry.
Hence, a preserving anchored representation
of the form $(\C,\Hi,\pi,V)$  
is equivalent to a pointed representation as introduced in 
Definition~5.1 in \cite{PaGNS}. 
If $(\C,\mI,\rho,W)$ is another preserving anchored representation, 
then a morphism $(\C,\Hi,\pi,V)\xrightarrow{(T,L)}(\C,\mI,\rho,W)$ of preserving anchored representations
consists of isometries $T:\C\to\C$ and $L:\Hi\to\mI$
satisfying (\ref{eq:intertwiner}) and (\ref{eq:morphismanchoredrepns}).  
$T$ must be of the form $T(\l)=\l e^{i\q}$ for all $\l\in\C$
for some $\q\in[0,2\pi).$ 
$L$ is an intertwiner of representations by (\ref{eq:intertwiner}). 
Let $\Xi:=W(1).$ The left diagram in 
(\ref{eq:morphismanchoredrepns}) entails 
$L(\W)=e^{i\q}\Xi.$ 
Since $V^*=\<\W,\;\cdot\;\>$ and $W^{*}=\<\Xi,\;\cdot\;\>,$ 
the right diagram in (\ref{eq:morphismanchoredrepns}) entails
$e^{i\q}\<\W,\;\cdot\;\>=\<\Xi,L(\;\cdot\;)\>$ as linear functionals on $\Hi.$ 
However, this second condition is implied by the first one in 
(\ref{eq:morphismanchoredrepns}) by Remark~\ref{rmk:unitarymorphismofpanrep}. Hence, this reproduces the morphisms of pointed representations in \cite{PaGNS} up to a phase. 
\ex

\blem
\label{lem:anrepA}
Let $(\mK,\Hi,\pi,V)\xrightarrow{(T,L)}(\mL,\mI,\rho,W)\xrightarrow{(S,M)}(\mM,\mJ,\s,X)$ 
be a pair of
composable morphisms of anchored representations on a $C^*$-algebra $\mA$. The composite 
of said morphisms is defined as $(S,M)\circ (T,L):=(ST,ML)$ 
and is a morphism of anchored representations. 
A similar statement is true for preserving anchored representations and their morphisms. 
The identity $(\mK,\Hi,\pi,V)\xrightarrow{(\id_{\mK},\id_{\Hi})}(\mK,\Hi,\pi,V)$
is the identity linear transformation on each Hilbert space.
The collection of all anchored representations on $\mA$ and their morphisms forms a category, which is 
denoted by $\mathbf{AnRep}(\mA).$ Similarly, the collection of all preserving anchored representations on $\mA$ and their morphisms forms a category, which is
denoted by $\mathbf{PAnRep}(\mA).$
Furthermore, $\mathbf{PAnRep}(\mA)$ is a subcategory of $\mathbf{AnRep}(\mA)$.
\elem

\blem
\label{lem:anrepf}
Let $f:\mA'\to\mA$ be a $^*$-homomorphism of $C^*$-algebras. 
The assignment
\be
\begin{split}
\mathbf{AnRep}(\mA)&\xrightarrow{\mathbf{AnRep}_{f}}\mathbf{AnRep}(\mA')\\
(\mK,\Hi,\pi,V)&\xmapsto{\qquad\quad}(\mK,\Hi,\pi\circ f,V)\\
\Big(\!(\mK,\Hi,\pi,V)\!\xrightarrow{\!(T,L)\!}\!(\mL,\mI,\rho,W)\!\Big)&\xmapsto{\quad\qquad}\Big(\!(\mK,\Hi,\pi\circ f,V)\!\xrightarrow{\!(T,L)\!}\!(\mL,\mI,\rho\circ f,W)\!\Big)
\end{split}
\ee
defines a functor.  Furthermore, it restricts to a functor
$\mathbf{PAnRep}(\mA)\xrightarrow{\mathbf{PAnRep}_{f}}\mathbf{PAnRep}(\mA')$.
\elem

\blem
\label{lem:anrepfunctor}
The assignment 
\be
\begin{split}
\CAlg^{\op}&\xrightarrow{\mathbf{\mathbf{AnRep}}}\tCAT\\
\mA&\xmapsto{\;\qquad\;\;\!}\mathbf{AnRep}(\mA)\\
\left(\mA'\xrightarrow{f}\mA\right)&\xmapsto{\;\qquad\;\;\!}\left(\mathbf{AnRep}(\mA)\xrightarrow{\mathbf{AnRep}_{f}}\mathbf{AnRep}(\mA')\right)
\end{split}
\ee
defines a functor. The same is true for $\CAlg^{\op}\xrightarrow{\mathbf{\mathbf{PAnRep}}}\tCAT$. 
\elem

\section{The restriction natural transformation}
\label{sec:rest}

The follow proposition illustrates how to construct an OCP map from an anchored representation.

\bn
\label{prop:rest}
Let $\mA$ be a $C^*$-algebra. The assignment
\be
\begin{split}
\mathbf{AnRep}(\mA)&\xrightarrow{\rest_{\mA}}\mathbf{OCP}(\mA)\\
(\mK,\Hi,\pi,V)&\xmapsto{\qquad\!}(\mK,\mathrm{Ad}_{V^*}\circ\pi)\\
\Big((\mK,\Hi,\pi,V)\xrightarrow{(T,L)}(\mL,\mI,\rho,W)\Big)&\xmapsto{\qquad\!}\Big((\mK,\mathrm{Ad}_{V^*}\circ\pi)\xrightarrow{T}(\mL,\mathrm{Ad}_{W^*}\circ\rho)\Big)
\end{split}
\ee
defines a functor. Furthermore, it restricts to a functor $\mathbf{PAnRep}(\mA)\xrightarrow{\rest_{\mA}}\mathbf{OpSt}(\mA).$ 
\en

Here, $(\mathrm{Ad}_{V^*}\circ\pi)(a):=V^*\pi(a)V$ for all $a\in\mA$ (cf.~Example~\ref{ex:cpmaps}). 

\bprf[Proof of Proposition~\ref{prop:rest}]
Let $(\mK,\Hi,\pi,V)$ be an anchored representation of $\mA.$ 
Then $\mathrm{Ad}_{V^*}\circ\pi$ is an OCP map because it is the composite of the CP map $\mathrm{Ad}_{V^*}$ 
and the $^*$-homomorphism $\pi$. 
Since $\mathrm{Ad}_{V^*}$ is unital if and only if $V$ is an isometry, 
$\mathrm{Ad}_{V^*}\circ\pi$ is an operator state when $(\mK,\Hi,\pi,V)$ is a preserving anchored representation.
Let $(\mK,\Hi,\pi,V)\xrightarrow{(T,L)}(\mL,\mI,\rho,W)$ be a morphism of anchored representations. 
In order for $(\mK,\mathrm{Ad}_{V^*}\circ\pi)\xrightarrow{T}(\mL,\mathrm{Ad}_{W^*}\circ\rho)$ 
to be a morphism of OCP maps, the diagram 
\be
\xy0;/r.25pc/:
(-12.5,7.5)*+{\mK}="K1";
(12.5,7.5)*+{\mK}="K2";
(-12.5,-7.5)*+{\mL}="L1";
(12.5,-7.5)*+{\mL}="L2";
{\ar"K1";"K2"^{\mathrm{Ad}_{V^*}(\pi(a))}};
{\ar"L1";"L2"_{\mathrm{Ad}_{W^*}(\rho(a))}};
{\ar"K1";"L1"_{T}};
{\ar"K2";"L2"^{T}};
\endxy
\ee
must commute for all $a\in\mA.$ Expanding out the definition of the adjoint action map 
provides the diagram 
\be
\xy0;/r.25pc/:
(-37.5,7.5)*+{\mK}="K1";
(-12.5,7.5)*+{\Hi}="H1";
(12.5,7.5)*+{\Hi}="H2";
(37.5,7.5)*+{\mK}="K2";
(-37.5,-7.5)*+{\mL}="L1";
(-12.5,-7.5)*+{\mI}="I1";
(12.5,-7.5)*+{\mI}="I2";
(37.5,-7.5)*+{\mL}="L2";
{\ar"K1";"H1"^{V}};
{\ar"H1";"H2"^{\pi(a)}};
{\ar"H2";"K2"^{V^*}};
{\ar"L1";"I1"_{W}};
{\ar"I1";"I2"_{\rho(a)}};
{\ar"I2";"L2"_{W^*}};
{\ar"K1";"L1"_{T}};
{\ar"K2";"L2"^{T}};
\endxy
.
\ee
This diagram commutes because the diagram
\be
\xy0;/r.25pc/:
(-37.5,7.5)*+{\mK}="K1";
(-12.5,7.5)*+{\Hi}="H1";
(12.5,7.5)*+{\Hi}="H2";
(37.5,7.5)*+{\mK}="K2";
(-37.5,-7.5)*+{\mL}="L1";
(-12.5,-7.5)*+{\mI}="I1";
(12.5,-7.5)*+{\mI}="I2";
(37.5,-7.5)*+{\mL}="L2";
{\ar"K1";"H1"^{V}};
{\ar"H1";"H2"^{\pi(a)}};
{\ar"H2";"K2"^{V^*}};
{\ar"L1";"I1"_{W}};
{\ar"I1";"I2"_{\rho(a)}};
{\ar"I2";"L2"_{W^*}};
{\ar"K1";"L1"_{T}};
{\ar"H1";"I1"_{L}};
{\ar"H2";"I2"^{L}};
{\ar"K2";"L2"^{T}};
\endxy
\ee
commutes due to Definition~\ref{defn:anchoredrepn}.%
\footnote{This is where commutativity of 
(\ref{eq:opstatesdiagram}) from Definition~\ref{defn:opstate} is used 
(cf.~Remark~\ref{rmk:morphismopstates}).
When working with morphisms of operator states, (\ref{eq:opstatemorphismother}) would be too strong, and its 
commutativity would not follow from Definition~\ref{defn:anchoredrepn}. 
}
\eprf

\blem
\label{lem:restnatural}
Let $f:\mA'\to\mA$ be a $^*$-homomorphism of $C^*$-algebras.
Then the diagram 
\be
\xy0;/r.25pc/:
(-20,7.5)*+{\mathbf{AnRep}(\mA)}="TL";
(20,7.5)*+{\mathbf{AnRep}(\mA')}="TR";
(-20,-7.5)*+{\mathbf{OCP}(\mA)}="BL";
(20,-7.5)*+{\mathbf{OCP}(\mA')}="BR";
{\ar"TL";"TR"^{\mathbf{AnRep}_{f}}};
{\ar"BL";"BR"_{\mathbf{OCP}_{f}}};
{\ar"TL";"BL"_{\rest_{\mA}}};
{\ar"TR";"BR"^{\rest_{\mA'}}};
\endxy
\ee
of functors commutes (on the nose). A similar statement holds for the subcategories obtained from $\mathbf{PAnRep}$ and $\mathbf{OpSt}$.
\elem

Lemma~\ref{lem:restnatural} states that $\rest$ is a natural transformation
\be
\xy0;/r.25pc/:
(-20,0)*+{\CAlg^{\op}}="1";
(20,0)*+{\tCAT}="2";
{\ar@/^1.5pc/"1";"2"^{\mathbf{OCP}}};
{\ar@/_1.5pc/"1";"2"_{\mathbf{AnRep}}};
{\ar@{=>}(0,-5);(0,5)_{\rest}};
\endxy
,
\ee
a special kind of oplax-natural transformation (cf.~Definition~\ref{defn:semipseudonaturaltransformation}).

\section{Stinespring's oplax-natural transformation}
\label{sec:Stinespring}

In the construction of a left adjoint to $\rest,$ some
preliminary facts will be needed.  

\blem
\label{lem:CPtoP}
Let $\vf:\mA\stoch\mB(\mK)$ be a completely positive map. 
Let $\vec{v}:=(v_{1},\dots,v_{n})\in\mK\oplus\cdots\oplus\mK$
denote a vector in the direct sum of $\mK$ with itself $n$ times. 
Then the assignment
\be
\begin{split}
s_{\vf,\vec{v}}:\mM_{n}(\mA)&\stoch\C\\
A&\xmapsto{\;\quad}\sum_{i,j=1}^{n}\<v_{i},\vf(a_{ij})v_{j}\>_{\mK}
\end{split}
\ee
is a positive linear functional. 
\elem

\bprf
Suppose $A\in\mM_{n}(\mA)$ is positive. 
Then, because $\vf$ is completely positive, 
$\vf_{n}(A)\ge0.$ Hence, 
\be
s_{\vf,\vec{v}}(A)
=\sum_{i,j=1}^{n}\<v_{i},\vf(a_{ij})v_{j}\>_{\mK}
=\<\vec{v},\vf_{n}(A)\vec{v}\>_{\mK\oplus\cdots\oplus\mK}\ge0.
\ee
Linearity of $s_{\vf,\vec{v}}$ follows from linearity of $\vf_{n}$ 
and linearity of the inner product in the right variable. 
\eprf

\blem
\label{lem:XmodN}
Let $X$ and $Y$ be topological vector spaces, let $X\xrightarrow{f}Y$ be a continuous linear map, and let $N\subseteq X$ be a closed vector subspace of $X.$ 
\begin{enumerate}[(a)]
\item
If $f(x)=0$ for all $x\in N,$ then there exists a unique continuous linear map $X/N\xrightarrow{g}Y$ such that 
\be
\xy0;/r.25pc/:
(-10,7.5)*+{X}="X";
(-10,-7.5)*+{X/N}="XN";
(10,7.5)*+{Y}="Y";
{\ar"X";"Y"^{f}};
{\ar@{->>}"X";"XN"};
{\ar"XN";"Y"_{g}};
\endxy
\ee
commutes. Here, $X/N$ is the quotient space of $X$ modulo $N$ and $X\twoheadrightarrow X/N$ is the quotient map.  
\item
If $M\subseteq Y$ is a closed subspace of $Y$ and $f(N)\subseteq M$, then there exists a unique continuous linear map $X/N\xrightarrow{h}Y/M$ such that 
\be
\xy0;/r.25pc/:
(-10,7.5)*+{X}="X";
(-10,-7.5)*+{X/N}="XN";
(10,7.5)*+{Y}="Y";
(10,-7.5)*+{Y/M}="YM";
{\ar"X";"Y"^{f}};
{\ar@{->>}"X";"XN"};
{\ar@{->>}"Y";"YM"};
{\ar"XN";"YM"_{h}};
\endxy
\ee
commutes.
\end{enumerate}
\elem

\bprf
For the first fact, see Theorem~1.41 and Exercise~9 in Chapter~1 of Rudin~\cite{Ru91}. 
The second fact is a consequence of the first. 
\eprf

\blem
\label{lem:bdd}
Let $M:\Hi\to\Hi$ a bounded positive 
operator on a Hilbert space $\Hi$. Then $M\le\lVert M\rVert\id_{\Hi}$ (in the sense that $\lVert M\rVert\id_{\Hi}-M$ is a positive 
operator). 
\elem

\bprf
This follows from the inequality $\<x,Mx\>_{\Hi}\le\lVert x\rVert_{\Hi}\lVert Mx\rVert_{\Hi}\le\lVert M\rVert\lVert x\rVert_{\Hi}^{2}$ for all $x\in\Hi$ by positivity of $M$, Cauchy--Schwarz, and the definition of the norm on $\mB(\Hi).$ 
\eprf

The following is the main result of this work. 

\bt
\label{thm:stinespring}
There exists a left adjoint $\mathbf{Stine}:\mathbf{OCP}\Rightarrow\mathbf{AnRep}$
to the natural transformation $\rest:\mathbf{AnRep}\Rightarrow\mathbf{OCP}$
\be
\xy0;/r.25pc/:
(-25,0)*+{\CAlg^{\op}}="1";
(25,0)*+{\tCAT}="2";
(0,0)*{\dashv};
{\ar@/^1.75pc/"1";"2"^{\mathbf{OCP}}};
{\ar@/_1.75pc/"1";"2"_{\mathbf{AnRep}}};
{\ar@{=>}(2.5,-5.5);(2.5,5.5)_{\rest}};
{\ar@{=>}(-2.5,5.5);(-2.5,-5.5)_{\mathbf{Stine}}};
\endxy
\ee
in the 2-category $\mathbf{Fun}(\CAlg^{\op},\tCAT).$ 
Furthermore, $\mathbf{Stine}:\mathbf{OCP}\Rightarrow\mathbf{AnRep}$ restricts to a left adjoint $\mathbf{Stine}:\mathbf{OpSt}\Rightarrow\mathbf{PAnRep}$ to $\rest:\mathbf{PAnRep}\Rightarrow\mathbf{OpSt}$.
\et

The 2-category $\mathbf{Fun}(\CAlg^{\op},\tCAT)$ is defined in Notation~\ref{not:functor2cat} and the definition of a 2-categorical adjunction is reviewed in Definition~\ref{defn:adjunction}. The consequences, universal properties, and applications of this theorem are discussed in Sections~\ref{sec:generalizations}~and~\ref{sec:comparetoGNS}. 

\bprf[Proof of Theorem~\ref{thm:stinespring}]
The proof will be split up into several steps. 
\begin{enumerate}[i.]
\item
\label{i}
For a fixed $C^*$-algebra $\mA,$ define the functor 
$\mathbf{Stine}_{\mA}:\mathbf{OCP}(\mA)\to\mathbf{AnRep}(\mA)$
on objects. 
\item
\label{ii}
For a fixed $C^*$-algebra $\mA,$ define the functor 
$\mathbf{Stine}_{\mA}:\mathbf{OCP}(\mA)\to\mathbf{AnRep}(\mA)$
on morphisms and prove functoriality. 
\item
\label{iii}
For a fixed $^*$-homomorphism $f:\mA'\to\mA,$ define the natural transformation 
$\mathbf{Stine}_{f}:\mathbf{Stine}_{\mA'}\circ\mathbf{OCP}_{f}\Rightarrow\mathbf{AnRep}_{f}\circ\mathbf{Stine}_{\mA}.$
\item
\label{iv}
Prove that $\mathbf{Stine}$ is an oplax-natural transformation 
(cf.~Definition~\ref{defn:semipseudonaturaltransformation}). 
\item
\label{v}
For a fixed $C^*$-algebra $\mA,$ 
construct the appropriate natural transformation 
$m_{\mA}:\mathbf{Stine}_{\mA}\circ\rest_{\mA}\Rightarrow\id_{\mathbf{AnRep}(\mA)}.$
\item
\label{vi}
Show that%
\footnote{The vertical concatenation is the vertical composition of
natural transformations. This notation is used in \cite{PaGNS} and \cite{Pa18}. It is also reviewed in Appendix~\ref{appendix} in the present work.
In particular, applying $m$ to a $C^*$-algebra $\mA$ gives a natural transformation 
$m_{\mA}:\mathbf{Stine}_{\mA}\circ\rest_{\mA}\Rightarrow\id_{\mathbf{AnRep}(\mA)}.$
}
$m:\begin{matrix}[0.7]\rest\\\mathbf{Stine}\end{matrix}\Rrightarrow\id_{\mathbf{AnRep}}$ is a modification 
(cf.~Definition~\ref{defn:modification})
between oplax-natural transformations in 
the 2-category $\mathbf{Fun}(\CAlg^{\op},\tCAT).$ 
\item
\label{vii}
Show that $\begin{matrix}[0.7]\mathbf{Stine}\\\rest\end{matrix}=\id_{\mathbf{OCP}}.$
\item
\label{viii}
Prove the zig-zag identities for adjunctions in 2-categories (cf.~Definition~\ref{defn:adjunction}), 
i.e.~complete the proof that $(\mathbf{Stine},\rest,\id,m)$
is an adjunction in $\mathbf{Fun}(\CAlg^{\op},\tCAT).$
\end{enumerate}
In all of the above steps, justifications for reducing to operator states and preserving anchored representations will be provided. 
In what follows, if a proof for any claim is not supplied, 
it is because the justification is analogous to the 
standard GNS construction arguments or it
follows easily from the definitions. The reader is 
referred to \cite{PaGNS} for more details. 
\begin{enumerate}[i.]
\item
\label{item:1}
The construction of an anchored representation from an OCP map will 
be Stinespring's construction (cf.~the proof of sufficiency of Theorem~1 in Stinespring~\cite{St55}). 
Let $\mA$ be a $C^*$-algebra and let $(\mK,\vf)$ be an OCP map on $\mA.$ 
Recall, this means $\vf:\mA\stoch\mB(\mK)$ is a CP map. 
Let $\mA\otimes\mK$ denote the vector space tensor product of $\mA$ with $\mK.$ 
In particular, elements of $\mA\otimes\mK$ are finite sums of tensor products of vectors
in $\mA$ and vectors in $\mK$
(in fact, all sums that follow are finite). 
The function
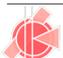
\be
\label{eq:sesquibilinearity}
\begin{split}
(\mA\times\mK)\times(\mA\times\mK)&\to\C\\
\big((a,v),(b,w)\big)&\xmapsto{\quad}\<v,\vf(a^{*}b)w\>_{\mK}
\end{split}
\ee
is conjugate bilinear in the first $\mA\times\mK$ factor and bilinear in the second $\mA\times\mK$ factor. Hence, by the universal property of the algebraic tensor product (cf.~Chapter~IV~Section~5 in Hungerford~\cite{Hu74}), the assignment
\be
\label{eq:sesquilinear}
\begin{split}
(\mA\otimes\mK)\times(\mA\otimes\mK)&\xrightarrow{\llangle\;\cdot\;,\;\cdot\;\rrangle_{\vf}}\C\\
\left(\sum_{i}a_{i}\otimes v_{i},\sum_{j}b_{j}\otimes w_{j}\right)&\xmapsto{\;\;\;\quad\quad\;\;\!}\sum_{i,j}\<v_{i},\vf(a_{i}^{*}b_{j})w_{j}\>_{\mK}
\end{split}
\ee
is well-defined, conjugate linear in the first variable, 
and linear in the second variable.
Furthermore, $\llangle\;\cdot\;,\;\cdot\;\rrangle_{\vf}$ 
satisfies%
\footnote{An overline here indicates complex conjugation.
This is not to be confused with the closure such as in (\ref{eq:stinespringspace}).}
\be
\label{eq:swapconjugate}
\llangle\z,\xi\rrangle_{\vf}=\ov{\llangle\xi,\z\rrangle_{\vf}}
\qquad\forall\;\z,\xi\in\mA\otimes\mK.
\ee
Since the matrix
\be
A:=
\begin{bmatrix}
a_{1}^{*}a_{1}&\cdots&a_{1}^{*}a_{n}\\
\vdots&&\vdots\\
a_{n}^{*}a_{1}&\cdots&a_{n}^{*}a_{n}\\
\end{bmatrix}
=
\begin{bmatrix}
a_{1}&\cdots&a_{n}\\
0&\cdots&0\\
\vdots&&\vdots\\
0&\cdots&0
\end{bmatrix}^*
\begin{bmatrix}
a_{1}&\cdots&a_{n}\\
0&\cdots&0\\
\vdots&&\vdots\\
0&\cdots&0
\end{bmatrix}
\ee
in $\mM_{n}(\mA)$ is positive for all $a_{1},\dots,a_{n}\in\mA,$ 
\be
\label{eq:nonnegativeondouble}
\llangle\xi,\xi\rrangle_{\vf}
=\sum_{i,j}\<v_{i},\vf(a_{i}^{*}a_{j})v_{j}\>_{\mK}
\ge0
\ee
by Lemma~\ref{lem:CPtoP} applied to $\xi=\sum_{i=1}^{n}a_{i}\otimes v_{i}\in\mA\otimes\mK.$
 
By the properties of $\llangle\;\cdot\;,\;\cdot\;\rrangle_{\vf},$
it follows that 
\be
\label{eq:CauchySchwarz}
|\llangle\xi,\zeta\rrangle_{\vf}|^{2}\le\llangle\xi,\xi\rrangle_{\vf}\llangle\zeta,\zeta\rrangle_{\vf}
\qquad\forall\;\xi,\z\in\mA\otimes\mK. 
\ee
This is the Cauchy--Schwarz inequality for such sesquilinear forms (cf.~Construction~3.1 in~\cite{PaGNS}).
Thus, $\llangle\;\cdot\;,\;\cdot\;\rrangle_{\vf}$ is a sesquilinear form whose associated seminorm endows $\mA\otimes\mK$ with the structure of a topological vector space. 
In general, $\llangle\;\cdot\;,\;\cdot\;\rrangle_{\vf}$ is not positive semi-definite. Hence, set 
\be
\label{eq:StinespringIdeal}
\mN_{\vf}:=\big\{\z\in\mA\otimes\mK\;:\;\llangle\z,\z\rrangle_{\vf}=0\big\}
\ee
to be its null-space. 
From the Cauchy--Schwarz inequality~(\ref{eq:CauchySchwarz}), 
it follows that 
\be
\label{eq:leftidealproperty}
\llangle\xi,\z\rrangle_{\vf}=0\quad\forall\;\xi\in\mA\otimes\mK,\;\z\in\mN_{\vf}.
\ee
Using this, one can show that 
$\mN_{\vf}$ is a closed vector subspace of $\mA\otimes\mK$ (it is closed since it is defined as the inverse image of $\{0\}$ under a continuous map). 

By the universal property of the tensor product, for each $a\in\mA,$ the map
\be
\label{eq:stinespringreppi}
\begin{split}
\mA\otimes\mK&\xrightarrow{\pi'_{\vf}(a)}\mA\otimes\mK\\
\sum_{i}a_{i}\otimes v_{i}&\xmapsto{\qquad}\sum_{i}aa_{i}\otimes v_{i}
\end{split}
\ee
is a well-defined linear transformation. 
If we write $\End(V)$ for the algebra of linear transformations
from a vector space $V$ to itself, then 
$\pi'_{\vf}:\mA\to\End(\mA\otimes\mK)$ defines a
representation of the algebra $\mA$ on the vector space 
$\mA\otimes\mK.$ 
Furthermore, for each $a\in\mA$ and $\xi=\sum_{i=1}^{n}a_{i}\otimes v_{i}\in\mA\otimes\mK,$ $\pi'_{\vf}(a)$ satisfies 
\be
\label{eq:Stinespringbounded}
\begin{split}
\llangle\pi'_{\vf}(a)\xi,\pi'_{\vf}(a)\xi\rrangle_{\vf}
&=\sum_{i,j=1}^{n}\big\<v_{i},\vf(a_{i}^*a^*aa_{j})v_{j}\big\>_{\mK}\;\;\;\text{by definition of $\pi'_{\vf}$ and $\llangle\;\cdot\;,\;\cdot\;\rrangle_{\vf}$}\\
&=s_{\vf,\vec{v}}
\left(
\begin{bmatrix}
a_{1}&\cdots&a_{n}\\
0&\cdots&0\\
\vdots&&\vdots\\
0&\cdots&0
\end{bmatrix}^*
\begin{bmatrix}
a^*a&&0\\
&\ddots&\\
0&&a^*a
\end{bmatrix}
\begin{bmatrix}
a_{1}&\cdots&a_{n}\\
0&\cdots&0\\
\vdots&&\vdots\\
0&\cdots&0
\end{bmatrix}
\right)\\
&\le
\left\lVert 
\begin{bmatrix}
a^*a&&0\\
&\ddots&\\
0&&a^*a
\end{bmatrix}
\right\rVert
s_{\vf,\vec{v}}
\left(
\begin{bmatrix}
a_{1}&\cdots&a_{n}\\
0&\cdots&0\\
\vdots&&\vdots\\
0&\cdots&0
\end{bmatrix}^*
\begin{bmatrix}
a_{1}&\cdots&a_{n}\\
0&\cdots&0\\
\vdots&&\vdots\\
0&\cdots&0
\end{bmatrix}
\right)\\
&=\lVert a^*a\rVert
s_{\vf,\vec{v}}
\left(
\begin{bmatrix}
a_{1}^*a_{1}&\cdots&a_{1}^*a_{n}\\
\vdots&&\vdots\\
a_{n}^*a_{1}&\cdots&a_{n}^*a_{n}\\
\end{bmatrix}
\right)
\\
&=\lVert a\rVert^{2}
\llangle\xi,\xi\rrangle_{\vf}.
\end{split}
\ee
In this calculation, $\vec{v}:=(v_{1},\dots,v_{n})\in\mK\oplus\cdots\oplus\mK$ and
the norm $\lVert a\rVert$ of $a$ is the one from the $C^*$-algebra $\mA$. 
The third line in (\ref{eq:Stinespringbounded}) follows from Lemma~\ref{lem:CPtoP} and the inequality $|\w(y^*xy)|\le\lVert x\rVert\w(y^*y)$ for all $x,y$ in a $C^*$-algebra and $\w$ a positive
linear functional on that $C^*$-algebra (see Proposition~2.1.5. part~(ii) in 
Dixmier~\cite{Di77} for a proof of this inequality). In this case, this inequality is
applied to the positive linear functional 
$\w:=s_{\vf,\vec{v}}:\mM_{n}(\mA)\stoch\C$ with $x$ positive
so that $|\w(y^*xy)|=\w(y^*xy).$
The fourth line in (\ref{eq:Stinespringbounded}) holds because 
the norm of $\mathrm{diag}(a^*a,\dots,a^*a)$ in $\mM_{n}(\mA)$ is
equal to $\lVert a^*a\rVert.$ 
This is because every injective $^{*}$-homomorphism of $C^*$-algebras is an isometry (cf.~Propositions~1.3.7~and~1.8.1 in Dixmier~\cite{Di77}). In this case, the $^{*}$-homomorphism is given by the function
$\mathcal{A}\to\mathcal{M}_{n}(\mA)$ sending $a\in\mA$ to $\mathrm{diag}(a,\dots,a).$
The last line of (\ref{eq:Stinespringbounded}) follows from the $C^*$-identity for $C^*$-algebras
and the definitions of  $\llangle\xi,\xi\rrangle_{\vf}$ and $s_{\vf,\vec{v}}.$ 
Thus, (\ref{eq:Stinespringbounded}) 
shows that $\pi'_{\vf}(a)$ is bounded/continuous. If we write $\mathcal{B}(V)$ for the algebra of bounded operators on a seminormed vector space $V$, then $\pi'_{\vf}(a)\in\mB(\mA\otimes\mK)$ for all $a\in\mA.$ 

Furthermore, (\ref{eq:Stinespringbounded}) 
shows that $\mN_{\vf}$ is an invariant subspace under the
$\pi'_{\vf}$ action, meaning $\pi'_{\vf}(a)\z\in\mathcal{N}_{\vf}$ for all $\z\in\mathcal{N}_{\vf}$ and $a\in\mA.$ Therefore, the quotient space
$(\mA\otimes\mK)/\mN_{\vf}$ has a well-defined action
$\pi_{\vf}:\mA\to\mB\big((\mA\otimes\mK)/\mN_{\vf}\big)$ by Lemma~\ref{lem:XmodN}.
If an element of $(\mA\otimes\mK)/\mN_{\vf}$ is denoted by $[\xi]_{\vf}$, or just $[\xi]$ when working with a fixed OCP map $\vf$, 
this induced action of $a\in\mA$ on $[\xi]$ is given by 
\be
\label{eq:stinespringrepaction}
\pi_{\vf}(a)[\xi]:=[\pi'_{\vf}(a)\xi].
\ee
By (\ref{eq:leftidealproperty}), 
the sesquilinear form $\llangle\;\cdot\;,\;\cdot\;\rrangle_{\vf}$
descends to a well-defined inner product
\be
\label{eq:inducedinnerproduct}
\begin{split}
(\mA\otimes\mK)/\mN_{\vf}\times(\mA\otimes\mK)/\mN_{\vf}&\xrightarrow{\<\;\cdot\;,\;\cdot\;\>_{\vf}}\C\\
([\xi],[\z])&\xmapsto{\;\qquad\;\;}\llangle\xi,\z\rrangle_{\vf}.
\end{split}
\ee
The fact that $\<\;\cdot\;,\;\cdot\;\>_{\vf}$
is positive definite follows from 
(\ref{eq:nonnegativeondouble}) and the definition of
$\mN_{\vf}$ in (\ref{eq:StinespringIdeal}).
Let
\be
\label{eq:stinespringspace}
\Hi_{\vf}:=\ov{(\mA\otimes\mK)/\mN_{\vf}}
\ee
denote the completion of this topological vector space
with respect to the inner product $\<\;\cdot\;,\;\cdot\;\>_{\vf}.$
Since $\pi_{\vf}(a)$ is a bounded linear operator 
on $(\mA\otimes\mK)/\mN_{\vf},$ 
it extends uniquely to a bounded linear operator, also denoted by $\pi_{\vf}(a),$ on $\Hi_{\vf}.$ 
Furthermore, since $\pi_{\vf}(a)\in\mB(\Hi_{\vf}),$ it has an adjoint $\pi_{\vf}(a)^*.$ 
This adjoint satisfies 
\be
\begin{split}
&\big\<\pi_{\vf}(a)^*[\xi],[\z]\big\>_{\vf}
=\big\<[\xi],\pi_{\vf}(a)[\z]\big\>_{\vf}
=\sum_{i,j}\big\<v_{i},\vf(a_{i}^*ab_{j})w_{j}\big\>_{\mK}\\
&\quad=\sum_{i,j}\Big\<v_{i},\vf\big((a^*a_{i})^*b_{j}\big)w_{j}\Big\>_{\mK}
=\sum_{i,j}\llangle a^*a_{i}\otimes v_{i},b_{j}\otimes w_{j}\rrangle_{\vf}
=\big\<\pi_{\vf}(a^*)[\xi],[\z]\big\>_{\vf}
\end{split}
\ee
for all $\xi=\big[\sum_{i}a_{i}\otimes v_{i}\big],\z=\big[\sum_{j}b_{j}\otimes w_{j}\big]\in(\mA\otimes\mK)/\mathcal{N}_{\vf}.$ 
By uniqueness of bounded adjoints, $\pi_{\vf}(a^*)=\pi_{\vf}(a)^*.$ 
Finally, note that $\pi_{\vf}(1_{\mA})=\id_{\Hi_{\vf}}.$ 
Thus, $\pi_{\vf}:\mA\to\mB(\Hi_{\vf})$ defines a $^*$-homomorphism.

Now, set 
\be
\label{eq:stinespringisometry}
\begin{split}
\mK&\xrightarrow{V_{\vf}}
\Hi_{\vf}\\
w&\xmapsto{\quad}[1_{\mA}\otimes w]. 
\end{split}
\ee
Then, 
\be
\label{eq:Vvfisometry}
\big\lVert V_{\vf}(w)\big\rVert_{\vf}^2
=\big\lVert[1_{\mA}\otimes w]\big\rVert^2_{\vf}
=\big\<w,\vf(1_{\mA}^*1_{\mA})w\big\>_{\mK}
=\big\<w,\vf(1_{\mA})w\big\>_{\mK}
\le
\big\lVert\vf(1_{\mA})\big\rVert\lVert w\rVert^2_{\mK}
\ee
for all $w\in\mK.$ 
The last inequality in (\ref{eq:Vvfisometry}) follows from Cauchy--Schwarz for $\<\;\cdot\;,\;\cdot\;\>_{\mK}$ and positivity of $\vf$. 
This proves that $V_{\vf}$ is bounded.
Note that if $\vf$ is unital, the inequality in (\ref{eq:Vvfisometry}) becomes an equality, which proves $V_{\vf}$ is an isometry.

This concludes Stinespring's construction 
of an anchored representation
from an OCP map, 
i.e.~the functor
\be
\label{eq:stineonobjects}
\begin{split}
\mathbf{OCP}(\mA)_{0}&\xrightarrow{\mathbf{Stine}_{\mA}}\mathbf{AnRep}(\mA)_{0}\\
(K,\vf)&\xmapsto{\;\;\qquad\;}(\mK,\Hi_{\vf},\pi_{\vf},V_{\vf})
\end{split}
\ee
on objects of $\mathbf{OCP}(\mA).$ It restricts to a well-defined map $\mathbf{OpSt}(\mA)_{0}\xrightarrow{\mathbf{Stine}_{\mA}}\mathbf{PAnRep}(\mA)_{0}$.

\item
Let $(\mK,\vf)\xrightarrow{T}(\mL,\psi)$ be a morphism of OCP maps.
By Definition~\ref{defn:opstate}, 
this means $T:\mK\to\mL$ is a bounded linear map and
$\psi(a) T=T\vf(a)$ for all $a\in\mA.$ 
Let $(\mK,\Hi_{\vf},\pi_{\vf},V_{\vf})$ and $(\mL,\Hi_{\psi},\pi_{\psi},V_{\psi})$
be the corresponding Stinespring anchored representations 
from the first step. 
Set
\be
\label{eq:stineT}
\mA\otimes\mK\xrightarrow{L'_{T}:=\id_{\mA}\otimes T}\mA\otimes\mL.
\ee
Then, for $\xi:=\sum_{i}a_{i}\otimes v_{i}\in\mA\otimes\mK,$ 
\be
\begin{split}
\llangle L'_{T}(\xi),L'_{T}(\xi)\rrangle_{\psi}
&=\sum_{i,j}\<T(v_{i}),\psi(a_{i}^*a_{j})T(v_{j})\>_{\mL}\quad\text{by (\ref{eq:stineT}) and (\ref{eq:sesquilinear})}\\
&=\sum_{i,j}\<T(v_{i}),T\vf(a_{i}^*a_{j})v_{j}\>_{\mL}\quad\text{by (\ref{eq:opstatesdiagram})}\\
&=\sum_{i,j}\<v_{i},T^*T\vf(a_{i}^*a_{j})v_{j}\>_{\mK}\\
&\le\lVert T^*T\rVert\llangle\xi,\xi\rrangle_{\vf}\quad\text{ by Lemma~\ref{lem:bdd}}.
\end{split}
\ee
Therefore, $L'_{T}(\mathcal{N}_{\vf})\subseteq\mathcal{N}_{\psi}$. Hence,  Lemma~\ref{lem:XmodN} implies $L'_{T}$ descends to a 
bounded 
linear map $L_{T}:{(\mA\otimes\mK)/\mN_{\vf}}\to{(\mA\otimes\mL)/\mN_{\psi}}$,
which itself extends uniquely to a bounded linear map $L_{T}:\Hi_{\vf}\to\Hi_{\psi}.$ 
Note that when $T$ is an isometry, $L_{T}$ is also an isometry. 
Commutativity of 
\be
\label{eq:stinemorphismTLT}
\xy0;/r.25pc/:
(-12.5,7.5)*+{\Hi_{\vf}}="H1";
(12.5,7.5)*+{\Hi_{\vf}}="H2";
(-12.5,-7.5)*+{\Hi_{\psi}}="I1";
(12.5,-7.5)*+{\Hi_{\psi}}="I2";
{\ar"H1";"H2"^{\pi_{\vf}(a)}};
{\ar"I1";"I2"_{\pi_{\psi}(a)}};
{\ar"H1";"I1"_{L_{T}}};
{\ar"H2";"I2"^{L_{T}}};
\endxy
\quad\forall\;a\in\mA
\quad\text{and}\quad
\xy0;/r.25pc/:
(-12.5,7.5)*+{\mK}="K";
(12.5,7.5)*+{\Hi_{\vf}}="H";
(-12.5,-7.5)*+{\mL}="L";
(12.5,-7.5)*+{\Hi_{\psi}}="I";
{\ar"K";"H"^{V_{\vf}}};
{\ar"L";"I"_{V_{\psi}}};
{\ar"K";"L"_{T}};
{\ar"H";"I"^{L_{T}}};
\endxy
\ee
follow directly from the definitions. However, 
commutativity of 
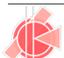
\be
\label{eq:stinemorphismdual}
\xy0;/r.25pc/:
(12.5,7.5)*+{\mK}="K";
(-12.5,7.5)*+{\Hi_{\vf}}="H";
(12.5,-7.5)*+{\mL}="L";
(-12.5,-7.5)*+{\Hi_{\psi}}="I";
{\ar"H";"K"^{V_{\vf}^*}};
{\ar"I";"L"_{V_{\psi}^*}};
{\ar"K";"L"^{T}};
{\ar"H";"I"_{L_{T}}};
\endxy
,
\ee
the last of the conditions in Definition~\ref{defn:opstate},
requires an argument. First, to find the formula for
$V^*_{\vf}:\Hi_{\vf}\to\mK,$ let 
$\left[\sum_{i}a_{i}\otimes v_{i}\right]_{\vf}\in(\mA\otimes\mK)/\mN_{\vf}$ 
and $w\in\mK.$ Then 
\be
\label{eq:dualofVphi}
\begin{split}
&\left\<V^*_{\vf}\left[\sum_{i}a_{i}\otimes v_{i}\right]_{\vf},w\right\>_{\mK}
=\left\<\left[\sum_{i}a_{i}\otimes v_{i}\right]_{\vf},V_{\vf}w\right\>_{\vf}\\
&\qquad=\left\<\left[\sum_{i}a_{i}\otimes v_{i}\right]_{\vf},[1_{\mA}\otimes w]_{\vf}\right\>_{\vf}
=\sum_{i}\<v_{i},\vf(a_{i}^*)w\>_{\mK}
=\left\<\sum_{i}\vf(a_{i})v_{i},w\right\>_{\mK}
\end{split}
\ee
entails the formula 
\be
\label{eq:Vdualformula}
V^*_{\vf}\left[\sum_{i}a_{i}\otimes v_{i}\right]_{\vf}
=\sum_{i}\vf(a_{i})v_{i}.
\ee 
Commutativity of the diagram (\ref{eq:stinemorphismdual}) then
follows from%
\footnote{This is where commutativity of (\ref{eq:opstatesdiagram}) 
in Definition~\ref{defn:opstate}
is needed (cf.~Remark~\ref{rmk:morphismopstates}).
When restricting to operator states,  
(\ref{eq:opstatemorphismother}) would be too weak, and 
the diagram (\ref{eq:stinemorphismdual})
might not commute. 
}
\be
\label{eq:commutativityofVphidualVpsudualT}
\xy0;/r.25pc/:
(-25,15)*+{V_{\psi}^*\left(L_{T}\left(\big[\sum_{i}a_{i}\otimes v_{i}\big]_{\vf}\right)\right)}="1";
(-40,0)*+{V_{\psi}^*\left(\big[\sum_{i}a_{i}\otimes T(v_{i})\big]_{\psi}\right)}="2";
(-25,-15)*+{\sum_{i}\psi(a_{i}) T(v_{i})}="3";
(25,-15)*+{\sum_{i}T\big(\vf(a_{i}) v_{i}\big)}="4";
(40,0)*+{T\big(\sum_{i}\vf(a_{i})v_{i}\big)}="5";
(25,15)*+{T\left(V_{\vf}^*\left(\big[\sum_{i}a_{i}\otimes v_{i}\big]_{\vf}\right)\right)}="6";
{\ar@{=}@/_1.0pc/"1";"2"_(0.7){\text{(\ref{eq:stineT})}}};
{\ar@{=}@/_1.0pc/"2";"3"_(0.35){\text{(\ref{eq:Vdualformula})}}};
{\ar@{=}@/_1.0pc/"3";"4"_{\text{(\ref{eq:opstatesdiagram})}}};
{\ar@{=}@/_1.0pc/"4";"5"};
{\ar@{=}@/_1.0pc/"5";"6"_(0.3){\text{(\ref{eq:Vdualformula})}}};
\endxy
.
\ee
This concludes the definition of the assignment 
\be
\label{eq:stineonmorphisms}
\begin{split}
\mathbf{OCP}(\mA)_{1}&\xrightarrow{\mathbf{Stine}_{\mA}}\mathbf{AnRep}(\mA)_{1}\\
\Big((\mK,\vf)\xrightarrow{T}(\mL,\psi)\Big)&\mapsto\Big((\mK,\Hi_{\vf},\pi_{\vf},V_{\vf})\xrightarrow{(T,L_{T})}(\mL,\Hi_{\psi},\pi_{\psi},V_{\psi})\Big)
\end{split}
\ee
on morphisms of OCP maps on $\mA.$ 
Note that $\mathbf{Stine}_{\mA}(\id_{\mathcal{K}})$ equals $\left(\id_{\mathcal{K}},\id_{\mathcal{H}_{\vf}}\right)$ for any OCP map $(\mK,\vf).$ 
Furthermore, for a composable pair 
$(\mK,\vf)\xrightarrow{T}(\mL,\psi)\xrightarrow{S}(\mM,\chi)$
of morphisms of OCP maps, 
the diagram
\be
\xy0;/r.25pc/:
(12.5,-7.5)*+{\Hi_{\vf}}="1";
(0,7.5)*+{\Hi_{\psi}}="2";
(-12.5,-7.5)*+{\Hi_{\chi}}="3";
{\ar"1";"2"_{L_{T}}};
{\ar"2";"3"_{L_{S}}};
{\ar"1";"3"^{L_{S T}}};
\endxy
\ee
commutes.
Hence, 
(\ref{eq:stineonobjects}) and (\ref{eq:stineonmorphisms})
define a functor $\mathbf{Stine}_{\mA}:\mathbf{OCP}(\mA)\to\mathbf{AnRep}(\mA)$, 
which restricts to a functor 
$\mathbf{OpSt}(\mA)\to\mathbf{PAnRep}(\mA)$, also denoted by $\mathbf{Stine}_{\mA}$.

\item
Let $f:\mA'\to\mA$ be a $^*$-homomorphism of $C^*$-algebras. 
Two diagrams associated with the constructions preceding this
are given by
\be
\label{eq:restandStine}
\xy0;/r.25pc/:
(-15,10)*+{\mathbf{OCP}(\mA)}="1";
(15,10)*+{\mathbf{AnRep}(\mA)}="2";
(-15,-10)*+{\mathbf{OCP}(\mA')}="3";
(15,-10)*+{\mathbf{AnRep}(\mA')}="4";
{\ar"2";"1"_{\rest_{\mA}}};
{\ar"4";"3"^{\rest_{\mA'}}};
{\ar"1";"3"_{\mathbf{OCP}_{f}}};
{\ar"2";"4"^{\mathbf{AnRep}_{f}}};
\endxy
\quad\text{and}\quad
\xy0;/r.25pc/:
(-15,10)*+{\mathbf{OCP}(\mA)}="1";
(15,10)*+{\mathbf{AnRep}(\mA)}="2";
(-15,-10)*+{\mathbf{OCP}(\mA')}="3";
(15,-10)*+{\mathbf{AnRep}(\mA')}="4";
{\ar"1";"2"^{\mathbf{Stine}_{\mA}}};
{\ar"3";"4"_{\mathbf{Stine}_{\mA'}}};
{\ar"1";"3"_{\mathbf{OCP}_{f}}};
{\ar"2";"4"^{\mathbf{AnRep}_{f}}};
\endxy
\;,
\ee
and similarly for $\mathbf{OpSt}$ and $\mathbf{AnRep}.$ 
The diagram of functors on the left commutes (on the nose) by Lemma~\ref{lem:restnatural}. 
However, the diagram on the right does not
(this is analogous to what happens in the GNS construction---cf.~Construction~3.3 and diagram~(3.21) in~\cite{PaGNS}).%
\footnote{As a simple example, let $\vf:\mA\stoch\mB(\mK)$ be an operator state with $\mathcal{N}_{\vf}=\{0\}$ and set $\mA'=\C.$ Then there is only a single (unital)
$^{*}$-homomorphism $f:\C\to\mA$ and $\vf$ gets pulled back to the map $\vf\circ f$ that sends $1\in\C$ to $\id_{\mK}.$ In this case, $\mathcal{N}_{\vf\circ f}=\{0\}$, which entails $\Hi_{\vf\circ f}\cong\C\otimes\mK\cong\mK$, while $\Hi_{\vf}\cong\mA\otimes\mK.$ 
}
Nevertheless, there is a natural transformation 
\be
\label{eq:Stinef}
\xy0;/r.25pc/:
(-15,10)*+{\mathbf{OCP}(\mA)}="1";
(15,10)*+{\mathbf{AnRep}(\mA)}="2";
(-15,-10)*+{\mathbf{OCP}(\mA')}="3";
(15,-10)*+{\mathbf{AnRep}(\mA')}="4";
{\ar"1";"2"^{\mathbf{Stine}_{\mA}}};
{\ar"3";"4"_{\mathbf{Stine}_{\mA'}}};
{\ar"1";"3"_{\mathbf{OCP}_{f}}};
{\ar"2";"4"^{\mathbf{AnRep}_{f}}};
{\ar@{=>}"3";"2"|-{\mathbf{Stine}_{f}}};
\endxy
\ee
defined as follows. 
Given an OCP map $(\mK,\vf)$ on $\mA,$ applying 
$\mathbf{OCP}_{f}$ followed by $\mathbf{Stine}_{\mA'}$ 
provides the Stinespring anchored representation
$(\mK,\Hi_{\vf\circ f},\pi_{\vf\circ f},V_{\vf\circ f})$
on $\mA'$,
while applying $\mathbf{Stine}_{\mA}$ followed by $\mathbf{AnRep}_{f}$ 
gives $(\mK,\Hi_{\vf},\pi_{\vf}\circ f,V_{\vf}).$ 
The morphism $\mathbf{Stine}_{f}(\mK,\vf)$ 
from the first to the latter 
is given
by $(\id_{\mK},L_{f}),$ where $L_{f}:\Hi_{\vf\circ f}\to\Hi_{\vf}$ is defined as the unique map associated to 
\be
\label{eq:Lf}
\mA'\otimes\mK\xrightarrow{L'_{f}:=f\otimes\id_{\mK}}\mA\otimes\mK
\ee
via Lemma~\ref{lem:XmodN}. This lemma applies because if
$\xi=\sum_{i}a'_{i}\otimes v_{i}\in\mA'\otimes\mK,$
then
\be
\label{eq:Lfisometry}
\begin{split}
\llangle L'_{f}(\xi),L'_{f}(\xi)\rrangle_{\vf}
&=\sum_{i,j}\Big\<v_{i},\vf\big(f(a'_{i})^*f(a'_{j})\big)v_{j}\Big\>_{\mK}\quad\text{by definition of $\llangle\;\cdot\;,\;\cdot\;\rrangle_{\vf}$}\\
&=\sum_{i,j}\Big\<v_{i},\vf\big(f({a'_{i}}^*a'_{j})\big)v_{j}\Big\>_{\mK}\quad\text{since $f$ is a $^*$-homomorphism}\\
&=\llangle\xi,\xi\rrangle_{\vf\circ f}\quad\text{by definition of $\llangle\;\cdot\;,\;\cdot\;\rrangle_{\vf\circ f}$}.
\end{split}
\ee 
In fact, this calculation shows $L_{f}:\Hi_{\vf\circ f}\to\Hi_{\vf}$ is an isometry.
The requirements of a morphism of anchored representations all hold as well. 
Furthermore, $\mathbf{Stine}_{f}$ is a natural transformation
because for any morphism $(\mK,\vf)\xrightarrow{T}(\mL,\psi)$ of
operator states on $\mA,$ the diagram 
\be
\xy0;/r.25pc/:
(-30,7.5)*+{(\mK,\Hi_{\vf\circ f},\pi_{\vf\circ f},V_{\vf\circ f})}="1";
(30,7.5)*+{(\mK,\Hi_{\vf},\pi_{\vf}\circ f,V_{\vf})}="2";
(-30,-7.5)*+{(\mL,\Hi_{\psi\circ f},\pi_{\psi\circ f},V_{\psi\circ f})}="3";
(30,-7.5)*+{(\mL,\Hi_{\psi},\pi_{\psi}\circ f,V_{\psi})}="4";
{\ar"1";"2"^{(\id_{\mK},L_{f})}};
{\ar"2";"4"^{(T,L_{T})}};
{\ar"1";"3"_{(T,L_{T})}};
{\ar"3";"4"_{(\id_{\mL},L_{f})}};
\endxy 
\ee
commutes. Finally, $\mathbf{Stine}_{f}:\mathbf{Stine}_{\mA'}\circ\mathbf{OCP}_{f}\Rightarrow\mathbf{AnRep}_{f}\circ\mathbf{Stine}_{\mA}$ restricts to a natural transformation
$\mathbf{Stine}_{f}:\mathbf{Stine}_{\mA'}\circ\mathbf{OpSt}_{f}\Rightarrow\mathbf{PAnRep}_{f}\circ\mathbf{Stine}_{\mA}$
since $\mathbf{Stine}_{f}(\mK,\vf)=(\id_{\mK},L_{f})$ consists of two isometries by (\ref{eq:Lfisometry}). 

\item
Oplax-naturality of $\mathbf{Stine}$ holds because
$\mathbf{Stine}_{\id_{\mA}}$ is the identity natural transformation
for every $C^*$-algebra $\mA$ 
and 
the two natural transformations
(after composition in the diagram on the right)
\be
\label{eq:Stineoplaxnaturalconditions}
\xy0;/r.25pc/:
(-15,10)*+{\mathbf{OCP}(\mA)}="1";
(15,10)*+{\mathbf{AnRep}(\mA)}="2";
(-15,-10)*+{\mathbf{OCP}(\mA'')}="3";
(15,-10)*+{\mathbf{AnRep}(\mA'')}="4";
{\ar"1";"2"^{\mathbf{Stine}_{\mA}}};
{\ar"3";"4"_{\mathbf{Stine}_{\mA''}}};
{\ar"1";"3"|-{\mathbf{OCP}_{f\circ f'}}};
{\ar"2";"4"|-{\mathbf{AnRep}_{f\circ f'}}};
{\ar@{=>}"3";"2"|-{\mathbf{Stine}_{f\circ f'}}};
\endxy
\quad\&\quad
\xy0;/r.25pc/:
(-15,20)*+{\mathbf{OCP}(\mA)}="1";
(15,20)*+{\mathbf{AnRep}(\mA)}="2";
(-15,0)*+{\mathbf{OCP}(\mA')}="3";
(15,0)*+{\mathbf{AnRep}(\mA')}="4";
(-15,-20)*+{\mathbf{OCP}(\mA'')}="5";
(15,-20)*+{\mathbf{AnRep}(\mA'')}="6";
{\ar"1";"2"^{\mathbf{Stine}_{\mA}}};
{\ar"3";"4"_{\mathbf{Stine}_{\mA'}}};
{\ar"1";"3"|-{\mathbf{OCP}_{f}}};
{\ar"2";"4"|-{\mathbf{AnRep}_{f}}};
{\ar@{=>}"3";"2"|-{\mathbf{Stine}_{f}}};
{\ar"5";"6"_{\mathbf{Stine}_{\mA''}}};
{\ar"3";"5"|-{\mathbf{OCP}_{f'}}};
{\ar"4";"6"|-{\mathbf{AnRep}_{f'}}};
{\ar@{=>}"5";"4"|-{\mathbf{Stine}_{f'}}};
\endxy
\ee
are equal
for every pair of composable $^*$-homomorphisms
$\mA''\xrightarrow{f'}\mA'\xrightarrow{f}\mA$. 

\item
Fix a $C^*$-algebra $\mA$ and let 
$(\mK,\Hi,\pi,V)$ be an anchored representation of $\mA.$
Applying the functor $\rest_{\mA}$ followed by $\mathbf{Stine}_{\mA}$ 
to this representation gives
\be
(\mK,\Hi,\pi,V)\xmapsto{\rest_{\mA}}(\mK,\mathrm{Ad}_{V^*}\circ\pi)
\xmapsto{\mathbf{Stine}_{\mA}}
(\mK,\Hi_{\mathrm{Ad}_{V^*}\circ\pi},\pi_{\mathrm{Ad}_{V^*}\circ\pi},V_{\mathrm{Ad}_{V^*}\circ\pi}).
\ee
Set $m'_{\pi,V}:\mA\otimes\mK\to\Hi$ to be the composite  
\be
\label{eq:mpiV}
\begin{split}
\mA\otimes\mK&\xrightarrow{\pi\otimes V}\mB(\Hi)\otimes\Hi\to\Hi\\
\sum_{i}a_{i}\otimes v_{i}&\xmapsto{\;\qquad\quad m'_{\pi,V}\quad\qquad\;}\sum_{i}\pi(a_{i})V(v_{i}),
\end{split}
\ee
where the second map is the canonical action of $\mB(\Hi)$ on $\Hi.$ 
If
$\xi=\sum_{i}a_{i}\otimes v_{i}\in\mA\otimes\mK,$
then
\be
\label{eq:mpiVisometry}
\big\<m'_{\pi,V}(\xi),m'_{\pi,V}(\xi)\big\>_{\Hi}
=\sum_{i,j}\big\<v_{i},V^*\pi(a_{i}^*a_{j})Vv_{j}\big\>_{\mK}\\
=\llangle\xi,\xi\rrangle_{\mathrm{Ad}_{V^*}\circ\pi}.
\ee
Thus, $m'_{\pi,V}$ defines an isometry $m_{\pi,V}:\Hi_{\mathrm{Ad}_{V^*}\circ\pi}\to\Hi$ by Lemma~\ref{lem:XmodN}. 
Note that $m_{\pi,V}$ is an isometry even though $(\mK,\Hi,\pi,V)$ need not be a \emph{preserving} anchored representation. 
Commutativity of the diagrams 
\be
\xy0;/r.25pc/:
(-17.5,7.5)*+{\Hi_{\mathrm{Ad}_{V^*}\circ\pi}}="H1";
(17.5,7.5)*+{\Hi_{\mathrm{Ad}_{V^*}\circ\pi}}="H2";
(-17.5,-7.5)*+{\Hi}="I1";
(17.5,-7.5)*+{\Hi}="I2";
{\ar"H1";"H2"^{\pi_{\mathrm{Ad}_{V^*}\circ\pi}(a)}};
{\ar"I1";"I2"_{\pi(a)}};
{\ar"H1";"I1"_{m_{\pi,V}}};
{\ar"H2";"I2"^{m_{\pi,V}}};
\endxy
\qquad\forall\;a\in\mA,
\ee
\be
\xy0;/r.25pc/:
(-15,7.5)*+{\mK}="K";
(15,7.5)*+{\Hi_{\mathrm{Ad}_{V^*}\circ\pi}}="H";
(-15,-7.5)*+{\mK}="L";
(15,-7.5)*+{\Hi}="I";
{\ar"K";"H"^{V_{\mathrm{Ad}_{V^*}\circ\pi}}};
{\ar"L";"I"_{V}};
{\ar"K";"L"_{\id_{\mK}}};
{\ar"H";"I"^{m_{\pi,V}}};
\endxy
\;\;,\quad\text{and }\quad
\xy0;/r.25pc/:
(15,7.5)*+{\mK}="K";
(-15,7.5)*+{\Hi_{\mathrm{Ad}_{V^*}\circ\pi}}="H";
(15,-7.5)*+{\mK}="L";
(-15,-7.5)*+{\Hi}="I";
{\ar"H";"K"^{V^*_{\mathrm{Ad}_{V^*}\circ\pi}}};
{\ar"I";"L"_{V^*}};
{\ar"K";"L"^{\id_{\mK}}};
{\ar"H";"I"_{m_{\pi,V}}};
\endxy
\ee
in the definition of a morphism 
$(\mK,\Hi_{\mathrm{Ad}_{V^*}\circ\pi},\pi_{\mathrm{Ad}_{V^*}\circ\pi},V_{\mathrm{Ad}_{V^*}\circ\pi})\xrightarrow{(\id_{\mK},m_{\pi,V})}(\mK,\Hi,\pi,V)$
of anchored representations
follow directly from the definitions
(for the last diagram, apply (\ref{eq:Vdualformula}) to 
$\vf:=\mathrm{Ad}_{V^*}\circ\pi$). 
Set $m_{\mA}$ to be the assignment 
\be
\label{eq:mA}
\begin{split}
\mathbf{AnRep}(\mA)_{0}&\xrightarrow{m_{\mA}}\mathbf{AnRep}(\mA)_{1}\\
(\mK,\Hi,\pi,V)&\xmapsto{\quad\;}(\id_{\mK},m_{\pi,V})
\end{split}
\ee
from objects of $\mathbf{AnRep}(\mA)$ to morphisms of
$\mathbf{AnRep}(\mA).$
Notice that $m_{\mA}$ restricts to a well-defined assignment 
$\mathbf{PAnRep}(\mA)_{0}\xrightarrow{m_{\mA}}\mathbf{PAnRep}(\mA)_{1}$ by (\ref{eq:mpiVisometry}). 
To see naturality of 
\be
\xy0;/r.25pc/:
(-20,-7.5)*+{\mathbf{AnRep}(\mA)}="1";
(0,7.5)*+{\mathbf{OCP}(\mA)}="2";
(20,-7.5)*+{\mathbf{AnRep}(\mA)}="3";
{\ar"1";"2"^{\rest_{\mA}}};
{\ar"2";"3"^{\mathbf{Stine}_{\mA}}};
{\ar"1";"3"_{\id_{\mathbf{AnRep}(\mA)}}};
{\ar@{=>}"2";(0,-6.5)^(0.6){m_{\mA}}};
\endxy
,
\ee
let
$(\mK,\Hi,\pi,V)\xrightarrow{(T,L)}(\mL,\mI,\rho,W)$
be a morphism of anchored representations. Then, the diagram
\be
\xy0;/r.25pc/:
(-30,7.5)*+{(\mK,\Hi_{\mathrm{Ad}_{V^*}\circ\pi},\pi_{\mathrm{Ad}_{V^*}\circ\pi},V_{\mathrm{Ad}_{V^*}\circ\pi})}="H1";
(30,7.5)*+{(\mK,\Hi,\pi,V)}="H2";
(-30,-7.5)*+{(\mL,\Hi_{\mathrm{Ad}_{W^*}\circ\rho},\pi_{\mathrm{Ad}_{W^*}\circ\rho},V_{\mathrm{Ad}_{W^*}\circ\rho})}="I1";
(30,-7.5)*+{(\mL,\mI,\rho,W)}="I2";
{\ar"H1";"H2"^(0.6){(\id_{\mK},m_{\pi,V})}};
{\ar"I1";"I2"_(0.6){(\id_{\mL},m_{\rho,W})}};
{\ar"H1";"I1"_{(T,L_{T})}};
{\ar"H2";"I2"^{(T,L)}};
\endxy
\ee
commutes by conditions 
(\ref{eq:intertwiner}) and 
(\ref{eq:morphismanchoredrepns})
in the definition of a morphism of anchored representations. 

\item
To see that the assignment sending a $C^*$-algebra $\mA$ to $m_{\mA}$ 
defines a modification%
\footnote{This composite of $\rest$ followed by $\mathbf{Stine}$ along the top
two double arrows 
is denoted by $\begin{matrix}[0.7]\rest\\\mathbf{Stine}\end{matrix}.$
}
\be
\label{eq:modification}
\xy0;/r.25pc/:
(-15,-7.5)*+{\mathbf{AnRep}}="1";
(0,7.5)*+{\mathbf{OCP}}="2";
(15,-7.5)*+{\mathbf{AnRep}}="3";
{\ar@{=>}"1";"2"^{\rest}};
{\ar@{=>}"2";"3"^{\mathbf{Stine}}};
{\ar@{=>}"1";"3"_{\id_{\mathbf{AnRep}}}};
{\ar@3{->}(0,3);(0,-6)_{m}};
\endxy
\ee
of oplax-natural transformations, 
for every morphism $f:\mA'\to\mA$ of $C^*$-algebras,
\be
\begin{split}
\hspace{-7mm}
\xy0;/r.23pc/:
(24,10)*+{\mathbf{AnRep}(\mA)}="1";
(-24,10)*+{\mathbf{AnRep}(\mA)}="2";
(24,-10)*+{\mathbf{AnRep}(\mA')}="3";
(-24,-10)*+{\mathbf{AnRep}(\mA')}="4";
(0,5)*+{\mathbf{OCP}(\mA)}="5";
(0,-15)*+{\mathbf{OCP}(\mA')}="6";
{\ar@/_0.5pc/"5";"1"^(0.505){\mathbf{Stine}_{\mA}}};
{\ar@/_0.5pc/"6";"3"_(0.60){\mathbf{Stine}_{\mA'}}};
{\ar"1";"3"|-{\mathbf{AnRep}_{f}}};
{\ar"2";"4"|-{\mathbf{AnRep}_{f}}};
{\ar@/_0.5pc/"2";"5"^(0.5125){\rest_{\mA}}};
{\ar@/_0.5pc/"4";"6"_(0.40){\rest_{\mA'}}};
{\ar"5";"6"|-{\mathbf{OCP}_{f}}};
{\ar@/^1.5pc/"2";"1"^{\id_{\mathbf{AnRep}(\mA)}}};
{\ar@{=>}"6"+(7,5);"1"+(-6,-6)|-{\mathbf{Stine}_{f}}};
{\ar@{=}"4"+(7,4);"5"+(-6,-4)|-{\id=\rest_{f}}};
{\ar@{=>}"5";(0,15)|-(0.4){m_{\mA}}};
\endxy
\hspace{-2mm}
=
\hspace{-2mm}
\xy0;/r.23pc/:
(-24,-10)*+{\mathbf{AnRep}(\mA')}="1";
(0,-15)*+{\mathbf{OpSt}(\mA')}="2";
(24,-10)*+{\mathbf{AnRep}(\mA')}="3";
(-24,10)*+{\mathbf{AnRep}(\mA)}="4";
(24,10)*+{\mathbf{AnRep}(\mA)}="5";
{\ar@/_0.5pc/"1";"2"_(0.40){\rest_{\mA'}}};
{\ar@/_0.5pc/"2";"3"_(0.60){\mathbf{Stine}_{\mA'}}};
{\ar"4";"1"|-{\mathbf{AnRep}_{f}}};
{\ar"5";"3"|-{\mathbf{AnRep}_{f}}};
{\ar@/^1.5pc/"4";"5"^{\id_{\mathbf{AnRep}(\mA)}}};
{\ar@/^1.5pc/"1";"3"|-{\id_{\mathbf{AnRep}(\mA')}}};
{\ar@{=}@/^0.5pc/"1"+(7,6);"5"+(-9,-4)|-{\id_{\mathbf{AnRep}_{f}}}};
{\ar@{=>}"2";(0,-6)|-(0.4){m_{\mA'}}};
\endxy
\!\!,
\end{split}
\raisetag{-3.75\baselineskip}
\ee
i.e.~for every object $(\mK,\Hi,\pi,V)$
of $\mathbf{AnRep}(\mA)$ with 
$\vf:=\mathrm{Ad}_{V^*}\circ\pi,$
the diagram
\be
\hspace{-3mm}
\xy0;/r.25pc/:
(-40,-7.5)*+{\big(\mK,\Hi_{\vf\circ f},\pi_{\vf\circ f},V_{\vf\circ f}\big)}="1";
(0,7.5)*+{\big(\mK,\Hi_{\vf},\pi_{\vf}\circ f,V_{\vf}\big)}="2";
(40,-7.5)*+{(\mK,\Hi,\pi\circ f,V)}="3";
{\ar"1";"2"^(0.4){(\id_{\mK},L_{f})=\mathbf{Stine}_{f}(\mK,\vf)\qquad}};
{\ar"2";"3"^(0.65){\;\;\;\;\;\;\;\;\qquad\mathbf{AnRep}_{f}\big(m_{\mA}
(\mK,\Hi,\pi,V)\big)=(\id_{\mK},m_{\pi,V})}};
{\ar"1";"3"_{m_{\mA'}(\mK,\Hi,\pi\circ f,V)=(\id_{\mK},m_{\pi\circ f,V})}};
\endxy
\ee
of morphisms of anchored representations of
$\mA'$ must commute. This follows directly from the definitions. 
Hence, $m$ is a modification, which also restricts to a modification
when working with $\mathbf{OpSt}$ and $\mathbf{PAnRep}$ by (\ref{eq:mpiVisometry}). 

\item
To see that
\be
\label{eq:restStine}
\xy0;/r.25pc/:
(-15,-7.5)*+{\mathbf{OCP}}="1";
(0,7.5)*+{\mathbf{AnRep}}="2";
(15,-7.5)*+{\mathbf{OCP}}="3";
{\ar@{=>}"1";"2"^{\mathbf{Stine}}};
{\ar@{=>}"2";"3"^{\rest}};
{\ar@{=>}"1";"3"_{\id_{\mathbf{OCP}}}};
\endxy
\ee
commutes, first fix a $C^*$-algebra $\mA.$ 
Commutativity in (\ref{eq:restStine}) requires that  
\be
\label{eq:restStineisid}
\xy0;/r.25pc/:
(-15,-7.5)*+{\mathbf{OCP}({\mA})}="1";
(0,7.5)*+{\mathbf{AnRep}({\mA})}="2";
(15,-7.5)*+{\mathbf{OCP}({\mA})}="3";
{\ar"1";"2"^{\mathbf{Stine}({\mA})}};
{\ar"2";"3"^{\rest_{\mA}}};
{\ar"1";"3"_{\id_{\mathbf{OCP}({\mA})}}};
\endxy
\ee
must commute. 
On objects, this translates to 
$\vf=\mathrm{Ad}_{V^*_{\vf}}\circ\pi_{\vf}$ for 
every OCP map $(\mK,\vf)$ on $\mA$, 
which follows from the definitions of $V_{\vf}$ and $\pi_{\vf}$.
Since a morphism $(\mK,\vf)\xrightarrow{T}(\mL,\psi)$ is unchanged
under $\rest_{\mA}\circ\mathbf{Stine}_{\mA},$  (\ref{eq:restStineisid}) commutes. 

Commutativity of (\ref{eq:restStine}) also requires that for 
every $^*$-homomorphism $f:\mA'\to\mA,$ 
\be
\begin{split}
\hspace{-4mm}
\xy0;/r.24pc/:
(24,10)*+{\mathbf{OCP}(\mA)}="1";
(-24,10)*+{\mathbf{OCP}(\mA)}="2";
(24,-10)*+{\mathbf{OCP}(\mA')}="3";
(-24,-10)*+{\mathbf{OCP}(\mA')}="4";
(0,5)*+{\mathbf{AnRep}(\mA)}="5";
(0,-15)*+{\mathbf{AnRep}(\mA')}="6";
{\ar@/_0.5pc/"5";"1"^(0.505){\mathbf{rest}_{\mA}}};
{\ar@/_0.5pc/"6";"3"_(0.60){\mathbf{rest}_{\mA'}}};
{\ar"1";"3"|-{\mathbf{OpSt}_{f}}};
{\ar"2";"4"|-{\mathbf{OpSt}_{f}}};
{\ar@/_0.5pc/"2";"5"^(0.45){\mathbf{Stine}_{\mA}}};
{\ar@/_0.5pc/"4";"6"_(0.40){\mathbf{Stine}_{\mA'}}};
{\ar"5";"6"|-{\mathbf{AnRep}_{f}}};
{\ar@/^1.5pc/"2";"1"^{\id_{\mathbf{OCP}(\mA)}}};
{\ar@{=}"6"+(7,5);"1"+(-6,-6)|-{\id}};
{\ar@{=>}"4"+(7,4);"5"+(-6,-4)|-{\mathbf{Stine}_{f}}};
{\ar@{=}"5";(0,15)|-(0.4){\id}};
\endxy
\hspace{-2mm}
=
\hspace{-2mm}
\xy0;/r.24pc/:
(-24,-10)*+{\mathbf{OCP}(\mA')}="1";
(0,-15)*+{\mathbf{AnRep}(\mA')}="2";
(24,-10)*+{\mathbf{OCP}(\mA')}="3";
(-24,10)*+{\mathbf{OCP}(\mA)}="4";
(24,10)*+{\mathbf{OCP}(\mA)}="5";
{\ar@/_0.5pc/"1";"2"_(0.40){\mathbf{Stine}_{\mA'}}};
{\ar@/_0.5pc/"2";"3"_(0.60){\rest_{\mA'}}};
{\ar"4";"1"|-{\mathbf{OCP}_{f}}};
{\ar"5";"3"|-{\mathbf{OCP}_{f}}};
{\ar@/^1.5pc/"4";"5"^{\id_{\mathbf{OCP}(\mA)}}};
{\ar@/^1.5pc/"1";"3"|-{\id_{\mathbf{OCP}(\mA')}}};
{\ar@{=}@/^0.5pc/"1"+(7,6);"5"+(-9,-4)|-{\id_{\mathbf{OCP}_{f}}}};
{\ar@{=}"2";(0,-6)|-(0.4){\id}};
\endxy
\!\!,
\end{split}
\raisetag{-3.75\baselineskip}
\ee
i.e.~to every OCP map $(\mK,\vf),$ the diagram 
\be
\xy0;/r.25pc/:
(-35,-7.5)*+{\big(\mK,\mathrm{Ad}_{V^*_{\vf\circ f}}\circ\pi_{\vf\circ f}\big)}="1";
(0,7.5)*+{\big(\mK,\mathrm{Ad}_{V^*_{\vf}}\circ\pi_{\vf}\circ f\big)}="2";
(35,-7.5)*+{(\mK,\vf\circ f)}="3";
{\ar"1";"2"^(0.4){\id_{\mK}=\rest_{\mA'}(\mathbf{Stine}_{f}(\mK,\vf))\;\qquad}};
{\ar"2";"3"^(0.65){\;\;\quad\mathbf{OCP}_{f}\left(\id_{\mK}\right)=\id_{\mK}}};
{\ar"1";"3"_{\id_{\mK}}};
\endxy
\ee
of morphisms of OCP maps on $\mA'$ must commute, 
which it clearly does. Commutativity of (\ref{eq:restStine}) with 
$\mathbf{OpSt}$ and $\mathbf{PAnRep}$ follows from this as well. 

\item
By the final remark in the appendix of \cite{PaGNS}, it suffices to prove
\be
\label{eq:firstzigzag}
\hspace{-4mm}
\xy0;/r.25pc/:
(0,22.5)*+{\mathbf{AnRep}(\mA)}="1";
(0,7.5)*+{\mathbf{OCP}(\mA)}="2";
(0,-7.5)*+{\mathbf{AnRep}(\mA)}="3";
(0,-22.5)*+{\mathbf{OCP}(\mA)}="4";
{\ar"1";"2"|-{\rest_{\mA}}};
{\ar"2";"3"|-{\mathbf{Stine}_{\mA}}};
{\ar"3";"4"|-{\rest_{\mA}}};
{\ar@/_4.5pc/"2";"4"_{\id_{\mathbf{OCP}(\mA)}}};
{\ar@/^4.5pc/"1";"3"^{\id_{\mathbf{AnRep}(\mA)}}};
{\ar@2{->}(-17.5,-7.5);"3"^(0.175){\id}};
{\ar@2{->}"2";(17.5,7.5)^(0.75){m_{\mA}}};
\endxy
\
=
\quad
\xy0;/r.25pc/:
(0,10)*+{\mathbf{AnRep}(\mA)}="1";
(0,-10)*+{\mathbf{OCP}(\mA)}="2";
{\ar@/_2.25pc/"1";"2"_{\rest_{\mA}}};
{\ar@/^2.25pc/"1";"2"^{\rest_{\mA}}};
{\ar@2{->}(-7.5,0);(7.5,0)^{\id_{\rest_{\mA}}}};
\endxy
\ee
and
\be
\label{eq:secondzigzag}
\hspace{-4mm}
\xy0;/r.25pc/:
(0,22.5)*+{\mathbf{OCP}(\mA)}="1";
(0,7.5)*+{\mathbf{AnRep}(\mA)}="2";
(0,-7.5)*+{\mathbf{OCP}(\mA)}="3";
(0,-22.5)*+{\mathbf{AnRep}(\mA)}="4";
{\ar"1";"2"|-{\mathbf{Stine}_{\mA}}};
{\ar"2";"3"|-{\rest_{\mA}}};
{\ar"3";"4"|-{\mathbf{Stine}_{\mA}}};
{\ar@/^4.5pc/"2";"4"^{\id_{\mathbf{AnRep}(\mA)}}};
{\ar@/_4.5pc/"1";"3"_{\id_{\mathbf{OCP}(\mA)}}};
{\ar@2{->}(-17.5,7.5);"2"^(0.175){\id}};
{\ar@2{->}"3";(17.5,-7.5)^(0.75){m_{\mA}}};
\endxy
\hspace{-2mm}
=
\quad
\xy0;/r.25pc/:
(0,10)*+{\mathbf{OCP}(\mA)}="1";
(0,-10)*+{\mathbf{AnRep}(\mA)}="2";
{\ar@/_2.25pc/"1";"2"_{\mathbf{Stine}_{\mA}}};
{\ar@/^2.25pc/"1";"2"^{\mathbf{Stine}_{\mA}}};
{\ar@2{->}(-7.5,0);(7.5,0)^{\id_{\mathbf{Stine}_{\mA}}}};
\endxy
\ee
for each object $\mA$ of $\CAlg^{\op}.$ 
The equality in (\ref{eq:firstzigzag}) follows from the equality 
\be
\rest_{\mA}\big(\underbrace{m_{\mA}(\mK,\Hi,\pi,V)}_{(\id_{\mK},m_{\pi,V})}\big)=\id_{\mK}
\ee
of morphisms of OCP maps from $(\mK,\mathrm{Ad}_{V^*}\circ\pi)$ to itself
for every anchored representation $(\mK,\Hi,\pi,V)$ of $\mA.$ 
To see the equality in (\ref{eq:secondzigzag}), 
consider an OCP map $(\mK,\vf).$ Applying 
$\mathbf{Stine}_{\mA}\circ\rest_{\mA}\circ\mathbf{Stine}_{\mA}$
to this gives $(\mK,\Hi_{\vf},\pi_{\vf},V_{\vf})$
because $\vf=\mathrm{Ad}_{V_{\vf}^*}\circ\pi_{\vf}.$  
In order for (\ref{eq:secondzigzag}) to hold, it should be the case that 
$m_{\mA}(\mK,\Hi_{\vf},\pi_{\vf},V_{\vf})=(\id_{\mK},\id_{\Hi_{\vf}})$ as morphisms
from $(\mK,\Hi_{\vf},\pi_{\vf},V_{\vf})$ to $(\mK,\Hi_{\vf},\pi_{\vf},V_{\vf})$
in $\mathbf{AnRep}(\mA).$ These are in fact equal because
for any element 
$\left[\sum_{i}a_{i}\otimes v_{i}\right]\in{(\mA\otimes\mK)/\mN_{\vf}},$
\be
m_{\pi_{\vf},V_{\vf}}\left(\left[\sum_{i}a_{i}\otimes v_{i}\right]\right)
=\sum_{i}\pi_{\vf}(a_{i}) V_{\vf}(v_{i})
=\sum_{i}\pi_{\vf}(a_{i}) [1_{\mA}\otimes v_{i}]
=\sum_{i}[a_{i}\otimes v_{i}].
\ee

\end{enumerate}
This proves that the quadruple $(\mathbf{Stine},\rest,\id,m)$ is
an adjunction in 
$\mathbf{Fun}(\CAlg^{\op},\tCAT)$, 
i.e.~$\mathbf{Stine}$ is left adjoint to $\rest.$ 
\eprf

\section{Stinespring dilations and their universal properties}
\label{sec:generalizations}

The left adjoint $\mathbf{Stine}$ to $\rest$ provides what is sometimes called
a \emph{minimal Stinespring representation/dilation} of an operator-valued completely positive map
(for comparison, see the discussion after Theorem~4.1 in Paulsen~\cite{Pa03}).
 
\bd
\label{defn:stinespringrepn}
A \define{Stinespring representation/dilation} of an OCP map
$(\mK,\vf)$ on $\mA$ is an anchored representation 
$(\mK,\Hi,\pi,V)$ of $\mA$ such that $\vf=\mathrm{Ad}_{V^*}\circ\pi.$ 
A Stinespring representation $(\mK,\Hi,\pi,V)$ is said to be \define{minimal}
iff 
\be
\pi(\mA)V(\mK):=\mathrm{span}\big\{\pi(a)V(v)\;:\;a\in\mA,\;v\in\mK\big\}
\ee
is dense in $\Hi.$ 
\ed

\br
\label{rmk:stinespringreps}
In terms of the functors and natural transformations introduced, 
$(\mK,\Hi,\pi,V)$ is a Stinespring representation of $(\mK,\vf)$ on $\mA$
if and only if $\rest_{\mA}(\mK,\Hi,\pi,V)=(\mK,\vf).$ 
Minimality will be addressed in Corollary~\ref{cor:minimalStinespring}. 
\er

The following corollaries explain the meaning of Theorem~\ref{thm:stinespring} more concretely. 

\bc
\label{cor:universalpropmorphismStine}
Let $(\mK,\vf)\xrightarrow{T}(\mathcal{L},\psi)$ be a morphism of
OCP maps on a $C^*$-algebra $\mA$ and let $(\mathcal{L},\mathcal{I},\rho,W)$ be any Stinespring representation of $\psi$. 
Then there exists a unique morphism $(\mK,\Hi_{\vf},\pi_{\vf},V_{\vf})\xrightarrow{\mathcal{T}}(\mathcal{L},\mathcal{I},\rho,W)$ such that $\rest_{\mA}(\mathcal{T})=T.$  
Furthermore, if $T$ is a morphism of operator states, then $\mathcal{T}$ is a morphism of preserving anchored representations. 
\ec

\bprf
Since $(\mathbf{Stine}_{\mA},\rest_{\mA},\id_{\mathbf{OpSt}(\mA)},m_{\mA})$ is an adjunction, existence and uniqueness follows from the universal property
of adjunctions in part~\ref{lem:adjunctioni} of Lemma~\ref{lem:universalpropertyadjunctions}. 
In the notation of that lemma, 
$c=(\mK,\vf),$ $d=(\mathcal{L},\mathcal{I},\rho,W),$ and $g=T.$ 
This unique morphism is given by 
\be
\mathcal{T}=\big(m_{\mA}(\mathcal{L},\mathcal{I},\rho,W)\big)\circ\mathbf{Stine}_{\mA}(T)=(\id_{\mL},m_{\rho,W})\circ(T,L_{T})=(T,m_{\rho,W} L_{T}).
\ee
Recall, $m_{\rho,W} L_{T}:\Hi_{\vf}\to\mathcal{I}$ is uniquely determined by 
the assignment
\be
\mA\otimes\mK/\mathcal{N}_{\vf}\ni\left[\sum_{i}a_{i}\otimes v_{i}\right]_{\vf}\mapsto\sum_{i}\rho(a_{i})W\big(T(v_{i})\big).
\ee
When $T$ is a morphism of operator states, both $T$ and $m_{\rho,W} L_{T}$ are isometries. Hence, $\mathcal{T}$ defines a morphism of preserving anchored representations. 
\eprf

The following is a special case of Corollary~\ref{cor:universalpropmorphismStine}. It will be used to provide the relationship to minimal Stinespring representations in Corollary~\ref{cor:minimalStinespring}. 

\bc
\label{cor:Stinespringinitial}
Let $(\mK,\vf)$ be an OCP map on a $C^*$-algebra $\mA$
and let $(\mK,\Hi,\pi,V)$ be any Stinespring representation of $\vf.$
Then there exists a unique morphism $(\mK,\Hi_{\vf},\pi_{\vf},V_{\vf})\xrightarrow{\mathcal{T}}(\mK,\Hi,\pi,V)$ such that $\rest_{\mA}(\mathcal{T})=\id_{\mK}.$  In fact, $\mathcal{T}=(\id_{\mK},m_{\pi,V})$,
and therefore consists of two isometries.
\ec

\bprf
This follows from Corollary~\ref{cor:universalpropmorphismStine} for $T=\id_{\mK}$. The fact that $\mathcal{T}$ consists of two isometries follows from the fact that $m_{\pi,V}$ is an isometry for any anchored representation $(\mK,\Hi,\pi,V)$ by (\ref{eq:mpiVisometry}). 
\eprf

\br
\label{rmk:Stinespringinitial}
If $(\mK,\vf)$ is an OCP map on $\mA,$ one can form a category of Stinespring representations of $(\mK,\vf).$ The objects of this category are Stinespring representations of $(\mK,\vf)$ and morphisms are of the form $(\mK,\Hi,\pi,V)\xrightarrow{(\id_{\mK},L)}(\mathcal{K},\mathcal{I},\rho,W)$.  
Corollary~\ref{cor:Stinespringinitial} then says that $\mathbf{Stine}_{\mA}(\mK,\vf)$ is an initial object in the category of Stinespring representations of $(\mK,\vf).$
\er

We now state the relationship between our adjunction and minimal Stinespring representations.

\bc
\label{cor:minimalStinespring}
Let $(\mK,\vf)$ be an OCP map on a $C^*$-algebra $\mA.$ 
Then $\mathbf{Stine}_{\mA}(\mK,\vf)\equiv(\mK,\Hi_{\vf},\pi_{\vf},V_{\vf})$ is a minimal Stinespring representation of $\vf.$ 
Conversely, given any minimal Stinespring representation $(\mK,\Hi,\pi,V)$ of $\vf,$ there exists a unique
isomorphism $(\mK,\Hi_{\vf},\pi_{\vf},V_{\vf})\to(\mK,\Hi,\pi,V)$ of anchored representations. In fact, this isomorphism consists of two unitaries. 
\ec

\bprf
{\color{white}{you found me!}}

\noindent
($\Rightarrow$)
By Remark~\ref{rmk:stinespringreps}, $(\mK,\Hi_{\vf},\pi_{\vf},V_{\vf})$ is a Stinespring representation due to part~\ref{vii} in the proof of Theorem~\ref{thm:stinespring}. By the construction in part~\ref{i}, 
\be
\mathrm{span}\big\{\pi_{\vf}(a)V_{\vf}(v)\;:\;a\in\mA,\;v\in\mK\big\}=\mathrm{span}\big\{[a\otimes v]_{\vf}\;:\;a\in\mA,\;v\in\mK\big\},
\ee
which is dense in $\Hi_{\vf}$ by definition (cf.~(\ref{eq:stinespringspace})). 
Hence, $(\mK,\Hi_{\vf},\pi_{\vf},V_{\vf})$ is minimal.

\noindent
($\Leftarrow$)
\sloppy
Conversely, let $(\mK,\Hi,\pi,V)$ be a minimal Stinespring representation
of $\vf.$ Then Corollary~\ref{cor:Stinespringinitial} provides the unique morphism 
$m_{\mA}(\mK,\Hi,\pi,V)\equiv(\id_{\mK},m_{\pi,V})$ from $(\mK,\Hi_{\vf},\pi_{\vf},V_{\vf})$ to $(\mK,\Hi,\pi,V)$, 
with $m_{\pi,V}$ an isometry. Furthermore, $m_{\pi,V}$ is unitary because the image of $\Hi_{\vf}$ under $m_{\pi,V}$ is exactly $\pi(\mA)V(\mK)$ by (\ref{eq:mpiV}) and equals $\Hi$ by the minimality assumption. 
\eprf

Theorem~\ref{thm:stinespring} entails the following with regards to morphisms and functoriality.

\bc
\label{cor:Stinespringmorphismsofalgebras}
Let $(\mK,\vf)$ be an OCP map on a $C^*$-algebra $\mA$ and let $f:\mA'\to\mA$ be a morphism of $C^*$-algebras. Let $(\mK,\Hi_{\vf\circ f},\pi_{\vf\circ f},V_{\vf\circ f})=\mathbf{Stine}_{\mA}(\mK,\vf\circ f)$ be the minimal Stinespring representation of $\vf\circ f$ and let $(\mK,\Hi_{\vf},\pi_{\vf}\circ f,V_{\vf})$ be the pull-back of the minimal Stinespring representation $\mathbf{Stine}_{\mA}(\mK,\vf)$ of $\vf$ under $f.$ 
Then 
$\mathbf{Stine}_{f}(\mK,\vf):(\mK,\Hi_{\vf\circ f},\pi_{\vf\circ f},V_{\vf\circ f})\to(\mK,\Hi_{\vf},\pi_{\vf}\circ f,V_{\vf})$ is the unique morphism in $\mathbf{AnRep}(\mA')$ such that $\rest_{\mA'}\left(\mathbf{Stine}_{f}(\mK,\vf)\right)=\id_{\mK}.$ In fact, 
\be
\label{eq:StinefismA}
\mathbf{Stine}_{f}(\mK,\vf)=m_{\mA'}\Big(\mathbf{AnRep}_{f}\big(\mathbf{Stine}_{\mA}(\mK,\vf)\big)\Big). 
\ee
Furthermore, 
the assignment $f\mapsto\mathbf{Stine}_{f}$ is functorial in the sense that 
\be
\label{eq:Stinefunctorial}
\mathbf{Stine}_{f\circ f'}(\mK,\vf)=\big(\mathbf{Stine}_{f}(\mK,\vf)\big)\circ\big(\mathbf{Stine}_{f'}(\mK,\vf\circ f)\big)
\ee
for all composable pairs $\mA''\xrightarrow{f'}\mA'\xrightarrow{f}\mA$ of morphisms of $C^*$-algebras and for all OCP maps $\vf:\mA\stoch\mB(\mK).$  
\ec

\bprf
Since $(\mK,\Hi_{\vf},\pi_{\vf}\circ f,V_{\vf})$ is a Stinespring representation of $(\mK,\vf\circ f),$ Corollary~\ref{cor:Stinespringinitial} says there is a unique morphism from $(\mK,\Hi_{\vf\circ f},\pi_{\vf\circ f},V_{\vf\circ f})$ to $(\mK,\Hi_{\vf},\pi_{\vf}\circ f,V_{\vf})$ that restricts to $\id_{\mK}$ under $\rest_{\mA'}.$ By construction (cf.~Equation~(\ref{eq:Stinef}) and what follows), $\mathbf{Stine}_{f}(\mK,\vf)$ is one such morphism. Similarly (cf.~Equations~(\ref{eq:mpiV})~and~(\ref{eq:mA})), 
$m_{\mA'}\Big(\mathbf{AnRep}_{f}\big(\mathbf{Stine}_{\mA}(\mK,\vf)\big)\Big)$ is another such morphism because the target of this morphism is 
\be
\mathbf{AnRep}_{f}\big(\mathbf{Stine}_{\mA}(\mK,\vf)\big)=(\mK,\Hi_{\vf},\pi_{\vf}\circ f,V_{\vf})
\ee
and the source of this morphism is 
\be
\mathbf{Stine}_{\mA'}\Big(\underbrace{\rest_{\mA'}\Big(\mathbf{AnRep}_{f}\big(\mathbf{Stine}_{\mA}(\mK,\vf)\big)\Big)}_{(\mK,\vf\circ f)}\Big)=(\mK,\Hi_{\vf\circ f},\pi_{\vf\circ f},V_{\vf\circ f}).
\ee
Equation~(\ref{eq:Stinefunctorial}) is precisely oplax-naturality of $\mathbf{Stine}$ from equality of the diagrams in~(\ref{eq:Stineoplaxnaturalconditions}).
\eprf

\br
Corollary~\ref{cor:Stinespringmorphismsofalgebras} also shows that the definition and oplaxness of $\mathbf{Stine}$ on morphisms of $C^*$-algebras is determined by (\ref{eq:StinefismA}),
which itself follows from the universal property of the adjunctions 
$(\mathbf{Stine}_{\mA},\rest_{\mA},\id_{\mathbf{OpSt}(\mA)},m_{\mA})$
over all $C^*$-algebras $\mA.$ 
This is discussed in more abstract form in Proposition~\ref{prop:adjunctionsinfunctorcategories} and Remark~\ref{rmk:simplifiedStinespringproof}. 
\er

\br
\label{rmk:Westerbaan}
In \cite{WeWe16}, Westerbaan and Westerbaan provide a universal property for Paschke dilations (a generalization of Stinespring dilations) of \emph{normal} completely positive maps $\vf:\mA\stoch\mB,$ where $\mA$ and $\mB$ are \emph{von Neumann algebras}. We refer the reader to their work and the references therein for any terminology not explained here. Briefly, a \define{Paschke dilation} of $\mA\xstoch{\vf}\mB$ is a triple $(\mathcal{P},\rho,\psi)$, which consists of a von Neumann algebra $\mathcal{P},$ a normal $^*$-homomorphism $\rho:\mA\to\mathcal{P},$ and a normal CP map $\psi:\mathcal{P}\stoch\mB$ such that $\vf=\psi\circ\rho$, satisfying the following universal property: for any other triple $(\mathcal{P}',\rho',\psi')$ with $\mathcal{P}'$ a von Neumann algebra, $\rho':\mA\to\mathcal{P}'$ a normal $^*$-homomorphism, and $\psi':\mathcal{P}'\stoch\mB$ a normal CP map such that $\vf=\psi'\circ\rho',$ there exists a unique normal CP map $\s:\mathcal{P}'\stoch\mathcal{P}$ such that the diagram
\be
\xy0;/r.25pc/:
(0,0)*+{\mathcal{P}}="P";
(0,-15)*+{\mathcal{P}'}="Pp";
(-20,15)*+{\mA}="A";
(20,15)*+{\mB}="B";
{\ar@{~>}"A";"B"^{\vf}|(0.95){\hole}};
{\ar"A";"P"^{\rho}};
{\ar@{~>}"P";"B"^{\psi}|(0.91){\hole}};
{\ar@{~>}@/_1.25pc/"Pp";"B"_{\psi'}|(0.949){\hole}};
{\ar@/_1.25pc/"A";"Pp"_{\rho'}};
{\ar@{~>}"Pp";"P"_{\s}|(0.87){\hole}};
\endxy
\ee
commutes. In Theorem~14 of \cite{WeWe16}, 
Westerbaan and Westerbaan prove
that a minimal Stinespring representation $(\mK,\Hi,\pi,V)$ of a \emph{normal} CP map $\vf:\mA\stoch\mB(\mK)$ induces a Paschke dilation given by the triple%
\footnote{Note that our conventions for $\mathrm{Ad}_{V^*}$ differ.}
 $(\mB(\Hi),\pi,\mathrm{Ad}_{V^*}).$ In Corollary~15 of \cite{WeWe16}, they essentially prove a converse to this statement. 
 
Although Westerbaan and Westerbaan provide a universal property for normal CPU maps $\vf:\mA\stoch\mB$ without choosing an embedding of $\mB$ into bounded operators on some Hilbert space as we have, our results do not require our maps to be normal nor do we require our domain $\mA$ to be a von Neumann algebra. As a result, we provide a universal property for minimal Stinespring dilations for \emph{all} $C^*$-algebras. Furthermore, our universal property highlights the relationship between morphisms of OCP maps and intertwiners of representations. This will have important implications, which will be discussed in Section~\ref{sec:comparetoGNS}. 
\er

There is also a close connection between our work and that of categorical quantum mechanics \cite{AbCo04}, \cite{Se07}. Proposition~\ref{prop:OCPandAnReparedagger} below shows that $\mathbf{OCP}({\mA})$ and $\mathbf{AnRep}({\mA})$ are dagger categories. This fact, together with our Stinespring adjunction, has some interesting consequences such as Theorem~\ref{thm:connectingStinespringrepns}, which states that there exists a unique minimal morphism of anchored representations between \emph{any} two Stinespring representations of the same OCP map.

\bd
\label{defn:daggercategory}
A \define{$^*$-category} (also called a \define{dagger category}) is a category $\mC$ together with a functor $^*:\mC^{\op}\to\mC$ satisfying 
\begin{enumerate}[i.]
\itemsep0pt
\item
$x^*=x$ for all objects $x$ in $\mC$,
\item
$(f^*)^*=f$ for all morphisms $f$ in $\mC.$ 
\end{enumerate}
\ed

Explicitly, functoriality of $^*$ says 
$\id_{x}^*=\id_{x}$ for all objects $x$ in $\mC$
and 
$(g\circ f)^*=f^*\circ g^*$ for all composable pairs of morphisms $f$ and $g$ in $\mC$.

\bn
\label{prop:OCPandAnReparedagger}
Let $\mA$ be a $C^*$-algebra. For each morphism 
$(\mK,\vf)\xrightarrow{T}(\mL,\psi)$ of OCP maps on $\mA,$ set
\be
\label{eq:dualofOCPmorphism}
\Big((\mK,\vf)\xrightarrow{T}(\mL,\psi)\Big)^*:=(\mK,\vf)\xleftarrow{T^*}(\mL,\psi),
\ee
and for each morphism $(\mK,\Hi,\pi,V)\xrightarrow{(T,L)}(\mL,\mI,\rho,W)$ of anchored representations on $\mA,$ set
\be
\label{eq:dualofAnRepmorphism}
\Big((\mK,\Hi,\pi,V)\xrightarrow{(T,L)}(\mL,\mI,\rho,W)\Big)^*:=(\mK,\Hi,\pi,V)\xleftarrow{(T^*,L^*)}(\mL,\mI,\rho,W).
\ee
Then (\ref{eq:dualofOCPmorphism}) is a morphism of OCP maps and (\ref{eq:dualofAnRepmorphism}) is a morphism of anchored representations. Furthermore, $\mathbf{OCP}({\mA})$ and $\mathbf{AnRep}({\mA})$ are $^*$-categories with respect to these assignments. 
\en

\bprf
That (\ref{eq:dualofOCPmorphism}) is a morphism of OCP maps   
has already been proved in the first part of the proof of Proposition~\ref{prop:morphismofstatesdirectsum}. 
An argument similar to it shows that $L^*$ is an intertwiner of representations for (\ref{eq:dualofAnRepmorphism}). 
Furthermore, taking the adjoints of the diagrams in (\ref{eq:morphismanchoredrepns}) shows that (\ref{eq:dualofAnRepmorphism}) is a morphism of anchored representations.
Finally, $\mathbf{OCP}({\mA})$ and $\mathbf{AnRep}({\mA})$ are $^*$-categories due to the definition and properties of the adjoint of a bounded linear map between Hilbert spaces. 
\eprf

Note that $\mathbf{OpSt}({\mA})$ and $\mathbf{PAnRep}({\mA})$ are \emph{not} $^*$-categories with the same $^*$ operation because the adjoint of an isometry need not be an isometry. The adjoint of an isometry is, however, a co-isometry, and hence a partial isometry (cf.~Chapter~15 in Halmos~\cite{Ha82}). 

\bd
\label{defn:partialisometriesorder}
Let $\mathcal{H}$ and $\mathcal{I}$ be two Hilbert spaces. 
The \define{initial space} of a bounded linear map $L:\mathcal{H}\to\mathcal{I}$ is the closed subspace $\ker(L)^{\perp}\subseteq\Hi.$
The map $L:\mathcal{H}\to\mathcal{I}$ is called a \define{co-isometry} iff $L^*$ is an isometry. 
It is called a \define{partial isometry} iff $L$ is an isometry when restricted to $\ker(L)^{\perp}\subseteq\Hi.$ 
\ed

The following lemma includes several properties of partial isometries that will be used in this work. 

\blem
\label{lem:partialisometries}
Let $L:\mI\to\mathcal{J}$ be a bounded linear map between Hilbert spaces. 
\begin{enumerate}[(a)]
\item
If $L$ is a partial isometry, it is an isometry when restricted to its initial space. 
\item
Isometries, co-isometries, projections, and unitary maps are partial isometries.
\item
The composite of a co-isometry followed by its adjoint is a partial isometry. 
\item
\label{item:tensor}
Let $\Hi$ be another Hilbert space (whose dimension is at least $1$). Then $\id_{\Hi}\hat{\otimes} L:\Hi\hat{\otimes}\mI\to\Hi\hat{\otimes}\mathcal{J}$ is a partial isometry if and only if $L$ is a partial isometry. Here, $\hat{\otimes}$ denotes the completed tensor product (cf.~Section~I.2.3 in Dixmier~\cite{Di81}). 
\item
If $L:\mI\to\mathcal{J}$ is a partial isometry, then there exists 
either an isometry or a co-isometry $U:\mI\to\mathcal{J}$ that agrees with $L$ on its initial space. 
\item
If $\mI$ and $\mathcal{J}$ are finite dimensional with the same dimension and $L:\mI\to\mathcal{J}$ is a partial isometry, then there exists a unitary $U:\mI\to\mathcal{J}$ that agrees with $L$ on its initial space. 
\end{enumerate}
\elem

\bprf
Most of these are adequately covered in Halmos~\cite{Ha82}, Halmos--McLaughlin~\cite{HaMc63}, and Hines--Braunstein~\cite{HiBr10}, with the exception of the forward implication in (\ref{item:tensor}). Suppose $\id_{\Hi}\hat{\otimes} L:\Hi\hat{\otimes}\mI\to\Hi\hat{\otimes}\mathcal{J}$ is a partial isometry. Then
\be
\id_{\Hi}\hat{\otimes} L
=(\id_{\Hi}\hat{\otimes} L)(\id_{\Hi}\hat{\otimes} L)^*(\id_{\Hi}\hat{\otimes} L)
=\id_{\Hi}\hat{\otimes} (LL^*L)
\ee
by Corollary~3 in Section~127 of \cite{Ha82} and the property of adjoints with respect to the completed tensor product. Hence, $\id_{\Hi}\hat{\otimes} (L-LL^*L)=0$, which holds if and only if $L=LL^*L$. Thus, $L$ is a partial isometry by this same corollary. 
\eprf

\br
\label{rmk:partialisometriesdonotcompose}
The composite of partial isometries need not be a partial isometry (cf.~Section~9.3 in Hines--Braunstein~\cite{HiBr10}). 
\er

Partial isometries have a natural partial ordering 
as described by Halmos and McLaughlin~\cite{HaMc63}. This partial ordering is slightly generalized below to include intertwiners of representations. 

\bd
\label{defn:partialorder}
Let $\Hi$ and $\mI$ be Hilbert spaces and let $L:\Hi\to\mI$ be a partial isometry. If $M:\Hi\to\mI$ is another partial isometry satisfying $M(x)=L(x)$ for all $x\in\ker(L)^{\perp},$ then $M$ is said to be an \define{extension} of $L$ and the notation $L\leqq M$ will be used.
If $\pi:\mA\to\mB(\Hi)$ and $\rho:\mA\to\mB(\mI)$ are representations of a $C^*$-algebra $\mA$ such that $L$ is an intertwiner from $(\Hi,\pi)$ to $(\mI,\rho)$, then $M$ is said to be an \define{intertwining extension} of $L$ iff $M$ is an intertwiner and $L\leqq M$. In this case, the notation $L\mathrel{\unlhd}M$ will be used when the representations are understood from the context. Whenever the notation $\leqq$ or $\mathrel{\unlhd}$ is used, it will be understood that the operators being compared are partial isometries. 
\ed

It is straightforward to check the following. 

\blem
\label{lem:partialorder}
The relations $\leqq$ and $\mathrel{\unlhd}$ from Definition~\ref{defn:partialorder} define partial orders on the set of partial isometries from $\Hi$ to $\mI$ and $(\Hi,\pi)$ to $(\mI,\rho),$ respectively. 
\elem

The relationship between $\mathrel{\unlhd}$ and morphisms of anchored representations will be described in Lemma~\ref{lem:intertwiningorders}. But first, we provide a useful result that compares \emph{any} two Stinespring representations of the same OCP map. 

\bt
\label{thm:connectingStinespringrepns}
Given any two Stinespring representations $(\mK,\Hi,\pi,V)$ and $(\mK,\mI,\rho,W)$ of an OCP map $(\mK,\vf)$ on a $C^*$-algebra $\mA,$ there exists a unique partial isometry $L:\Hi\to\mI$ such that
\begin{enumerate}[i.]
\item
$(\mK,\Hi,\pi,V)\xrightarrow{(\id_{\mK},L)}(\mK,\mI,\rho,W)$ is a morphism of anchored representations, 
\item
$L$ restricted to $\overline{\pi(\mA)V(\mK)}\subseteq\Hi$ is a unitary intertwiner onto $\overline{\rho(\mA)W(\mK)}\subseteq\mI,$ and
\item
for any other bounded linear map $M:\Hi\to\mI$ satisfying the first condition, $M$ agrees with $L$ on the intersection of their initial spaces (and hence also satisfies the second condition). 
\end{enumerate}
\et

The third condition is a uniqueness condition guaranteeing there exists a unique minimal partial isometry (minimal in the sense of the partial order on partial isometries) $L$ for which $(\id_{\mK},L)$ is a morphism of anchored representations. This is justified in Lemma~\ref{lem:intertwiningorders}. 

\bprf[Proof of Theorem~\ref{thm:connectingStinespringrepns}]
By Corollary~\ref{cor:Stinespringinitial}, we have a unique pair of morphisms of anchored representations consisting of isometries of the form
\be
\xy0;/r.25pc/:
(0,7.5)*+{(\mK,\Hi_{\vf},\pi_{\vf},V_{\vf})}="S";
(-20,-7.5)*+{(\mK,\Hi,\pi,V)}="H";
(20,-7.5)*+{(\mK,\mI,\rho,W)}="I";
{\ar"S";"H"_{(\id_{\mK},m_{\pi,V})}};
{\ar"S";"I"^{(\id_{\mK},m_{\rho,W})}};
\endxy
.
\ee
By Proposition~\ref{prop:OCPandAnReparedagger}, 
the composite $(\mK,\Hi,\pi,V)\xrightarrow{(\id_{\mK},m_{\rho,W}m_{\pi,V}^*)}(\mK,\mI,\rho,W)$ defines a morphism of anchored representations that restricts to $\id_{\mK}$. Since $m_{\pi,V}$ and $m_{\rho,W}$ are isometries, $L:=m_{\rho,W}m_{\pi,V}^*$ is a partial isometry by Lemma~\ref{lem:partialisometries}. Therefore, $\Hi$ and $\mI$ can be decomposed as
\be
\label{eq:invariantsubspaces}
\Hi=m_{\pi,V}(\Hi_{\vf})\oplus m_{\pi,V}(\Hi_{\vf})^{\perp}
\quad\text{ and }\quad
\mI=m_{\rho,W}(\Hi_{\vf})\oplus m_{\rho,W}(\Hi_{\vf})^{\perp}, 
\ee
where each term in the direct sums is an invariant subspace under the appropriate actions of $\pi$ or $\rho$. To see this for the orthogonal complement, take the dual of (\ref{eq:intertwiner}). 
By construction, 
$L$ restricts to a unitary intertwiner from $m_{\pi,V}(\Hi_{\vf})$ to $m_{\rho,W}(\Hi_{\vf}).$ 
But by the definition of $m$ induced from (\ref{eq:mpiV}), 
\be
\label{eq:mpiVHvfpiVK}
m_{\pi,V}(\Hi_{\vf})=\overline{\pi(A)V(\mK)}
\quad\text{ and }\quad
m_{\rho,W}(\Hi_{\vf})=\overline{\rho(A)W(\mK)},
\ee 
which proves the first two claims.  

Now, let $M:\Hi\to\mI$ be a bounded linear map satisfying the first condition and let $x=\sum_{i}\pi(a_{i})V(v_{i})\in\pi(\mA)V(\mK).$ Then
\be
M(x)=\sum_{i}M\big(\pi(a_{i})V(v_{i})\big)=\sum_{i}\rho(a_{i})M\big(V(v_{i})\big)=\sum_{i}\rho(a_{i})W(v_{i})=L(x)
\ee
since $M$ is linear, since $M\pi(a)=\rho(a)M$ for all $a\in\mA,$ and because $MV=W$. By continuity, $M$ agrees with $L$ on $\overline{\pi(\mA)V(\mK)},$ the initial space of $L.$ 
\eprf

\blem
\label{lem:intertwiningorders}
In terms of the same notation as in Theorem~\ref{thm:connectingStinespringrepns}, set 
\be
\mathcal{P}(L):=\left\{\Hi\xrightarrow{M}\mI\;:\;L\leqq M\mbox{ and }\;\left((\mK,\Hi,\pi,V)\!\!\xrightarrow{(\id_{\mK},M)}\!\!(\mK,\mI,\rho,W)\right)\in\mathbf{AnRep}(\mA)_{1}\right\},
\ee
where the bounded linear maps $M$ are all understood to be partial isometries. 
Then 
\be
\label{eq:mypartialorder}
\mathcal{P}(L)=\left\{\Hi\xrightarrow{M}\mI\;:\;L\mathrel{\unlhd} M\right\}.
\ee
\elem

\bprf
The containment $\supseteq$ in (\ref{eq:mypartialorder}) holds by the diagram~(\ref{eq:intertwiner}) from Definition~\ref{defn:anchoredrepn}. 
To see the reverse containment,
let $M:\Hi\to\mI$ be a partial isometry such that $L\mathrel{\unlhd}M.$ Then commutativity of the diagrams 
\be
\label{eq:Lppartial}
\xy0;/r.25pc/:
(-10,0)*+{\mK}="K";
(10,7.5)*+{\Hi}="H";
(10,-7.5)*+{\mI}="I";
{\ar"K";"H"^{V}};
{\ar"K";"I"_{W}};
{\ar"H";"I"^{M}};
\endxy
\quad\text{ and }\quad
\xy0;/r.25pc/:
(10,0)*+{\mK}="K";
(-10,7.5)*+{\Hi}="H";
(-10,-7.5)*+{\mI}="I";
{\ar"H";"K"^{V^*}};
{\ar"I";"K"_{W^*}};
{\ar"H";"I"_{M}};
\endxy
\ee
follows from (\ref{eq:mpiVHvfpiVK}), the decomposition (\ref{eq:invariantsubspaces}), and the fact that $(\mK,\Hi,\pi,V)\xrightarrow{(\id_{\mK},L)}(\mK,\mI,\rho,W)$ is a morphism of anchored representations. Finally, the intertwining condition (\ref{eq:intertwiner}) is included in the definition of $L\mathrel{\unlhd}M.$ Hence, $\subseteq$ holds as well.
\eprf

\br
Although maximal elements for the partial order $\leqq$ on partial isometries are either isometries or co-isometries (cf.~\cite{HaMc63}), this is not true in general of maximal elements for the partial order $\mathrel{\unlhd}$ on intertwining partial isometries.
This point will be addressed in Remark~\ref{rmk:droppingequivalencedirectsum}. 
\er

\section{Examples and applications}
\label{sec:comparetoGNS}

The purification postulate has been used by Chiribella, D'Ariano, and Perinotti to classify finite-dimensional quantum theories among all operational probabilistic theories (OPTs) \cite{CDP10}, \cite{Ch14}. We do not need to review OPTs here, but will instead provide a definition of a purification of a process, our formulation of the purification postulate, and the standard finite-dimensional purification postulate of \cite{CDP10}. Our version of the purification postulate isolates some key assumptions made by \cite{CDP10} that are implicit from the tensor network (diagrammatic) perspective. We prove our purification postulate using our Stinespring adjunction and show how it reduces to the standard finite-dimensional one. 
We will first use our Stinespring adjunction to reproduce a Gelfand--Naimark--Segal (GNS) adjunction for states \cite{PaGNS}.

\bx
\label{ex:GNS}
Let $\mK:=\C.$ By Example~\ref{ex:states}, an operator state $\w:\mA\stoch\mB(\C)\cong\C$ is a state. Applying $\mathbf{Stine}_{\mA}$ to $\w$ provides an anchored representation $(\C,\Hi_{\w},\pi_{\w},V_{\w}).$ In this case, $\mA\otimes\C\cong\mA$ so that $\mathcal{N}_{\w}\cong\{a\in\mA\;:\;\w(a^*a)=0\}$ under this isomorphism. This agrees with the null-space from the usual GNS construction. Hence, the completion $\Hi_{\w}:=\overline{\mA/\mathcal{N}_{\w}}$ and the representation $\pi_{\w}:\mA\to\mB(\Hi_{\w})$ agree with the usual GNS Hilbert space and representation. By Example~\ref{ex:pRep}, $V_{\w}:\C\to\Hi_{\w}$ produces a unit vector $\W_{\w}:=V_{\w}(1)\in\Hi_{\w}$ and $V_{\w}^*=\<\W_{\w},\;\cdot\;\>:\Hi_{\w}\to\C$. Hence, $\mathrm{Ad}_{V^*_{\w}}:\mB(\Hi_{\w})\stoch\mB(\C)$ can be identified with its evaluation at $1$ and is equivalently described by the state $\<\W_{\w},\;\cdot\;\W_{\w}\>:\mB(\Hi_{\w})\stoch\C.$ Applying $\mathbf{rest}_{\mA}$ entails $\w=\<\W_{\w},\pi(\;\cdot\;)\W_{\w}\>$, so that $\w$ has been represented by a pure state. If $(\C,\Hi,\pi,V)$ is another anchored representation of $\w,$ set $\W:=V(1)\in\Hi.$ By Corollary~\ref{cor:Stinespringinitial}, there is a unique morphism of anchored representations 
of the form
$(\C,\Hi_{\w},\pi_{\w},V_{\w})\xrightarrow{(\id_{\C},m_{\pi,V})}(\C,\Hi,\pi,V)$, where $m_{\pi,V}$ is given by 
\be
\begin{split}
\Hi_{\w}&\xrightarrow{m_{\pi,V}}\Hi\\
[a]_{\w}&\xmapsto{\;\;\quad\;\;}\pi(a)\W, 
\end{split}
\ee
which agrees with the modification from Construction 5.21 in \cite{PaGNS}. 
The universal property of Stinespring's adjunction thus reproduces the minimality of the GNS construction in the sense that it reproduces the smallest cyclic representation of $\mA$ on which $\w$ can be realized as a pure state. 
\ex

\br
\label{rmk:GNS}
The GNS adjunction in \cite{PaGNS} was defined in terms of a category of pointed representations instead of anchored representations. This remark further explains this slight (and only technical) distinction. It is not needed for the present work but is included for completeness. There are natural transformations
$\S:\pRep\Rightarrow\mathbf{PAnRep}$
and $\U:\St\Rightarrow\mathbf{OpSt}$ defined as follows. 
First, for every $C^*$-algebra $\mA,$
the category $\pRep(\mA)$ has objects 
$(\Hi,\pi,\W)$ with $\Hi$ a Hilbert space, $\pi:\mA\to\mB(\Hi)$ 
a $C^*$-algebra representation, and $\W\in\Hi$ a unit vector. 
A morphism $(\Hi,\pi,\W)\xrightarrow{L}(\mI,\rho,\Xi)$ is an
isometric intertwiner of representations such that $L(\W)=\Xi.$ 
The category $\St(\mA),$ on the other hand, is the discrete category of
states on $\mA,$ i.e.~the objects are positive unital linear maps 
$\w:\mA\stoch\C$ and there are no non-identity morphisms.  
To every $C^*$-algebra $\mA,$ set
\be
\begin{split}
\pRep(\mA)&\xrightarrow{\S_{\mA}}\mathbf{PAnRep}(\mA)\\
(\Hi,\pi,\W)&\xmapsto{\quad}(\C,\Hi,\pi,V_{\W})\\
\Big((\Hi,\pi,\W)\xrightarrow{L}(\mI,\rho,\Xi)\Big)&\xmapsto{\quad}\Big((\C,\Hi,\pi,V_{\W})\xrightarrow{(\id_{\C},L)}(\C,\mI,\rho,V_{\Xi})\Big).
\end{split}
\ee
Here, $V_{\W}:\C\to\Hi$ is the map that sends $\l\in\C$ to $\l\W.$ 
It is not difficult to show that $\S_{\mA}$ is a functor.
Also, set 
\be
\begin{split}
\St(\mA)&\xrightarrow{\U_{\mA}}\mathbf{OpSt}(\mA)\\
\w&\xmapsto{\quad}(\C,\tilde{\w}),
\end{split}
\ee
where $\tilde{\w}:\mA\stoch\mB(\C)$ is defined by 
$\mA\ni a\mapsto\tilde{\w}(a)=\w(a)\;\cdot\;,$
i.e.~multiplication by $\w(a)$ on the Hilbert space $\C.$ 
Since $\St(\mA)$ has only identity morphisms, this specifies the functor
$\U_{\mA}.$ 
Examples~\ref{ex:states}~and~\ref{ex:pRep} show that the functors
$\S_{\mA}$ and $\U_{\mA}$ are faithful but not full. 
The only reason these functors are not full is that
the categories $\mathbf{OpSt}(\mA)$ and $\mathbf{PAnRep}(\mA)$
contain more data in their morphisms. However, 
the only added information for a morphism of \emph{states} and 
pointed representations
is a phase factor, which is a symmetry that can safely be ignored 
in the discussion of the GNS construction.

Given a $^*$-homomorphism $f:\mA'\to\mA,$ 
the equalities 
\be
\mathbf{PAnRep}(f)\circ\S_{\mA}=\S_{\mA'}\circ\pRep(f)
\quad\text{and}\quad
\mathbf{OpSt}(f)\circ\U_{\mA}=\U_{\mA'}\circ\St(f)
\ee
also hold. Furthermore, in the diagram%
\footnote{The dashed arrows are used to depict the three-dimensional nature of the diagram in (\ref{eq:GNSStinediagram}). They do not have an alternative mathematical meaning.}
\be
\label{eq:GNSStinediagram}
\xy0;/r.25pc/:
(-25,10)*+{\CAlg^{\op}}="1";
(22,5)*+{\tCAT}="2";
(-25,-10)*+{\CAlg^{\op}}="3";
(22,-15)*+{\tCAT}="4";
{\ar@/^1.75pc/"1";"2"^{\mathbf{OpSt}}};
{\ar@/_1.75pc/"1";"2"|-{\mathbf{AnRep}}};
{\ar@{=>}(3.5,1.5);(8.5,11.5)_{\rest}};
{\ar@{=>}(-5.5,13.5);(-10.5,3.5)_{\mathbf{Stine}}};
{\ar@{-->}@/^1.75pc/"3";"4"|-{{\transparent{0.75}\St}}};
{\ar@/_1.75pc/"3";"4"_{\pRep}};
{\ar@{==>}(3.5,-19.0);(8.5,-9.0)_{\rest}};
{\ar@{==>}(-5.5,-7.0);(-10.5,-17.0)_{\pGNS}};
{\ar@{=}"1";"3"};
{\ar@{=}"2";"4"};
{\ar@{=>}(-3,-18);(-3,-2)_{\S}};
{\ar@{==>}(1,-4);(1,13)|-(0.265){\large\phantom{X}}^{\U}};
\endxy
,
\ee
although the equality 
\be
\begin{matrix}\rest\\\U\end{matrix}=\begin{matrix}\S\\\rest\end{matrix}
\ee
holds, there is only an invertible modification
\be
\label{eq:modificationStineGNS}
\begin{matrix}\U\\\mathbf{Stine}\end{matrix}
\Rrightarrow
\begin{matrix}\pGNS\\\S\end{matrix}
\ee
since the two composites are not \emph{exactly} equal but are canonically isomorphic. 
Indeed, for a fixed $C^*$-algebra $\mA,$ the resulting natural
isomorphism 
\be
\xy0;/r.25pc/:
(-20,7.5)*+{\mathbf{PAnRep}(\mA)}="TL";
(20,7.5)*+{\mathbf{OpSt}(\mA)}="TR";
(-20,-7.5)*+{\pRep(\mA)}="BL";
(20,-7.5)*+{\St(\mA)}="BR";
{\ar"BR";"TR"_{\U_{\mA}}};
{\ar"TR";"TL"_{\mathbf{Stine}_{\mA}}};
{\ar"BR";"BL"^{\pGNS_{\mA}}};
{\ar"BL";"TL"^{\S_{\mA}}};
{\ar@{=>}"BL";"TR"_{\cong}};
\endxy
\ee
is defined by its evaluation on a state $\w:\mA\stoch\C$ by the morphism
\be
\big(\C,\Hi_{\w},\pi_{\w},V_{\W_{\w}}\big)
\xrightarrow{(\id_{\C},L_{\w})}
\big(\C,\Hi_{\tilde{\w}},\pi_{\tilde{\w}},V_{\tilde{\w}}\big),
\ee
where $L_{\w}:\Hi_{\w}\to\Hi_{\tilde{\w}}$ is defined as the unique extension of
\be
\mA/\mN_{\w}\ni[a]\mapsto[a\otimes1]\in(\mA\otimes\C)/\mN_{\tilde{\w}}.
\ee
Similar calculations to the above show that this map is bounded and
extends to a unitary intertwiner. Hence, $(\id_{\C},L_{\w})$ defines
an isomorphism in the category $\mathbf{PAnRep}(\mA).$
The appropriate diagram also commutes when one considers
a $^*$-homomorphism $f:\mA'\to\mA.$ 

This tells us our Stinespring adjunction reduces to the GNS
adjunction by restricting to the images of $\S$ and $\U.$ 
For example, if $\w:\mA\stoch\C$ is a state, one can construct a pointed representation of it via GNS and view that pointed representation as a preserving anchored representation. Similarly, one can view $\w$ as an operator state and apply Stinespring's construction to obtain another preserving anchored representation. The invertible modification in (\ref{eq:modificationStineGNS}) says that these two preserving anchored representations are canonically isomorphic. Working out the details of Stinespring's adjunction in this present work has highlighted the importance of not restricting morphisms to be isometries so that our universal property is more robust and improves on our previous GNS adjunction in several respects. 
\er

\bd
Let $\mA$ be a $C^*$-algebra. A \define{purification} of a state $\w:\mA\stoch\C$ consists of a Hilbert space $\Hi,$ a $^*$-homomorphism $\pi:\mA\to\mB(\Hi),$ and a unit vector $\W\in\Hi$ such that 
\be
\<\W,\pi(a)\W\>_{\Hi}=\w(a)\qquad\forall\;a\in\mA. 
\ee
Hence, a purification of $\w$ will be written as a triple $(\Hi,\pi,\W).$ 
In other words, a purification of a state is Stinespring representation of that state.
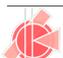
\ed

\bd
Let $\mA$ be a $C^*$-algebra and $\mK$ a Hilbert space. A \define{purification} of an OCP map $\vf:\mA\stoch\mB(\mK)$ consists of a Hilbert space $\Hi,$ a $^*$-homomorphism $\pi:\mA\to\mB(\Hi),$ and a bounded linear map $V:\mK\to\Hi$ such that 
\be
V^*\pi(a)V=\vf(a)\qquad\forall\;a\in\mA. 
\ee
In other words, a purification of an OCP map is a Stinespring representation of that OCP map.
\ed

The following version of the purification postulate might seem unfamiliar, but we show that it is equivalent to the usual purification postulate (from Section~VII.A. and~VII.B. in~\cite{CDP10}) when the algebras are finite-dimensional matrix algebras. This equivalence is worked out in detail in Lemma~\ref{lem:finitedimensionalStinespringdilation} and Theorem~\ref{thm:ChiribellaPurification}. 

\begin{postulate}
[The Purification Postulate for processes]
\label{thm:purificationpostulateprocess}
Let $\mA$ be a $C^*$-algebra and $\mK$ a Hilbert space. Every OCP map $\vf:\mA\stoch\mB(\mK)$ has a purification $(\mK,\Hi,\pi,V).$ Furthermore, given any other purification $(\mK,\mI,\rho,W)$ for which $(\Hi,\pi)$ is unitarily equivalent to $(\mI,\rho),$ there exists a unitary intertwiner $\Hi\xrightarrow{U}\mI$ such that $UV=W.$ 
\end{postulate}

Stinespring's theorem guarantees the existence of purifications. The existence of the unitary intertwiner in these postulates is referred to as the \emph{essential uniqueness of purifications}. 
The purification postulate for processes implies the one for states by setting $\mK=\C.$ We will prove the purification postulate for processes on finite-dimensional matrix algebras in Theorem~\ref{thm:ChiribellaPurification} and finite-dimensional $C^*$-algebras in Corollary~\ref{cor:fdcalgpurification} after a few lemmas.
Our version of the purification postulate is formulated without using traces or tensor products since these may be absent or ambiguous for general $C^*$-algebras. This is partially achieved by using completely positive \emph{unital} maps instead of the more common completely positive \emph{trace-preserving} maps. In the finite-dimensional setting, these are equivalent and correspond to the Heisenberg and Schr\"odinger pictures, respectively. Completely positive unital maps are used to map observables to observables while their duals (adjoints with respect to the Hilbert--Schmidt inner product), completely positive trace-preserving maps, are used to map density matrices to density matrices. However, the category of 
$C^*$-algebras does not have duals and therefore does not have symmetric purifications as defined in \cite{SSC18}.
Since states supersede density matrices in the infinite-dimensional setting and the notion of a trace is not always available, it is sometimes more convenient to work within the Heisenberg picture. 

\blem
\label{lem:finitedimensionalStinespringdilation}
Fix positive integers $k$ and $n$. 
Set $\mK:=\C^{k}$, let $\vf:\mM_{n}(\C)\stoch\mB(\C^{k})\cong\mathcal{M}_{k}(\C)$ be an OCP map, and let $(\mK,\Hi,\pi,V)$ be a finite-dimensional Stinespring representation of $(\mK,\vf).$ Then there exists a $p\in\N$ and a unitary map $R:\Hi\to\C^{n}\otimes\C^{p}$ such that
\be
R\pi(A)R^*=A\otimes\mathds{1}_{p}
\quad\text{ and }\quad
\vf(A)=(RV)^*(A\otimes\mathds{1}_{p})(RV)
\qquad\forall\;A\in\mM_{n}(\C). 
\ee
In other words, there exists an isomorphism $(\mK,\Hi,\pi,V)\xrightarrow{(\id_{\mK},R)}(\mK,\C^{n}\otimes\C^{p},i,RV)$ of anchored representations, where $R$ is unitary and where $i$ is defined by
\be
\label{eq:itakesAtoAtensorid}
\begin{split}
i:\mM_{n}(\C)&\to\mM_{n}(\C)\otimes\mM_{p}(\C)\cong\mB(\C^{n}\otimes\C^{p})\\
A\;\;\;\;&\xmapsto{\quad}\;\;A\otimes\mathds{1}_{p}.
\end{split}
\ee
\elem

\bprf
This follows from a general fact about unital $^*$-homomorphisms between finite-dimensional matrix algebras (cf.~Section~1.1.2 in Fillmore~\cite{Fi96}) and our definitions.
\eprf

Lemma~\ref{lem:finitedimensionalStinespringdilation} provides a characterization of the form of Stinespring dilations of OCP maps between matrix algebras. This allows a visual calculus to be implemented via circuit diagrams \cite{CDP10}, \cite{SSC18}. We set
\be
\vcenter{\hbox{\rotatebox{90}{$\xrightarrow{\qquad\;\text{time}\qquad\;}$}}}
\quad
\vcenter{\hbox{\rotatebox{-90}{$\xrightarrow{\quad\text{composition}\quad}$}}}
\qquad\quad
\vcenter{\hbox{
\begin{tikzpicture}
\node[draw,align=center,inner sep=7pt] (vf) at (0,0){$\vf$};
\draw (0,1.25) to node[left]{$\mA$} (vf);
\draw (vf) to node[left]{$\mB$} (0,-1.25);
\end{tikzpicture}
}}
:=
\xy0;/r.25pc/:
(0,7.5)*+{\mA}="A";
(0,-7.5)*+{\mB}="B";
{\ar@{~>}"A";"B"^{\vf}};
\endxy
\quad,\quad
\vcenter{\hbox{
\begin{tikzpicture}[circuit ee IEC]
\node[ground,align=center,rotate=90] (g) at (0,0) {};
\draw (g) to node[left]{$\mA$} (0,-1.25);
\end{tikzpicture}
}}
:=
\xy0;/r.25pc/:
(0,7.5)*+{\C}="C";
(0,-7.5)*+{\mA}="A";
{\ar@{~>}"C";"A"^{!}};
\endxy
\quad,
\ee
\be
\vcenter{\hbox{
\begin{tikzpicture}
\node[draw,align=center,inner sep=7pt] (vf) at (0,0){$\vf$};
\draw (0,1.25) to node[left]{$\mA$} (vf);
\draw (vf) to node[left]{$\mB$} (0,-1.25);
\end{tikzpicture}
}}
\vcenter{\hbox{
\begin{tikzpicture}
\node[draw,align=center,inner sep=5pt] (vf) at (0,0){\;$\vf'$};
\draw (0,1.25) to node[left]{$\mA'$} (vf);
\draw (vf) to node[left]{$\mB'$} (0,-1.25);
\end{tikzpicture}
}}
:=
\vcenter{\hbox{
\begin{tikzpicture}
\node[draw,align=center,inner sep=5pt] (vf) at (0,0){$\vf\otimes\vf'$};
\draw (0,1.25) to node[left]{$\mA\otimes\mA'$} (vf);
\draw (vf) to node[left]{$\mB\otimes\mB'$} (0,-1.25);
\end{tikzpicture}
}}
\quad,\quad
\text{ and }
\quad
\vcenter{\hbox{
\begin{tikzpicture}
\node[draw,align=center,inner sep=7pt] (vf) at (0,0){$\vf$};
\node[draw,align=center,inner sep=7pt] (ps) at (0,-2){$\psi$};
\draw (0,1.25) to node[left]{$\mA$} (vf);
\draw (vf) to node[left]{$\mB$} (ps);
\draw (ps) to node[left]{$\mC$} (0,-3.25);
\end{tikzpicture}
}}
:=
\vcenter{\hbox{
\begin{tikzpicture}
\node[draw,align=center,inner sep=7pt] (vf) at (0,0){$\psi\circ\vf$};
\draw (0,1.25) to node[left]{$\mA$} (vf);
\draw (vf) to node[left]{$\mC$} (0,-1.25);
\end{tikzpicture}
}}
.
\ee
The map $\C\xstoch{!}\mA$ is the unique unital linear map from $\C$ to $\mA.$ 
Note that the algebras here are finite dimensional 
so the tensor product is the standard one. By choosing this convention, the direction of time is up. The direction of time is consistent with \cite{SSC18}, while the direction of composition is opposite, because the dual maps on density matrices are used in \cite{SSC18}. 
Using the shorthand $\mathcal{M}_{n}:=\mathcal{M}_{n}(\C)$ (and similarly for other dimensions), the relationship between our dilation (from Lemma~\ref{lem:finitedimensionalStinespringdilation}) and the one of \cite[Definition~45]{CDP10} is then given by
\be
\vcenter{\hbox{
\begin{tikzpicture}
\node[draw,align=center,inner sep=7pt] (vf) at (0,0){$\vf$};
\draw (0,1.25) to node[left]{$\mathcal{M}_{n}$} (vf);
\draw (vf) to node[left]{$\mathcal{M}_{k}$} (0,-1.25);
\end{tikzpicture}
}}
=
\vcenter{\hbox{
\begin{tikzpicture}
\node[draw,align=center,inner sep=7pt] (pi) at (0,0){$\pi$};
\node[draw,align=center,inner sep=7pt] (V) at (0,-2){$\mathrm{Ad}_{V^*}$};
\draw (0,1.25) to node[left]{$\mathcal{M}_{n}$} (pi);
\draw (pi) to node[left]{$\mB(\Hi)$} (V);
\draw (V) to node[left]{$\mathcal{M}_{k}$} (0,-3.25);
\end{tikzpicture}
}}
=
\vcenter{\hbox{
\begin{tikzpicture}
\node[draw,align=center,inner sep=7pt] (pi) at (0,0){$\pi$};
\node[draw,align=center,inner sep=7pt] (R) at (0,-2){$\mathrm{Ad}_{R}$};
\node[draw,align=center,inner sep=7pt] (Ri) at (0,-4){$\mathrm{Ad}_{R^*}$};
\node[draw,align=center,inner sep=7pt] (V) at (0,-6){$\mathrm{Ad}_{V^*}$};
\draw (0,1.25) to node[left]{$\mathcal{M}_{n}$} (pi);
\draw (pi) to node[left]{$\mB(\Hi)$} (R);
\draw ([xshift=-0.5cm]R.south) -- node[left]{$\mathcal{M}_{n}$} ([xshift=-0.5cm]Ri.north);
\draw ([xshift=0.5cm]R.south) to node[left]{$\mathcal{M}_{p}$} ([xshift=0.5cm]Ri.north);
\draw (Ri) to node[left]{$\mB(\Hi)$} (V);
\draw (V) to node[left]{$\mathcal{M}_{k}$} (0,-7.25);
\end{tikzpicture}
}}
=
\vcenter{\hbox{
\begin{tikzpicture}[circuit ee IEC]
\node[ground,align=center,rotate=90] (g) at (0.5,0) {};
\node[draw,align=center,inner sep=7pt] (RV) at (0,-2){$\mathrm{Ad}_{(RV)^*}$};
\draw ([xshift=-0.5cm]0,1.25) to node[left]{$\mathcal{M}_{n}$} ([xshift=-0.5cm]RV.north);
\draw(g) to node[left]{$\mathcal{M}_{p}$} ([xshift=0.5cm]RV.north);
\draw (RV) to node[left]{$\mathcal{M}_{k}$} (0,-3.25);
\end{tikzpicture}
}}
\ee
due to Lemma~\ref{lem:finitedimensionalStinespringdilation}. 
To prove the purification postulate for processes, we need to recall a few standard facts about the commutant (cf.~Chapters~1~and~2 in Dixmier~\cite{Di81} and Chapter~4 in Fillmore~\cite{Fi96}). 

\bd
\label{defn:commutant}
Let $S\subseteq\mA$ be a subset of a $C^*$-algebra $\mA.$ 
The \define{commutant} of $S$ inside $\mA$ is
the unital algebra
\be
S':=\{a\in\mA\;:\;as=sa\;\forall\;s\in S\}.
\ee
Since the commutant depends on the embedding algebra, 
$S'$ will often be written as $S'\subseteq\mA.$%
\footnote{Writing $S'\subseteq\mA$ will also avoid confusion with
the notation used for different $C^*$-algebras $\mA,\mA',\mA''$ earlier.}
\ed

\br
\label{rmk:commutantclosed}
If a subset $S\subseteq\mA$ is $^*$-closed 
(meaning $a\in S$ implies $a^*\in S$), then 
$S'$ is a unital $^*$-algebra. 
In fact, $S'\subseteq\mA$ is a $C^*$-subalgebra of $\mA$
since $S'=\bigcap_{s\in S}\{s\}'$ is the intersection of the kernels of the commutators $[s,\;\cdot\;]:\mA\to\mA,$ which are all closed (cf.~Chapter~2 in Topping~\cite{To71}).
\er

\bx
\label{ex:commutants}
If $\Hi$ is a finite-dimensional Hilbert space, then $\mB(\Hi)'\subseteq\mB(\Hi)$ is $\C\id_{\Hi},$ all scalar multiples of the identity in $\Hi.$ 
If $\mK$ is another finite-dimensional Hilbert space, then the commutant $\big(\mB(\Hi)\otimes\id_{\mK}\big)'\subseteq\mB(\Hi)\otimes\mB(\mK)$ is $\id_{\Hi}\otimes\mB(\mK),$ all bounded operators of the form ${\id_{\Hi}\otimes B}$ with $B\in\mB(\mK).$ If $\pi:\mathcal{M}_{n}(\C)\to\mB(\Hi)$ and $R:\Hi\to\C^{n}\otimes\C^{p}$ are as in Lemma~\ref{lem:finitedimensionalStinespringdilation}, then $\pi(\mA)'\subseteq\mB(\Hi)$ is $R^*\big(\mathds{1}_{n}\otimes\mathcal{M}_{p}(\C)\big)R,$ operators of the form $R^*(\mathds{1}_{n}\otimes B)R$ with $B\in\mathcal{M}_{p}(\C).$ Given a positive integer $t$ and non-negative integers $n_{1},\dots,n_{t},c_{1},\dots,c_{t},$ with $\sum_{j}n_{j},\sum_{j}c_{j}>0$, set $m=\sum_{j=1}^{t}n_{j}c_{j}$ and let 
\be
\bigboxplus_{j=1}^{t}\Big(\mathcal{M}_{n_{j}}(\C)\otimes\mathcal{M}_{c_{j}}(\C)\Big)
\ee 
be the $C^*$-subalgebra of $\mathcal{M}_{m}(\C)$ consisting of all linear combinations of  block diagonal matrices of the form 
\be
\bigboxplus_{j=1}^{t}(A_{j}\otimes B_{j}):=
\begin{bmatrix}
A_{1}\otimes B_{1}&&0\\
&\ddots&\\
0&&A_{t}\otimes B_{t}
\end{bmatrix}
\ee
(if any of the integers $n_{1},\dots,n_{t},c_{1},\dots,c_{t},$ are zero, terms corresponding to them are excluded from the above matrix).
Then, the commutant of $\bigboxplus_{j=1}^{t}\big(\mathcal{M}_{n_{j}}(\C)\otimes\mathds{1}_{c_{j}}\big)$ inside $\mathcal{M}_{m}(\C)$ is $\bigboxplus_{j=1}^{t}\big(\mathds{1}_{n_{j}}\otimes\mathcal{M}_{c_{j}}(\C)\big).$
The notation $\bigboxplus$ is used for `internal' direct sum to distinguish it from the `external' direct sum $\bigoplus$. 
\ex

\bt
\label{thm:ChiribellaPurification}
The purification postulate for OCP maps  
(Postulate~\ref{thm:purificationpostulateprocess}) holds whenever all algebras are finite-dimensional matrix algebras and all Hilbert spaces are finite dimensional. 
\et

\bprf
Purifications exist by Stinespring's theorem (the associated anchored representation is finite dimensional by minimality). For the essential uniqueness, fix positive integers $k$ and $n$. 
Set $\mK:=\C^{k}$, let $\vf:\mM_{n}(\C)\stoch\mB(\C^{k})\cong\mathcal{M}_{k}(\C)$ be an OCP map, and let $(\mK,\Hi,\pi,V)$ and $(\mK,\mI,\rho,W)$ be two finite-dimensional Stinespring representations of $(\mK,\vf)$ with unitarily equivalent representations $\pi:\mathcal{M}_{n}(\C)\to\mB(\Hi)$ and $\rho:\mathcal{M}_{n}(\C)\to\mB(\mI)$. 
By Theorem~\ref{thm:connectingStinespringrepns}, there exists a partial isometry $L:\Hi\to\mI$ such that 
$(\mK,\Hi,\pi,V)\xrightarrow{(\id_{\mK},L)}(\mK,\mI,\rho,W)$ is a morphism of anchored representations. 
In what follows, we will prove there exists a unitary $U$ such that $L\mathrel{\unlhd}U.$ By Lemma~\ref{lem:finitedimensionalStinespringdilation}, there exist integers $p$ and $q$ together with unitary maps $R:\Hi\to\C^{n}\otimes\C^{p}$ and $S:\mI\to\C^{n}\otimes\C^{q}$ such that
\be
\begin{gathered}
R\pi(A)R^*=A\otimes\mathds{1}_{p}, \qquad S\rho(A)S^*=A\otimes\mathds{1}_{q},
\quad\text{ and }\;\\
(RV)^*(A\otimes\mathds{1}_{p})(RV)=\vf(A)=(SW)^*(A\otimes\mathds{1}_{q})(SW)
\end{gathered}
\ee
for all $A\in\mM_{n}(\C).$ Since $\pi$ and $\rho$ are unitarily equivalent, $p=q$. 
These facts imply $(\mK,\C^{n}\otimes\C^{p},i,SW)\xrightarrow{(\id_{\mK},RL^*S^*)}(\mK,\C^{n}\otimes\C^{p},i,RV)$ is a morphism of anchored representations, where $RL^*S^*:\C^{n}\otimes\C^{p}\to\C^{n}\otimes\C^{p}$ is a partial isometry and $i$ is the map from (\ref{eq:itakesAtoAtensorid}). 
By Example~\ref{ex:commutants} and Lemma~\ref{lem:partialisometries}, there exists a unique partial isometry $P:\C^{p}\to\C^{p}$ such that $RL^*S^*=\mathds{1}_{n}\otimes P$. By the last fact in Lemma~\ref{lem:partialisometries}, there exists a unitary $M:\C^{p}\to\C^{p}$ such that $P\leqq M$ (cf.~Definition~\ref{defn:partialorder}). Hence $\mathds{1}_{n}\otimes P\leqq\mathds{1}_{n}\otimes M$. By Example~\ref{ex:commutants}, $\mathds{1}_{n}\otimes M$ is in $i\big(\mathcal{M}_{n}(\C)\big)'\subseteq\mathcal{M}_{n}(\C)\otimes\mathcal{M}_{p}(\C).$ Therefore, $\mathds{1}_{n}\otimes P\mathrel{\unlhd}\mathds{1}_{n}\otimes M$ with respect to the representation $i$. Hence, 
$(\mK,\C^{n}\otimes\C^{p},i,SW)\xrightarrow{(\id_{\mK},\mathds{1}_{n}\otimes M)}(\mK,\C^{n}\otimes\C^{p},i,RV)$ is a morphism of anchored representations by Lemma~\ref{lem:intertwiningorders}. Setting $U:=S^*(\mathds{1}_{n}\otimes M^*)R,$ we obtain the required unitary for essential uniqueness of purifications since it provides an isomorphism 
$(\mK,\Hi,\pi,V)\xrightarrow{(\id_{\mK},U)}(\mK,\mI,\rho,W)$
of anchored representations.
\eprf

The situation in the previous proof is summarized by the commutative diagram 
\be
\label{eq:unitaryforpurification}
\xy0;/r.25pc/:
(20,0)*+{\mM_{n}(\C)}="A";
(-20,0)*+{\mM_{k}(\C)}="B";
(0,7.5)*+{\mB(\Hi)}="H";
(0,-7.5)*+{\mB(\mI)}="I";
(0,20)*+{\mM_{n}(\C)\otimes\mM_{p}(\C)}="top";
(0,-20)*+{\mM_{n}(\C)\otimes\mM_{p}(\C)}="bot";
{\ar@{~>}"H";"B"_{\qquad\mathrm{Ad}_{V^*}}|(0.66){\hole}};
{\ar@{~>}"I";"B"^{\qquad\;\;\mathrm{Ad}_{W^*}}|(0.66){\hole}};
{\ar"A";"H"_{\pi}};
{\ar"A";"I"^{\rho}};
{\ar"I";"H"|-{\mathrm{Ad}_{U^*}}};
{\ar"H";"top"|-(0.4){\mathrm{Ad}_{R}}};
{\ar"bot";"I"|-(0.4){\mathrm{Ad}_{S^*}}};
{\ar@/_1.5pc/"A";"top"_(0.35){i}};
{\ar@/^1.5pc/"A";"bot"^(0.35){i}};
{\ar@{~>}@/_1.5pc/"top";"B"_(0.65){\mathrm{Ad}_{(RV)^*}}|(0.929){\hole}};
{\ar@{~>}@/^1.5pc/"bot";"B"^(0.65){\mathrm{Ad}_{(SW)^*}}|(0.928){\hole}};
\endxy
.
\ee
Because $R,U,$ and $S$ are unitary, commutativity of this diagram is equivalent to the conditions required of morphisms of anchored representations. 
In terms of circuit diagrams, we have
\be
\vcenter{\hbox{
\begin{tikzpicture}[circuit ee IEC]
\node[ground,align=center,rotate=90] (g) at (0.5,0) {};
\node[draw,align=center,inner sep=7pt] (RV) at (0,-2){$\mathrm{Ad}_{(RV)^*}$};
\draw ([xshift=-0.5cm]0,1.25) to node[left]{$\mathcal{M}_{n}$} ([xshift=-0.5cm]RV.north);
\draw(g) to node[left]{$\mathcal{M}_{p}$} ([xshift=0.5cm]RV.north);
\draw (RV) to node[left]{$\mathcal{M}_{k}$} (0,-3.25);
\end{tikzpicture}
}}
=
\vcenter{\hbox{
\begin{tikzpicture}
\node[draw,align=center,inner sep=7pt] (vf) at (0,0){$\vf$};
\draw (0,1.25) to node[left]{$\mathcal{M}_{n}$} (vf);
\draw (vf) to node[left]{$\mathcal{M}_{k}$} (0,-1.25);
\end{tikzpicture}
}}
=
\vcenter{\hbox{
\begin{tikzpicture}[circuit ee IEC]
\node[ground,align=center,rotate=90] (g) at (0.5,0) {};
\node[draw,align=center,inner sep=7pt] (RV) at (0,-2){$\mathrm{Ad}_{(SW)^*}$};
\draw ([xshift=-0.5cm]0,1.25) to node[left]{$\mathcal{M}_{n}$} ([xshift=-0.5cm]RV.north);
\draw(g) to node[left]{$\mathcal{M}_{p}$} ([xshift=0.5cm]RV.north);
\draw (RV) to node[left]{$\mathcal{M}_{k}$} (0,-3.25);
\end{tikzpicture}
}}
\ee
for our two purifications of $\vf$. 
The middle part of the diagram (\ref{eq:unitaryforpurification}) relates the representations and the Stinespring dilations from our point of view. Namely, we obtain a unitary $U:\mathcal{H}\to\mathcal{I}$ such that 
\be
\vcenter{\hbox{
\begin{tikzpicture}
\node[draw,align=center,inner sep=7pt] (vf) at (0,0){$\rho$};
\node[draw,align=center,inner sep=7pt] (ps) at (0,-2){$\mathrm{Ad}_{U^*}$};
\draw (0,1.25) to node[left]{$\mathcal{M}_{n}$} (vf);
\draw (vf) to node[left]{$\mB(\mI)$} (ps);
\draw (ps) to node[left]{$\mB(\Hi)$} (0,-3.25);
\end{tikzpicture}
}}
=
\vcenter{\hbox{
\begin{tikzpicture}
\node[draw,align=center,inner sep=7pt] (vf) at (0,0){$\pi$};
\draw (0,1.25) to node[left]{$\mathcal{M}_{n}$} (vf);
\draw (vf) to node[left]{$\mB(\Hi)$} (0,-1.25);
\end{tikzpicture}
}}
\quad\text{ and }\quad
\vcenter{\hbox{
\begin{tikzpicture}
\node[draw,align=center,inner sep=7pt] (vf) at (0,0){$\mathrm{Ad}_{W^*}$};
\draw (0,1.25) to node[left]{$\mB(\mI)$} (vf);
\draw (vf) to node[left]{$\mathcal{M}_{k}$} (0,-1.25);
\end{tikzpicture}
}}
=
\vcenter{\hbox{
\begin{tikzpicture}
\node[draw,align=center,inner sep=7pt] (vf) at (0,0){$\mathrm{Ad}_{U^*}$};
\node[draw,align=center,inner sep=7pt] (ps) at (0,-2){$\mathrm{Ad}_{V^*}$};
\draw (0,1.25) to node[left]{$\mB(\mI)$} (vf);
\draw (vf) to node[left]{$\mB(\Hi)$} (ps);
\draw (ps) to node[left]{$\mathcal{M}_{k}$} (0,-3.25);
\end{tikzpicture}
}}
\quad.
\ee
Combining the three parts gives
\be
\vcenter{\hbox{
\begin{tikzpicture}
\node[draw,align=center,inner sep=7pt] (S) at (0,0){$\mathrm{Ad}_{S^*}$};
\node[draw,align=center,inner sep=7pt] (U) at (0,-2){$\mathrm{Ad}_{U^*}$};
\node[draw,align=center,inner sep=7pt] (R) at (0,-4){$\mathrm{Ad}_{R}$};
\draw ([xshift=-0.5cm]0,1.25) to node[left]{$\mathcal{M}_{n}$} ([xshift=-0.5cm]S.north);
\draw ([xshift=0.5cm]0,1.25) to node[left]{$\mathcal{M}_{p}$} ([xshift=0.5cm]S.north);
\draw (S) to node[left]{$\mB(\Hi)$} (U);
\draw (U) to node[left]{$\mB(\mI)$} (R);
\draw ([xshift=-0.5cm]R.south) to node[left]{$\mathcal{M}_{n}$} ([xshift=-0.5cm]0,-5.25);
\draw ([xshift=0.5cm]R.south) to node[left]{$\mathcal{M}_{p}$} ([xshift=0.5cm]0,-5.25);
\end{tikzpicture}
}}
=
\vcenter{\hbox{
\begin{tikzpicture}
\node[draw,align=center,inner sep=7pt] (O) at (0.5,2){$\mathrm{Ad}_{M}$};
\draw ([xshift=-0.5cm]0,3.25) to node[left]{$\mathcal{M}_{n}$} ([xshift=-0.5cm]vf.north);
\draw (O) to node[left]{$\mathcal{M}_{p}$} (0.5,0.5);
\draw (0.5,3.25) to node[left]{$\mathcal{M}_{p}$} (O);
\end{tikzpicture}
}}
\implies
\vcenter{\hbox{
\begin{tikzpicture}
\node[draw,align=center,inner sep=7pt] (vf) at (0,0){$\mathrm{Ad}_{(SW)^*}$};
\draw ([xshift=-0.5cm]0,1.25) to node[left]{$\mathcal{M}_{n}$} ([xshift=-0.5cm]vf.north);
\draw ([xshift=0.5cm]0,1.25) to node[left]{$\mathcal{M}_{p}$} ([xshift=0.5cm]vf.north);
\draw (vf) to node[left]{$\mathcal{M}_{k}$} (0,-1.25);
\end{tikzpicture}
}}
=
\vcenter{\hbox{
\begin{tikzpicture}
\node[draw,align=center,inner sep=7pt] (vf) at (0,0){$\mathrm{Ad}_{(RV)^*}$};
\node[draw,align=center,inner sep=7pt] (O) at (0.5,2){$\mathrm{Ad}_{M}$};
\draw ([xshift=-0.5cm]0,3.25) to node[left]{$\mathcal{M}_{n}$} ([xshift=-0.5cm]vf.north);
\draw (O) to node[left]{$\mathcal{M}_{p}$} ([xshift=0.5cm]vf.north);
\draw (0.5,3.25) to node[left]{$\mathcal{M}_{p}$} (O);
\draw (vf) to node[left]{$\mathcal{M}_{k}$} (0,-1.25);
\end{tikzpicture}
}}
\quad,
\ee
which agrees precisely with Postulate~1 in Section~VII.A. in~\cite{CDP10} (after setting $k=1$). 

\br
The unitary $U$ in the essential uniqueness of Postulate~\ref{thm:purificationpostulateprocess} is not necessarily unique if it exists. Indeed, there could be many unitary extensions $\mathds{1}_{n}\otimes M$ of the partial isometry $\mathds{1}_{n}\otimes P$ from the proof of Theorem~\ref{thm:ChiribellaPurification}. 
\er

\br
If one drops the assumptions that the finite-dimensional representations $(\Hi,\pi)$ and $(\mI,\rho)$ are unitarily equivalent in Postulate~\ref{thm:purificationpostulateprocess} and Theorem~\ref{thm:ChiribellaPurification}, then one obtains a slight generalization of the purification postulate. In this more general case, there exists an intertwining map $U:\Hi\to\mI$ that satisfies $UV=W$ and $U$ is either an isometry or a co-isometry. The proof of this is similar to the proof of Theorem~\ref{thm:ChiribellaPurification} except that $p$ need not equal $q$. 
The only partial isometry $O$ that satisfies $(A\otimes\mathds{1}_{p})O=O(A\otimes\mathds{1}_{q})$ for all $A\in\mathcal{M}_{n}(\C)$ must be of the form $\mathds{1}_{n}\otimes P$ for some partial isometry $P:\C^{q}\to\C^{p}$. Hence, 
one obtains a unique partial isometry $P:\C^{q}\to\C^{p}$ such that $RL^* S^*=\mathds{1}_{n}\otimes P$. 
All intertwining extensions of this partial isometry must therefore also be of this form. Therefore, there exists either an isometry or a co-isometry 
$M:\C^{q}\to\C^{p}$ satisfying $\mathds{1}_{n}\otimes P\mathrel{\unlhd}\mathds{1}_{n}\otimes M$ by Lemma~\ref{lem:partialisometries}. Then $U:=S^*(\mathds{1}_{n}\otimes M)^*R$ is the required partial isometry. Note that if $\dim\Hi=\dim\mI$ is assumed, then the representations are automatically unitarily equivalent since maximal partial isometries between finite-dimensional Hilbert spaces of equal dimension are unitary. 
\er

Our formulation of the purification postulate is also valid for arbitrary finite-dimensional $C^*$-algebras. 

\bc
\label{cor:fdcalgpurification}
The purification postulate for OCP maps 
(Postulate~\ref{thm:purificationpostulateprocess}) holds whenever all algebras are finite-dimensional $C^*$-algebras and all Hilbert spaces are finite dimensional. 
\ec

\bprf
Existence follows from Stinespring's theorem as before. What follows is a proof of the essential uniqueness of purifications. 
By the discussions preceding this, it suffices to consider the case 
\be
\mK:=\C^{k}
\quad\text{ and }\quad
\mA:=\bigoplus_{j=1}^{t}\mathcal{M}_{n_{j}}(\C),
\ee
where $t\in\N$ and $n_{j}\in\N$ for all $j\in\{1,\dots,t\}.$ Given an OCP map $\vf:\mA\stoch\mB(\mK),$ it also suffices to consider two Stinespring representations of the form $(\mK,\C^{m},\pi,V)$ and $(\mK,\C^{m},\pi,W)$, where $m\in\N$ and the representation $\pi:\mA\to\mathcal{M}_{m}(\C)$ is of the form
\be
\pi\left(\bigoplus_{j=1}^{t}A_{j}\right)=
\bigboxplus_{j=1}^{t}(A_{j}\otimes\mathds{1}_{c_{j}})
,
\ee 
where the $c_{1},\dots,c_{t}$ are non-negative integers such that $m=\sum_{j=1}^{t}n_{j}c_{j}.$ By similar arguments to those implemented in the proof of Theorem~\ref{thm:ChiribellaPurification}, there exist partial isometries (for the non-zero $c_{j}$) $L_{j}:\C^{c_{j}}\to\C^{c_{j}}$ such that $\left(\id_{\mK},L:=\bigboxplus_{j=1}^{t}(\mathds{1}_{n_{j}}\otimes L_{j})\right)$ is a morphism of anchored representations from $(\mK,\C^{m},\pi,V)$ to $(\mK,\C^{m},\pi,W)$. These can be extended to unitaries $U_{j}:\C^{c_{j}}\to\C^{c_{j}}$ by finite dimensionality. Hence, $\left(\id_{\mK},U:=\bigboxplus_{j=1}^{t}(\mathds{1}_{n_{j}}\otimes U_{j})\right)$ is an isomorphism of anchored representations.
\eprf

\br
\label{rmk:droppingequivalencedirectsum}
If one drops the assumptions that the finite-dimensional representations $(\Hi,\pi)$ and $(\mI,\rho)$ are unitarily equivalent in Postulate~\ref{thm:purificationpostulateprocess} and Corollary~\ref{cor:fdcalgpurification}, then $U=\bigboxplus_{j=1}^{t}(\mathds{1}_{n_{j}}\otimes U_{j})$ from the proof of Corollary~\ref{cor:fdcalgpurification} is replaced by an internal block sum of  partial isometries that have been extended to isometries or co-isometries. In particular, even though $U$ is maximal with respect to $\mathrel{\unlhd}$, it need not be an isometry nor a co-isometry. In more detail, let $(\mK,\C^{m},\pi,V)$ and $(\mK,\C^{n},\rho,W)$ be the two Stinespring representations, where 
\be
\begin{split}
\pi\left(\bigoplus_{j=1}^{t}A_{j}\right)=\bigboxplus_{j=1}^{t}(A_{j}\otimes\mathds{1}_{c_{j}})
,\quad&\quad
\rho\left(\bigoplus_{j=1}^{t}A_{j}\right)=\bigboxplus_{j=1}^{t}(A_{j}\otimes\mathds{1}_{d_{j}}),
\\
m=\sum_{j=1}^{t}n_{j}c_{j},
\quad\qquad\;\;&\text{and}\qquad\qquad
n=\sum_{j=1}^{t}n_{j}d_{j}.
\end{split}
\ee
One can show $L$ (from Theorem~\ref{thm:connectingStinespringrepns}) must be of the form $\bigboxplus_{j=1}^{t}(\mathds{1}_{n_{j}}\otimes P_{j})$. 
In fact, all intertwining extensions of $L$ must also be of this form. 
By extending such intertwining partial isometries, one obtains an isometry or a co-isometry $U_{j}:\C^{c_{j}}\to\C^{d_{j}}$ if $c_{j}\le d_{j}$ or $c_{j}\ge d_{j}$, respectively. 
Hence, if $c_{j}-d_{j}$ changes sign as $j$ varies, then $U=\bigboxplus_{j=1}^{t}(\mathds{1}_{n_{j}}\otimes U_{j})$ is neither an isometry nor a co-isometry. 
\er

\acknowledgments{
The majority of this work was completed when the author was an Assistant Research Professor at the University of Connecticut. 
The author thanks Benjamin Russo for
many helpful discussions about 
completely positive maps. 
The author thanks Edward Effros~\cite{Ef78} for inspiration. 
Finally, and most importantly, the author thanks two anonymous reviewers and the editorial board of \emph{Compositionality}, who have 
provided numerous constructive comments and suggestions that have
greatly improved this work. }


\appendix{

\section{Index of notation}
\label{index}
\begin{longtable}{c|c|c|c}
\hline
Notation & Name/description & Location & Page \\
\hline
$\N$& natural numbers (excludes $0$)& Not~\ref{not:intro}& \pageref{not:intro}\\
\hline
\begin{tabular}{c}$\Hi,\mI,\mJ,$\\$\mK,\mL,\mM$\end{tabular}& Hilbert spaces& Not~\ref{not:intro}& \pageref{not:intro}\\
\hline
$\mB(\Hi)$&bounded operators on $\Hi$&Not~\ref{not:intro}& \pageref{not:intro}\\
\hline
$\mA,\mA',\mA''$&(unital) $C^*$-algebra&Not~\ref{not:intro}& \pageref{not:intro}\\
\hline
$\<\;\cdot\;,\;\cdot\;\>_{\mK}$&inner product on $\mK$&Not~\ref{not:intro}&\pageref{not:intro}\\
\hline
$\mM_{n}(\mA)$&$n\times n$ matrices with coeffs in $\mA$&Not~\ref{not:intro}& \pageref{not:intro}\\
\hline
$\vf,\psi,\chi$&\begin{tabular}{c}positive or completely positive maps\\operator states when unital\end{tabular}&\begin{tabular}{cc}Def'n~\ref{defn:cpmaps}\\Def'n~\ref{defn:opstate}\end{tabular}&\begin{tabular}{cc}\pageref{defn:cpmaps}\\\pageref{defn:opstate}\end{tabular}\\
\hline
$\vf_{n}$&the $n$-ampliation of $\vf$&Def'n~\ref{defn:cpmaps}&\pageref{defn:cpmaps}\\
\hline
\begin{tabular}{c}PU, CP,\\CPU\end{tabular}&\begin{tabular}{c}positive unital, completely positive,\\completely positive unital\end{tabular}&Def'n~\ref{defn:cpmaps}&\pageref{defn:cpmaps}\\
\hline
$\mathrm{Ad}_{T}$&the adjoint action map for $T$&Ex~\ref{ex:cpmaps}&\pageref{ex:cpmaps}\\
\hline
$(\mK,\vf)$&\begin{tabular}{c}operator-valued CP (OCP)\\map or operator state $\vf:\mA\stoch\mB(\mK)$\end{tabular}&Def'n~\ref{defn:opstate}&\pageref{defn:opstate}\\
\hline
$\tr$&un-normalized trace&Ex~\ref{ex:tracial}&\pageref{ex:tracial}\\
\hline
$\t$&tracial map&Ex~\ref{ex:tracial}&\pageref{ex:tracial}\\
\hline
\begin{tabular}{c}$\mathbf{OCP}$\\$\mathbf{OpSt}$\end{tabular}&\begin{tabular}{c}OCP maps functor\\operator states functor\end{tabular}&\begin{tabular}{c}Lem~\ref{lem:opstA}, \\\ref{lem:opstf}, \ref{lem:opstfunctor}\end{tabular}&\begin{tabular}{c}\pageref{lem:opstA},\\ \pageref{lem:opstf}, \pageref{lem:opstfunctor}\end{tabular}\\
\hline
$(\mK,\Hi,\pi,V)$&anchored representation on a $C^*$-algebra&Def'n~\ref{defn:anchoredrepn}&\pageref{defn:anchoredrepn}\\
\hline
\begin{tabular}{c}$\mathbf{AnRep}$\\$\mathbf{PAnRep}$\end{tabular}&\begin{tabular}{c}anchored representation functor\\preserving anchored representation functor\end{tabular}&\begin{tabular}{c}Lem~\ref{lem:anrepA}, \\\ref{lem:anrepf}, \ref{lem:anrepfunctor}\end{tabular}&\begin{tabular}{c}\pageref{lem:anrepA},\\\pageref{lem:anrepf}, \pageref{lem:anrepfunctor}\end{tabular}\\
\hline
$\mathbf{rest}$&restriction natural transformation&Prop~\ref{prop:rest}&\pageref{prop:rest}\\
\hline
$s_{\vf,\vec{v}}$&a certain linear functional $\mM_{n}(\mA)\stoch\C$&Lem~\ref{lem:CPtoP}&\pageref{lem:CPtoP}\\
\hline
$\mathbf{Stine}$&Stinespring oplax-natural transformation&Thm~\ref{thm:stinespring}&\pageref{thm:stinespring}\\
\hline
$\otimes$&algebraic tensor product of vector spaces&\begin{tabular}{c}Item~\ref{item:1} in\\Thm~\ref{thm:stinespring}\end{tabular}&\pageref{item:1}\\
\hline
$\llangle\;\cdot\;,\;\cdot\;\rrangle_{\vf}$&\begin{tabular}{c}sesquilinear form on $\mA\otimes\mK$\\from $\vf:\mA\stoch\mB(\mK)$\end{tabular}&Eqn~(\ref{eq:sesquilinear})&\pageref{eq:sesquilinear}\\
\hline
$\mN_{\vf}$&null-space associated to $\vf$&Eqn~(\ref{eq:StinespringIdeal})&\pageref{eq:StinespringIdeal}\\
\hline
$\pi_{\vf}$&\begin{tabular}{c}Stinespring representation of a $C^*$-algebra\\from an OCP map $(\mK,\vf)$\end{tabular}&\begin{tabular}{c}Eqn~(\ref{eq:stinespringreppi})\\Eqn~(\ref{eq:stinespringrepaction})\end{tabular}&\begin{tabular}{c}\pageref{eq:stinespringreppi},\\\pageref{eq:stinespringrepaction}\end{tabular}\\
\hline
$[\xi]_{\vf}$ or $[\xi]$&elements of $(\mA\otimes\mK)/\mathcal{N}_{\vf}$&Eqn~(\ref{eq:stinespringrepaction})&\pageref{eq:stinespringrepaction}\\
\hline
$\<\;\cdot\;,\;\cdot\;\>_{\vf}$&\begin{tabular}{c}induced inner product on\\$(\mA\otimes\mK)/\mN_{\vf}$ and $\Hi_{\vf}$ from $\llangle\;\cdot\;,\;\cdot\;\rrangle_{\vf}$\end{tabular}&Eqn~(\ref{eq:inducedinnerproduct})&\pageref{eq:inducedinnerproduct}\\
\hline
$\Hi_{\vf}$&\begin{tabular}{c}Stinespring Hilbert space\\$\ov{(\mA\otimes\mK)/\mN_{\vf}}$ associated to $\vf$\end{tabular}&Eqn~(\ref{eq:stinespringspace})&\pageref{eq:stinespringspace}\\
\hline
$V_{\vf}$&Stinespring isometry $V_{\vf}:\mK\to\Hi_{\vf}$&Eqn~(\ref{eq:stinespringisometry})&\pageref{eq:stinespringisometry}\\
\hline
$(\mK,\Hi_{\vf},\pi_{\vf},V_{\vf})$&\begin{tabular}{c}$\mathbf{Stine}_{\mA}(\mK,\vf),$ Stinespring anchored\\rep'n from OCP map $(\mK,\vf)$\end{tabular}&Eqn~(\ref{eq:stineonobjects})&\pageref{eq:stineonobjects}\\
\hline
$(T,L_{T})$&\begin{tabular}{c}$\mathbf{Stine}_{\mA}(T),$ Stinespring morphism\\from OCP map morphism $T$\end{tabular}&Eqn~(\ref{eq:stineonmorphisms})&\pageref{eq:stineonmorphisms}\\
\hline
$L_{f}$&\begin{tabular}{c}$\mathbf{Stine}_{f}(\mK,\vf),$ Stinespring isometry\\for a $^*$-homomorphism $f:\mA'\to\mA$\end{tabular}&Eqn~(\ref{eq:Lf})&\pageref{eq:Lf}\\
\hline
$m_{\pi,V}$&\begin{tabular}{c}Stinespring isometry associated to an\\ anchored representation $(\mK,\Hi,\pi,V)$\end{tabular}&Eqn~(\ref{eq:mpiV})&\pageref{eq:mpiV}\\
\hline
$m_{\mA}$&Stinespring natural transformation on $\mA$&Eqn~(\ref{eq:mA})&\pageref{eq:mA}\\
\hline
$m$&Stinespring modification&Eqn~(\ref{eq:modification})&\pageref{eq:modification}\\
\hline
$\hat{\otimes}$&\begin{tabular}{c}completed tensor product\\for Hilbert spaces\end{tabular}&Lem~\ref{lem:partialisometries}&\pageref{lem:partialisometries}\\
\hline
$L\leqq M$&extension of partial isometry&Def'n~\ref{defn:partialorder}&\pageref{defn:partialorder}\\
\hline
$L\mathrel{\unlhd}M$&intertwining extension of partial isometry&Def'n~\ref{defn:partialorder}&\pageref{defn:partialorder}\\
\hline
$S'\subseteq\mA$&the commutant of $S$ in $\mA$&Def'n~\ref{defn:commutant}&\pageref{defn:commutant}\\
\hline
$\bigboxplus$&
`internal' direct sum&Ex~\ref{ex:commutants}&\pageref{ex:commutants}\\
\hline
\end{longtable}

\section{2-categorical preliminaries}
\label{appendix}

We briefly recall the definitions of oplax-natural transformations and modifications.
In addition, we include the universal property associated with adjunctions
because it is used in explaining the Stinespring adjunction more concretely. 
For details on 2-categories and their pasting diagrams, we refer the reader to B{\'e}nabou's original work
\cite{Be} as well as Kelly and Street's review \cite{KeSt74}.
For a more introductory take emphasizing string diagrams, see \cite{Pa18}.
For other details on oplax-natural transformations and modifications, we refer the reader to Section 7.5 in Borceux~\cite{Bo94}.

\bd
\label{defn:semipseudonaturaltransformation}
Let $\mathcal{C}$ and $\mathcal{D}$ be two
(strict)
2-categories and let $F, G : \mathcal{C} \to \mathcal{D}$ 
be two (strict) functors. An
\define{oplax-natural transformation}
$\rho$ from $F$ to $G$, written as 
$\rho:F\Rightarrow G,$ consists of
\begin{enumerate}[i.]
\item
a function $\rho : C_{0} \to D_{1}$ assigning a 
1-morphism in $\mD$ to each object $x$ in $\mC$ in the following manner 
\be
\xy0;/r.15pc/:
(-30,0)*+{x}="1";
(30,10)*+{F(x)}="2";
(30,-10)*+{G(x)}="3";
{\ar^{\rho(x)} "2";"3"};
{\ar@{|->}^{\rho} "1"+(20,0);(10,0)};
\endxy
\ee

\item
and a function $\rho : C_{1} \to D_{2}$ assigning 
a 2-morphism in $\mD$ to each 1-morphism $y \xleftarrow{\a} x$ in $\mC$ in the following manner 
\be
\label{eq:oplaxrho}
\xy0;/r.15pc/:
(-30,0)*+{x}="1x";
(-60,0)*+{y}="1y";
{\ar_{\a} "1x";"1y"};
(30,10)*+{F(y)}="2y";
(30,-10)*+{G(y)}="3y";
(60,10)*+{F(x)}="2x";
(60,-10)*+{G(x)}="3x";
{\ar^{\rho(x)} "2x";"3x"};
{\ar_{\rho(y)} "2y";"3y"};
{\ar_{F(\a)} "2x";"2y"};
{\ar^{G(\a)} "3x";"3y"};
{\ar@{=>}^{\rho(\a)} "2y";"3x"};
{\ar@{|->}^{\rho} "1x"+(20,0);(10,0)};
\endxy
.
\ee
\end{enumerate}

These data must satisfy the following conditions:

\begin{enumerate}[(a)]
\item
For every object $x$ in $\mC,$ 
\be
\label{eq:oplaxidentity}
\rho(\id_{x})=\id_{\rho(x)}.
\ee

\item
For every pair $(z\xleftarrow{\a}y,y\xleftarrow{\b}x)$ 
of composable 1-morphisms in $\mathcal{C},$ 
\be
\label{eq:oplaxnaturalitycomposition}
\xy0;/r.20pc/:
(-30,10)*+{F(z)}="2z";
(-30,-10)*+{G(z)}="3z";
(0,10)*+{F(y)}="2y";
(0,-10)*+{G(y)}="3y";
(30,10)*+{F(x)}="2x";
(30,-10)*+{G(x)}="3x";
{\ar|-{\rho(x)} "2x";"3x"};
{\ar|-{\rho(y)} "2y";"3y"};
{\ar|-{\rho(z)} "2z";"3z"};
{\ar_{F(\b)} "2x";"2y"};
{\ar^{G(\b)} "3x";"3y"};
{\ar_{F(\a)} "2y";"2z"};
{\ar^{G(\a)} "3y";"3z"};
{\ar@{=>}^{\rho(\b)} "2y";"3x"};
{\ar@{=>}^{\rho(\a)} "2z";"3y"};
\endxy
\qquad
=
\qquad
\xy0;/r.20pc/:
(-15,10)*+{F(z)}="2z";
(-15,-10)*+{G(z)}="3z";
(15,10)*+{F(x)}="2x";
(15,-10)*+{G(x)}="3x";
{\ar|-{\rho(x)} "2x";"3x"};
{\ar|-{\rho(z)} "2z";"3z"};
{\ar_{F(\a\b)} "2x";"2z"};
{\ar^{G(\a\b)} "3x";"3z"};
{\ar@{=>}^{\rho(\a\b)} "2z";"3x"};
\endxy
.
\ee
\item
For every 2-morphism 
\be
\label{eq:2morphismSigma}
\xymatrix{
y &&
\ar@/_2pc/[ll]_{\a}="8"
\ar@/^2pc/[ll]^{\g}="9"
\ar@{}"8";"9"|(.2){\,}="10"
\ar@{}"8";"9"|(.8){\,}="11"
\ar@{=>}"10";"11"^{\S}
x 
}
,
\ee 
the equality 
\be
\label{eq:oplaxnaturalSigma}
\xy0;/r.25pc/:
(-12.5,7.5)*+{F(y)}="Fy";
(-12.5,-7.5)*+{G(y)}="Gy";
(12.5,7.5)*+{F(x)}="Fx";
(12.5,-7.5)*+{G(x)}="Gx";
{\ar@/_1.25pc/"Fx";"Fy"_{F(\a)}};
{\ar@/^1.25pc/"Fx";"Fy"|-{F(\g)}};
{\ar@/^1.25pc/"Gx";"Gy"^{G(\g)}};
{\ar"Fy";"Gy"_{\rho(y)}};
{\ar"Fx";"Gx"^{\rho(x)}};
{\ar@{=>}@/_0.75pc/"Fy"+(3,-5);"Gx"_(0.4){\rho(\g)}};
{\ar@{=>}(0,12);(0,4)|-{F(\S)}};
\endxy
\qquad
=
\qquad
\xy0;/r.25pc/:
(-12.5,7.5)*+{F(y)}="Fy";
(-12.5,-7.5)*+{G(y)}="Gy";
(12.5,7.5)*+{F(x)}="Fx";
(12.5,-7.5)*+{G(x)}="Gx";
{\ar@/_1.25pc/"Fx";"Fy"_{F(\a)}};
{\ar@/_1.25pc/"Gx";"Gy"|-{G(\a)}};
{\ar@/^1.25pc/"Gx";"Gy"^{G(\g)}};
{\ar"Fy";"Gy"_{\rho(y)}};
{\ar"Fx";"Gx"^{\rho(x)}};
{\ar@{=>}@/^0.75pc/"Fy";"Gx"^{\rho(\a)}};
{\ar@{=>}(0,-4);(0,-12)|-{G(\S)}};
\endxy
.
\ee
holds. 
\end{enumerate}
\ed

\br
We use the prefix ``oplax'' because for a \define{lax-natural transformation} (cf.~Definition~7.5.2 in Borceux~\cite{Bo94}), the source and targets of the 2-morphism in (\ref{eq:oplaxrho}) are switched. Note that Equations~(\ref{eq:oplaxnaturalitycomposition})~and~(\ref{eq:oplaxnaturalSigma}) must be appropriately modified for a lax-natural transformation. For a \define{pseudo-natural transformation}, the 2-morphisms in (\ref{eq:oplaxrho}) are (vertically) invertible (cf.~Definition~\ref{defn:verticalcompnattrans}). The vertical inverse of $\rho(\a)$ will be written as $\overline{\rho(\a)}.$ By the uniqueness of inverses, if a pseudo-natural transformation is oplax, then its vertical inverse is lax. This fact is used in Proposition~\ref{prop:adjunctionsinfunctorcategories}. If $\rho(\a)$ is the identity for all 1-morphisms $\a,$ then $\rho$ is called a \define{natural transformation}. 
\er

The definition of a modification does not
change if one uses oplax-natural
transformations instead of
pseudo-natural transformations. 

\bd
\label{defn:modification}
Let $\mathcal{C}$ and $\mathcal{D}$ be two 
2-categories, $F,G:\mathcal{C}\to\mathcal{D}$ 
be two 2-functors, and $\rho,\s:F\Rightarrow G$ 
be two oplax-natural transformations. A 
\define{modification} $m$ 
from $\s$ to $\rho$, written as 
$m:\s\Rrightarrow\rho$ and drawn as 
\be
\xy0;/r.15pc/:
(0,0)*+{\mathcal{C}}="2";
(-40,0)*+{\mathcal{D}}="3";
{\ar@/_2.25pc/_{F} "2";"3"};
{\ar@/^2.25pc/^{G} "2";"3"};
{\ar@2{->}(-30,7.5);(-30,-7.5)_{\r}};
{\ar@{=>}(-10,7.5);(-10,-7.5)^{\s}};
{\ar@3{->}(-13,0);(-27,0)_{m}};
\endxy
,
\ee
consists of a function $m:C_{0}\to D_{2}$ 
assigning a 2-morphism in $\mD$ to each object $x$ in $\mC$ in the 
following manner 
\be
\xy0;/r.15pc/:
(-30,0)*+{x}="1";
(50,15)*+{F(x)}="2";
(50,-15)*+{G(x)}="3";
{\ar@/_1.75pc/_{\rho(x)}"2";"3"};
{\ar@/^1.75pc/^{\s(x)}"2";"3"};
{\ar@2{->}(60,0);(40,0)_{m(x)}};
{\ar@{|->}^{m}"1"+(20,0);(10,0)};
\endxy
.
\ee
This assignment must satisfy the condition that 
for every 1-morphism $y\xleftarrow{\a}x,$ 
\be
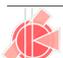
\label{eq:coherencemodification}
\xy0;/r.25pc/:
(-12.5,10)*+{F(y)}="Fy";
(-12.5,-10)*+{G(y)}="Gy";
(12.5,10)*+{F(x)}="Fx";
(12.5,-10)*+{G(x)}="Gx";
{\ar"Fx";"Fy"_{F(\a)}};
{\ar"Gx";"Gy"^{G(\a)}};
{\ar@/_1.75pc/_{\s(y)}"Fy";"Gy"};
{\ar@/^1.75pc/|-{\rho(y)}"Fy";"Gy"};
{\ar@/^1.75pc/|-{\rho(x)}"Fx";"Gx"};
{\ar@2{->}(-18.5,0);(-8.5,0)^{m(y)}};
{\ar@2{->}@/^0.75pc/"Fy";"Gx"^(0.6){\rho(\a)}};
\endxy
\qquad
=
\qquad
\xy0;/r.25pc/:
(-12.5,10)*+{F(y)}="Fy";
(-12.5,-10)*+{G(y)}="Gy";
(12.5,10)*+{F(x)}="Fx";
(12.5,-10)*+{G(x)}="Gx";
{\ar"Fx";"Fy"_{F(\a)}};
{\ar"Gx";"Gy"^{G(\a)}};
{\ar@/_1.75pc/|-{\s(y)}"Fy";"Gy"};
{\ar@/_1.75pc/|-{\s(x)}"Fx";"Gx"};
{\ar@/^1.75pc/^{\rho(x)}"Fx";"Gx"};
{\ar@2{->}(8.5,0);(18.5,0)^{m(x)}};
{\ar@2{->}@/_0.75pc/"Fy";"Gx"_(0.35){\s(\a)}};
\endxy
.
\ee
\ed

\br
If one has a modification between lax-natural transformations, the diagram in 
(\ref{eq:coherencemodification}) is modified appropriately. 
This will be used in Proposition~\ref{prop:adjunctionsinfunctorcategories}. 
\er

The composition of oplax-natural transformations
and modifications are not changed as a result
of these alterations to the usual definitions.

\bd
\label{defn:verticalcompnattrans}
The \define{vertical composite} of oplax-natural
transformations is denoted using vertical concatenation as in
\be
\xy0;/r.15pc/:
(25,0)*+{\mathcal{C}}="2";
(-25,0)*+{\mathcal{D}}="3";
{\ar@/_2.25pc/_{F} "2";"3"};
{\ar@/^2.25pc/^{H} "2";"3"};
{\ar@2{->}(0,13.5);(0,-13.5)_{\begin{smallmatrix}\r\\\s\end{smallmatrix}}};
\endxy
\quad:=\quad
\xy0;/r.15pc/:
(25,0)*+{\mathcal{C}}="2";
(-25,0)*+{\mathcal{D}}="3";
{\ar@/_2.25pc/_{F} "2";"3"};
{\ar|-(0.5){G} "2";"3"};
{\ar@/^2.25pc/^{H} "2";"3"};
{\ar@2{->}(0,13.5);(0,2.5)_{\r}};
{\ar@{=>}(0,-2.5);(0,-13.5)_{\s}};
\endxy
\ee
and is defined by the assignments
\be
\xy0;/r.15pc/:
(-30,0)*+{x}="1";
(30,20)*+{F(x)}="2";
(30,0)*+{G(x)}="3";
(30,-20)*+{H(x)}="4";
{\ar^{\rho(x)} "2";"3"};
{\ar^{\s(x)} "3";"4"};
{\ar@{|->}^{\begin{smallmatrix}\r\\\s\end{smallmatrix}} "1"+(20,0);(10,0)};
\endxy
\ee
for each object $x$ in $\mC$ and 
\be
\xy0;/r.15pc/:
(-30,0)*+{x}="1x";
(-60,0)*+{y}="1y";
{\ar_{\a} "1x";"1y"};
(30,20)*+{F(y)}="2y";
(30,0)*+{G(y)}="3y";
(30,-20)*+{H(y)}="4y";
(70,20)*+{F(x)}="2x";
(70,0)*+{G(x)}="3x";
(70,-20)*+{H(x)}="4x";
{\ar@{|->}^{\rho} "1x"+(20,0);(10,0)};
{\ar_{F(\a)} "2x";"2y"};
{\ar^{\rho(x)} "2x";"3x"};
{\ar_{\rho(y)} "2y";"3y"};
{\ar|-{G(\a)} "3x";"3y"};
{\ar@{=>}^{\rho(\a)} "2y";"3x"};
{\ar^{\s(x)} "3x";"4x"};
{\ar_{\s(y)} "3y";"4y"};
{\ar^{H(\a)} "4x";"4y"};
{\ar@{=>}^{\s(\a)} "3y";"4x"};
\endxy
\ee
for each morphism $y \xleftarrow{\a} x$ in $\mC.$ 
\ed

\bd
\label{eq:verticalcompmod}
The \define{vertical composite} of modifications is
denoted using vertical concatenation as in 
\be
\xy0;/r.15pc/:
(27,0)*+{\mathcal{C}}="2";
(-27,0)*+{\mathcal{D}}="3";
{\ar@/_2.25pc/_{F} "2";"3"};
{\ar@/^2.25pc/^{H} "2";"3"};
{\ar@2{->}(-15,7.5);(-15,-7.5)_{\begin{smallmatrix}\r\\\l\end{smallmatrix}}};
{\ar@{=>}(15,7.5);(15,-7.5)^{\begin{smallmatrix}\s\\\t\end{smallmatrix}}};
{\ar@3{->}(12,0);(-12,0)_{\begin{smallmatrix}m\\n\end{smallmatrix}}};
\endxy
\quad:=\quad
\xy0;/r.15pc/:
(25,0)*+{\mathcal{C}}="2";
(-25,0)*+{\mathcal{D}}="3";
{\ar@/_2.25pc/_{F} "2";"3"};
{\ar|-(0.5){G} "2";"3"};
{\ar@/^2.25pc/^{H} "2";"3"};
{\ar@2{->}(-8,12.5);(-8,1.5)_{\r}};
{\ar@2{->}(8,12.5);(8,1.5)^{\s}};
{\ar@3{->}(6,7);(-6,7)_{m}};
{\ar@{=>}(-8,-1.5);(-8,-12.5)_{\l}};
{\ar@{=>}(8,-1.5);(8,-12.5)^{\t}};
{\ar@3{->}(6,-7);(-6,-7)^{n}};
\endxy
\ee
and is defined by the assignment
\be
\xy0;/r.15pc/:
(-30,0)*+{x}="1";
(40,30)*+{F(x)}="2";
(40,0)*+{G(x)}="3";
(40,-30)*+{H(x)}="4";
{\ar@/^1.75pc/^{\s(x)} "2";"3"};
{\ar@/_1.75pc/_{\rho(x)} "2";"3"};
{\ar@/^1.75pc/^{\t(x)} "3";"4"};
{\ar@/_1.75pc/_{\l(x)} "3";"4"};
{\ar@2{->}(50,15);(30,15)_{m(x)}};
{\ar@2{->}(50,-15);(30,-15)_{n(x)}};
{\ar@{|->}^{\begin{smallmatrix}m\\ n\end{smallmatrix}} "1"+(20,0);(10,0)};
\endxy
\ee
for each object $x$ in $\mC.$ 
\ed

\begin{notation}
\label{not:functor2cat}
Let  $\mC$ and $\mD$ be two (strict) 2-categories. Let
$\mathbf{Fun}(\mC,\mathcal{D})$ 
be the 2-category whose objects are (strict) functors from 
$\mC$ to $\mD,$ 1-morphisms are oplax-natural transformations,
and 2-morphisms are modifications.
The compositions are defined as above (modifications have one additional composition,  which we have not defined). 
\end{notation}

\bd
\label{defn:adjunction}
Let $\mC$ be a (strict) 2-category.
An \define{adjunction} in $\mC$
consists of a pair of objects $x,y$ in $\mC,$
a pair of morphisms 
\be
\xy0;/r.25pc/:
(-10,0)*+{x}="x";
(10,0)*+{y}="y";
{\ar@<1ex>"x";"y"^{f}};
{\ar@<1ex>"y";"x"^{g}};
\endxy
\ee
and a pair of 2-morphisms
\be
\xy0;/r.25pc/:
(10,7.5)*+{x}="x1";
(-10,7.5)*+{x}="x2";
(0,-7.5)*+{y}="y";
{\ar"x1";"x2"_{\id_{x}}};
{\ar"x1";"y"^{f}};
{\ar"y";"x2"^{g}};
{\ar@{=>}(0,7);"y"_(0.35){\h}};
\endxy
\aand
\xy0;/r.25pc/:
(10,-7.5)*+{y}="y1";
(-10,-7.5)*+{y}="y2";
(0,7.5)*+{x}="x";
{\ar"y1";"y2"^{\id_{y}}};
{\ar"y1";"x"_{g}};
{\ar"x";"y2"_{f}};
{\ar@{=>}"x";(0,-7)_(0.65){\e}};
\endxy
\ee
satisfying
\be
\label{eq:zigzag1}
\xy0;/r.25pc/:
(22.5,0)*+{x}="x1";
(7.5,0)*+{y}="y1";
(-7.5,0)*+{x}="x2";
(-22.5,0)*+{y}="y2";
{\ar@/_2.25pc/"x1";"x2"_{\id_{x}}};
{\ar"x1";"y1"|-{f}};
{\ar"y1";"x2"|-{g}};
{\ar"x2";"y2"|-{f}};
{\ar@/^2.25pc/"y1";"y2"^{\id_{y}}};
{\ar@{=>}"y1"+(0,8.5);"y1"^(0.35){\h}};
{\ar@{=>}"x2";"x2"+(0,-8.5)^(0.55){\e}};
\endxy
\quad
=
\quad
\xy0;/r.25pc/:
(10,0)*+{x}="x";
(-10,0)*+{y}="y";
{\ar@/^1.5pc/"x";"y"^{f}};
{\ar@/_1.5pc/"x";"y"_{f}};
{\ar@{=>}(0,5);(0,-5)_{\id_{f}}};
\endxy
\ee
and
\be
\label{eq:zigzag2}
\xy0;/r.25pc/:
(22.5,0)*+{y}="y1";
(7.5,0)*+{x}="x1";
(-7.5,0)*+{y}="y2";
(-22.5,0)*+{x}="x2";
{\ar@/_2.25pc/"x1";"x2"_{\id_{x}}};
{\ar"y1";"x1"|-{g}};
{\ar"x1";"y2"|-{f}};
{\ar"y2";"x2"|-{g}};
{\ar@/^2.25pc/"y1";"y2"^{\id_{y}}};
{\ar@{=>}"y2"+(0,8.5);"y2"^(0.35){\h}};
{\ar@{=>}"x1";"x1"+(0,-8.5)^(0.55){\e}};
\endxy
\quad
=
\quad
\xy0;/r.25pc/:
(10,0)*+{y}="y";
(-10,0)*+{x}="x";
{\ar@/^1.5pc/"y";"x"^{g}};
{\ar@/_1.5pc/"y";"x"_{g}};
{\ar@{=>}(0,5);(0,-5)_{\id_{g}}};
\endxy
.
\ee
Conditions~(\ref{eq:zigzag1})~and~(\ref{eq:zigzag2}) are known as the
\define{zig-zag identities}. 
An adjunction as above is typically
written as a quadruple $(f,g,\h,\e)$
and we say $f$ is \define{left adjoint} 
to $g$ and write $f\dashv g.$ 
However, the notation 
\be
\xy0;/r.25pc/:
(-10,0)*+{x}="x";
(-10,5.25)*{\overset{\id_{x}}{\rotatebox{180}{\Large$\circlearrowleft$}}};
(10,0)*+{y}="y";
(10,-5.25)*{\underset{\id_{y}}{\rotatebox{0}{\Large$\circlearrowleft$}}};
{\ar@<1ex>"x";"y"^{f}};
{\ar@{}@<1ex>"x";"y"}="top";
{\ar@<1ex>"y";"x"^{g}};
{\ar@{}@<1ex>"y";"x"}="bot";
{\ar@{}(0,2);(0,-2)|-{\eta\bot\e}};
\endxy
\ee
may occasionally be employed to more clearly indicate all of the data
in the definition of an adjunction. 
\ed

The usual notion of an adjunction is one where the 2-category is that of 
categories, functors, and natural transformations. One may express 
adjunctions in terms of a universal property in this case. 

\blem
\label{lem:universalpropertyadjunctions}
Let $\mC,\mD$ be categories, let $F,G:\mC\to\mD$ be functors, 
and let $\eta:\id_{\mC}\Rightarrow G\circ F$ and $\e:F\circ G\Rightarrow\id_{\mD}$ be natural transformations so that $(F,G,\eta,\e)$ is an adjunction. Then, the following universal properties hold.
\begin{enumerate}[i.]
\item
\label{lem:adjunctioni}
For every $c\in\mC,$ $d\in\mD,$ and $c\xrightarrow{g}G(d)$ there exists a unique $F(c)\xrightarrow{f}d$ such that 
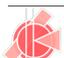
\be
\xy0;/r.25pc/:
(10,7.5)*+{c}="c";
(-10,7.5)*+{G(d)}="Gd";
(0,-7.5)*+{G(F(c))}="GFc";
{\ar"c";"Gd"_{g}};
{\ar"c";"GFc"^{\eta_{c}}};
{\ar"GFc";"Gd"^{G(f)}};
{\ar@{=}(0,4);(0,-1)};
\endxy
\ee
\item
\label{lem:adjunctionii}
For every $c\in\mC,$ $d\in\mD,$ and $F(c)\xrightarrow{f}d$ there exists a unique $c\xrightarrow{g}G(d)$ such that 
\be
\xy0;/r.25pc/:
(10,7.5)*+{F(c)}="Fc";
(-10,7.5)*+{d}="d";
(0,-7.5)*+{F(G(d))}="FGd";
{\ar"Fc";"d"_{f}};
{\ar"Fc";"FGd"^{F(g)}};
{\ar"FGd";"d"^{\e_{d}}};
{\ar@{=}(0,4);(0,-1)};
\endxy
\ee
\end{enumerate}
\elem

\bprf
This is an equivalent definition of an adjunction (cf.~Chapter~IV Section~1 in Mac~Lane~\cite{Ma98}). Given $c\xrightarrow{g}G(d),$ the morphism $F(c)\xrightarrow{f}d$ is given by $f:=\e_{d}\circ F(g).$ Conversely, given $F(c)\xrightarrow{f}d$, the morphism $c\xrightarrow{g}G(d)$ is given by $g:=G(f)\circ\eta_{c}.$ 
\eprf

\bn
\label{prop:adjunctionsinfunctorcategories}
Let  $\mC$ and $\mD$ be two (strict) 2-categories. 
Let 
\be
\label{eq:adjunctionwhole}
\xy0;/r.25pc/:
(-10,0)*+{F}="x";
(-10,5.25)*{\overset{\id_{F}}{\rotatebox{180}{\Large$\circlearrowleft$}}};
(10,0)*+{G}="y";
(10,-5.25)*{\underset{\id_{G}}{\rotatebox{0}{\Large$\circlearrowleft$}}};
{\ar@<1ex>"x";"y"^{\s}};
{\ar@{}@<1ex>"x";"y"}="top";
{\ar@<1ex>"y";"x"^{\rho}};
{\ar@{}@<1ex>"y";"x"}="bot";
{\ar@{}(0,2);(0,-2)|-{\eta\bot\e}};
\endxy
\ee
be an adjunction in the 2-category $\mathbf{Fun}(\mC,\mathcal{D}).$
Then 
\be
\label{eq:adjunctiononeachobject}
\xy0;/r.30pc/:
(-10,0)*+{F(x)}="x";
(-10,5.25)*{\overset{\id_{F(x)}}{\rotatebox{180}{\Large$\circlearrowleft$}}};
(10,0)*+{G(x)}="y";
(10,-5.25)*{\underset{\id_{G(x)}}{\rotatebox{0}{\Large$\circlearrowleft$}}};
{\ar@<1.25ex>"x";"y"^{\s(x)}};
{\ar@{}@<1.25ex>"x";"y"}="top";
{\ar@<1.25ex>"y";"x"^{\rho(x)}};
{\ar@{}@<1.25ex>"y";"x"}="bot";
{\ar@{}(0,2);(0,-2)|-{\eta(x)\bot\e(x)}};
\endxy
\ee
is an adjunction in $\mathcal{D}$ for all objects $x$ in $\mC.$ 
Conversely, let $F,G:\mC\to\mD$ be two functors, let 
$\rho:G\Rightarrow F$ be a pseudo-natural transformation, and let
$\s:\mC_{0}\to\mD_{1}$ and $\eta,\e:\mC_{0}\to\mD_{2}$ be assignments such that (\ref{eq:adjunctiononeachobject}) is an adjunction 
for all $x$ in $\mC.$ Then there exists an extension of $\s$ to an oplax-natural transformation for which $\eta$ together with $\e$ define modifications and (\ref{eq:adjunctionwhole}) is an adjunction in $\mathbf{Fun}(\mC,\mathcal{D}).$ Furthermore, $\s$ is unique up to a canonical isomorphism. 
\en

\br
That $\rho$ be a pseudo-natural transformation and not just a lax- or oplax-natural transformation is explicitly used in the proof. 
\er

\bprf[Proof of Proposition~\ref{prop:adjunctionsinfunctorcategories}] 
{\color{white}{you found me!}}

\noindent
($\Rightarrow$)
The forward direction was proved in the final Remark in \cite{PaGNS}. 

\noindent
($\Leftarrow$)
For the reverse direction, $\s_{\g}:\s_{y}\circ F(\g)\Rightarrow G(\g)\circ\s_{x}$ must be constructed for each morphism $y\xleftarrow{\g}x$ in $\mC.$ It is cumbersome to do this using globular diagrams, so we implement string diagrams to simplify the proof (see \cite{Pa18} for an introduction to string diagrams). By convention, string diagrams will be read from top to bottom and from right to left. Define $\s_{\g}$ by
\be
\label{eq:sigma}
\xy0;/r.25pc/:
(-15,0)*+{G(y)}="Gy";
(0,10)*+{F(y)}="Fy";
(15,0)*+{F(x)}="Fx";
(0,-10)*+{G(x)}="Gx";
{\ar"Fx";"Fy"_{F(\g)}};
{\ar"Fy";"Gy"_{\s_{y}}};
{\ar"Fx";"Gx"^{\s_{x}}};
{\ar"Gx";"Gy"^{G(\g)}};
{\ar@{=>}(0,4);(0,-4)^{\s_{\g}}};
\endxy
\equiv
\vcenter{\hbox{
\begin{tikzpicture}
\newcommand\dist{1.33};
\node[draw,align=center,rounded corners=.4cm,inner sep=8pt] (sg) at (0,0){$\s_{\g}$};
\draw[shorten >= -1pt+\pgflinewidth] (-\dist,{1.5*\dist}) to (-\dist,\dist) node[right]{$\s_{y}$} to [out=-90,in=165] (sg);
\draw[shorten >= -1pt+\pgflinewidth] (\dist,{1.5*\dist}) to (\dist,\dist) node[left]{$F(\g)$} to [out=-90,in=15] (sg);
\draw[shorten <= -1pt+\pgflinewidth]  (sg) to [out=-165,in=90] (-\dist,-\dist) to node[right]{$G(\g)$} (-\dist,{-1.5*\dist});
\draw[shorten <= -1pt+\pgflinewidth] (sg) to [out=-15,in=90] (\dist,-\dist) to node[left]{$\s_{x}$} (\dist,{-1.5*\dist});
\end{tikzpicture}
}}
:=
\vcenter{\hbox{
\begin{tikzpicture}
\newcommand\dist{0.675};
\node[draw,align=center,rounded corners=.4cm,inner sep=7pt] (rg) at (0,0){$\overline{\r_{\g}}$};
\node[draw,align=center,rounded corners=.4cm,inner sep=8pt] (hx) at ({2*\dist},{2*\dist}){$\eta_{x}$};
\node[draw,align=center,rounded corners=.4cm,inner sep=7.5pt] (ey) at ({-2*\dist},{-2*\dist}){$\e_{y}$};
\draw[shorten >= -1pt+\pgflinewidth] (-\dist,{3*\dist}) to [out=-90,in=90] node[right]{$F(\g)$} (-\dist,\dist) to [out=-90,in=165] (rg);
\draw[shorten >= -1pt+\pgflinewidth] (hx) to [out=180,in=90] (\dist,\dist) to [out=-90,in=15] node[right]{$\rho_{x}$} (rg);
\draw[shorten <= -1pt+\pgflinewidth]  (rg) to [out=-165,in=90] node[left]{$\rho_{y}$} (-\dist,-\dist) to [out=-90,in=0] (ey);
\draw[shorten <= -1pt+\pgflinewidth] (rg) to [out=-15,in=90] (\dist,-\dist) to [out=-90,in=90] node[left]{$G(\g)$} (\dist,{-3*\dist});
\draw[shorten <= -1pt+\pgflinewidth] (ey) [out=180,in=-90] to ({-3*\dist},{-\dist}) to [out=90,in=-90] node[right]{$\s_{y}$} ({-3*\dist},{3*\dist});
\draw[shorten <= -1pt+\pgflinewidth] (hx) [out=0,in=90] to ({3*\dist},{\dist}) to [out=-90,in=90] node[left]{$\s_{x}$} ({3*\dist},{-3*\dist});
\end{tikzpicture}
}}
,
\ee
where $\overline{\rho_{\g}}$ is the vertical inverse of $\rho_{\g}.$ 
To verify the oplax-naturality of $\s,$ consider a composable pair $z\xleftarrow{\de}y\xleftarrow{\g}x$ of morphisms in $\mC.$ 
Then, 
\be
\begin{split}
&
\vcenter{\hbox{
\begin{tikzpicture}[scale=0.95, every node/.style={scale=0.95}]
\newcommand\dist{1.0};
\node[draw,align=center,rounded corners=.37cm,inner sep=8pt] (sg) at (\dist,-\dist){\!$\s_{\g}$\!};
\node[draw,align=center,rounded corners=.38cm,inner sep=9pt] (sd) at (-\dist,\dist){\!$\s_{\de}$\!};
\draw[shorten >= -1pt+\pgflinewidth] ({2*\dist},{2*\dist}) to [out=-90,in=90] node[left]{$F(\g)$} ({2*\dist},0) to [out=-90,in=15] (sg);
\draw[shorten >= -1pt+\pgflinewidth] (0,{2*\dist}) node[right]{$F(\de)$} to [out=-90,in=15] (sd);
\draw[shorten >= -1pt+\pgflinewidth] ({-2*\dist},{2*\dist}) node[right]{$\s_{z}$} to [out=-90,in=165] (sd);
\draw[shorten <= -1pt+\pgflinewidth]  (sd) to [out=-15,in=90] (0,0) node[right]{$\s_{y}$} to [out=-90,in=165] (sg);
\draw[shorten <= -1pt+\pgflinewidth] (sg) to [out=-15,in=90] ({2*\dist},{-2*\dist}) node[left]{$\s_{x}$};
\draw[shorten <= -1pt+\pgflinewidth] (sg) to [out=-165,in=90] (0,{-2*\dist}) node[left]{$G(\g)$};
\draw[shorten <= -1pt+\pgflinewidth] (sd) to [out=-165,in=90] ({-2*\dist},0) to [out=-90,in=90] node[right]{$G(\de)$} ({-2*\dist},{-2*\dist});
\end{tikzpicture}
}}
\overset{\text{(\ref{eq:sigma})}}{=\joinrel=\joinrel=}
\vcenter{\hbox{
\begin{tikzpicture}[scale=0.95, every node/.style={scale=0.95}]
\newcommand\dist{0.675};
\node[draw,align=center,rounded corners=.4cm,inner sep=8pt] (rd) at ({-3*\dist},0){$\overline{\r_{\de}}$};
\node[draw,align=center,rounded corners=.4cm,inner sep=8pt] (hy) at ({-1*\dist},{2*\dist}){$\eta_{y}$};
\node[draw,align=center,rounded corners=.4cm,inner sep=9pt] (ez) at ({-5*\dist},{-2*\dist}){$\e_{z}$};
\node[draw,align=center,rounded corners=.4cm,inner sep=8pt] (rg) at ({3*\dist},0){$\overline{\r_{\g}}$};
\node[draw,align=center,rounded corners=.4cm,inner sep=9pt] (hx) at ({5*\dist},{2*\dist}){$\eta_{x}$};
\node[draw,align=center,rounded corners=.4cm,inner sep=8pt] (ey) at ({1*\dist},{-2*\dist}){$\e_{y}$};
\draw[shorten >= -1pt+\pgflinewidth] ({-4*\dist},{3*\dist}) to [out=-90,in=90] node[right]{$F(\de)$} (-4*\dist,\dist) to [out=-90,in=165] (rd);
\draw[shorten >= -1pt+\pgflinewidth] (hy) to [out=180,in=90] (-2*\dist,\dist) to [out=-90,in=15] node[right]{$\rho_{y}$} (rd);
\draw[shorten <= -1pt+\pgflinewidth]  (rd) to [out=-165,in=90] node[left]{$\rho_{z}$} ({-4*\dist},-\dist) to [out=-90,in=0] (ez);
\draw[shorten <= -1pt+\pgflinewidth] (rd) to [out=-15,in=90] ({-2*\dist},-\dist) to [out=-90,in=90] node[left]{$G(\de)$} ({-2*\dist},{-3*\dist});
\draw[shorten <= -1pt+\pgflinewidth] (ez) [out=180,in=-90] to ({-6*\dist},{-\dist}) to [out=90,in=-90] node[right]{$\s_{z}$} ({-6*\dist},{3*\dist});
\draw[shorten >= -1pt+\pgflinewidth] (hy) [out=0,in=90] to ({0},{\dist}) to [out=-90,in=90] node[left]{$\s_{y}$} ({0},{-\dist}) to [out=-90,in=180] (ey);
\draw[shorten >= -1pt+\pgflinewidth] ({2*\dist},{3*\dist}) to [out=-90,in=90] node[right]{$F(\g)$} ({2*\dist},\dist) to [out=-90,in=165] (rg);
\draw[shorten >= -1pt+\pgflinewidth] (hx) to [out=180,in=90] ({4*\dist},\dist) to [out=-90,in=15] node[right]{$\rho_{x}$} (rg);
\draw[shorten <= -1pt+\pgflinewidth]  (rg) to [out=-165,in=90] node[left]{$\rho_{y}$} ({2*\dist},-\dist) to [out=-90,in=0] (ey);
\draw[shorten <= -1pt+\pgflinewidth] (rg) to [out=-15,in=90] ({4*\dist},-\dist) to [out=-90,in=90] node[left]{$G(\g)$} ({4*\dist},{-3*\dist});
\draw[shorten <= -1pt+\pgflinewidth] (hx) [out=0,in=90] to ({6*\dist},{\dist}) to [out=-90,in=90] node[left]{$\s_{x}$} ({6*\dist},{-3*\dist});
\end{tikzpicture}
}}
\\
&\overset{\text{(\ref{eq:zigzag2})}}{=\joinrel=\joinrel=}
\vcenter{\hbox{
\begin{tikzpicture}[scale=0.8, every node/.style={scale=0.8}]
\newcommand\dist{0.675};
\node[draw,align=center,rounded corners=.35cm,inner sep=8pt] (rg) at ({1*\dist},{1*\dist}){$\overline{\r_{\g}}$};
\node[draw,align=center,rounded corners=.33cm,inner sep=8pt] (hx) at ({3*\dist},{3*\dist}){\!$\eta_{x}$\!};
\node[draw,align=center,rounded corners=.35cm,inner sep=8pt] (rd) at ({-1*\dist},{-1*\dist}){$\overline{\r_{\de}}$};
\node[draw,align=center,rounded corners=.29cm,inner sep=7pt] (ez) at ({-3*\dist},{-3*\dist}){$\e_{z}$};
\draw[shorten >= -1pt+\pgflinewidth] (0,{4*\dist}) to [out=-90,in=90] node[right]{$F(\g)$} (0,{2*\dist}) to [out=-90,in=165] (rg);
\draw[shorten >= -1pt+\pgflinewidth] (hx) to [out=180,in=90] ({2*\dist},{2*\dist}) to [out=-90,in=15] node[right]{$\rho_{x}$} (rg);
\draw[shorten <= -1pt+\pgflinewidth]  (rg) to [out=-165,in=90] node[left]{$\rho_{y}$} (0,0) to [out=-90,in=0] (rd);
\draw[shorten <= -1pt+\pgflinewidth] (rg) to [out=-15,in=90] ({2*\dist},0) to [out=-90,in=90] node[left]{$G(\g)$} ({2*\dist},{-4*\dist});
\draw[shorten >= -1pt+\pgflinewidth] (rd) [out=180,in=-90] to ({-2*\dist},{0}) to [out=90,in=-90] node[right]{$F(\de)$} ({-2*\dist},{4*\dist});
\draw[shorten >= -1pt+\pgflinewidth] (hx) [out=0,in=90] to ({4*\dist},{2*\dist}) to [out=-90,in=90] node[left]{$\s_{x}$} ({4*\dist},{-4*\dist});
\draw[shorten <= -1pt+\pgflinewidth]  (rd) to [out=-165,in=90] node[left]{$\rho_{z}$} ({-2*\dist},{-2*\dist}) to [out=-90,in=0] (ez);
\draw[shorten <= -1pt+\pgflinewidth] (rd) to [out=-15,in=90] ({0},{-2*\dist}) to [out=-90,in=90] node[left]{$G(\de)$} ({0},{-4*\dist});
\draw[shorten >= -1pt+\pgflinewidth] (ez) [out=180,in=-90] to ({-4*\dist},{-2*\dist}) to [out=90,in=-90] node[right]{$\s_{z}$} ({-4*\dist},{4*\dist});
\end{tikzpicture}
}}
\underset{\text{for $\rho$}}{\overset{\text{(\ref{eq:oplaxnaturalitycomposition})}}{=\joinrel=\joinrel=}}
\vcenter{\hbox{
\begin{tikzpicture}[scale=0.90, every node/.style={scale=0.90}]
\newcommand\dist{0.675};
\node[draw,align=center,rounded corners=.37cm,inner sep=7pt] (rdg) at (0,0){\!$\overline{\r_{\de\g}}$\!};
\node[draw,align=center,rounded corners=.37cm,inner sep=8pt] (hx) at ({2*\dist},{2*\dist}){\!$\eta_{x}$\!};
\node[draw,align=center,rounded corners=.35cm,inner sep=8pt] (ez) at ({-2*\dist},{-2*\dist}){$\e_{z}$};
\draw[shorten >= -1pt+\pgflinewidth] (-\dist,{3*\dist}) to [out=-90,in=90] node[right]{$F(\de\g)$} (-\dist,\dist) to [out=-90,in=165] (rdg);
\draw[shorten >= -1pt+\pgflinewidth] (hx) to [out=180,in=90] (\dist,\dist) to [out=-90,in=15] node[right]{$\rho_{x}$} (rdg);
\draw[shorten <= -1pt+\pgflinewidth]  (rdg) to [out=-165,in=90] node[left]{$\rho_{z}$} (-\dist,-\dist) to [out=-90,in=0] (ez);
\draw[shorten <= -1pt+\pgflinewidth] (rdg) to [out=-15,in=90] (\dist,-\dist) to [out=-90,in=90] node[left]{$G(\de\g)$} (\dist,{-3*\dist});
\draw[shorten <= -1pt+\pgflinewidth] (ez) [out=180,in=-90] to ({-3*\dist},{-\dist}) to [out=90,in=-90] node[right]{$\s_{z}$} ({-3*\dist},{3*\dist});
\draw[shorten <= -1pt+\pgflinewidth] (hx) [out=0,in=90] to ({3*\dist},{\dist}) to [out=-90,in=90] node[left]{$\s_{x}$} ({3*\dist},{-3*\dist});
\end{tikzpicture}
}}
\overset{\text{(\ref{eq:sigma})}}{=\joinrel=\joinrel=}\!
\vcenter{\hbox{
\begin{tikzpicture}[scale=0.90, every node/.style={scale=0.90}]
\newcommand\dist{1.15};
\node[draw,align=center,rounded corners=.37cm,inner sep=8pt] (sdg) at (0,0){\!$\s_{\de\g}$\!};
\draw[shorten >= -1pt+\pgflinewidth] (-\dist,{1.5*\dist}) to (-\dist,\dist) node[right]{$\s_{z}$} to [out=-90,in=165] (sdg);
\draw[shorten >= -1pt+\pgflinewidth] (\dist,{1.5*\dist}) to (\dist,\dist) node[left]{$F(\de\g)$} to [out=-90,in=15] (sdg);
\draw[shorten <= -1pt+\pgflinewidth]  (sdg) to [out=-165,in=90] (-\dist,-\dist) to node[right]{$G(\de\g)$} (-\dist,{-1.5*\dist});
\draw[shorten <= -1pt+\pgflinewidth] (sdg) to [out=-15,in=90] (\dist,-\dist) to node[left]{$\s_{x}$} (\dist,{-1.5*\dist});
\end{tikzpicture}
}}
\end{split}
\ee
This proves (\ref{eq:oplaxnaturalitycomposition}) for $\s.$ 
For $x\xleftarrow{\id_{x}}x,$ one obtains 
\be
\vcenter{\hbox{
\begin{tikzpicture}[scale=0.8, every node/.style={scale=0.8}]
\newcommand\dist{1.35};
\node[draw,align=center,rounded corners=.3cm,inner sep=8pt] (sg) at (0,0){\!\!$\s_{\id_{x}}$\!\!};
\draw[shorten >= -1pt+\pgflinewidth] (-\dist,{1.5*\dist}) to (-\dist,\dist) node[right]{$\s_{x}$} to [out=-90,in=165] (sg);
\draw[shorten >= -1pt+\pgflinewidth,dashed] (\dist,{1.5*\dist}) to (\dist,\dist) node[left]{$\id_{F(x)}$} to [out=-90,in=15] (sg);
\draw[shorten <= -1pt+\pgflinewidth,dashed]  (sg) to [out=-165,in=90] (-\dist,-\dist) to node[right]{$\id_{G(x)}$} (-\dist,{-1.5*\dist});
\draw[shorten <= -1pt+\pgflinewidth] (sg) to [out=-15,in=90] (\dist,-\dist) to node[left]{$\s_{x}$} (\dist,{-1.5*\dist});
\end{tikzpicture}
}}
\overset{\text{(\ref{eq:sigma})}}{=\joinrel=\joinrel=}
\vcenter{\hbox{
\begin{tikzpicture}[scale=0.8, every node/.style={scale=0.8}]
\newcommand\dist{0.675};
\node[draw,align=center,rounded corners=.25cm,inner sep=5pt] (rg) at (0,0){\!$\overline{\r_{\id_{x}}}$\!};
\node[draw,align=center,rounded corners=.32cm,inner sep=8pt] (hx) at ({2*\dist},{2*\dist}){\!$\eta_{x}$\!};
\node[draw,align=center,rounded corners=.27cm,inner sep=7pt] (ey) at ({-2*\dist},{-2*\dist}){\!$\e_{x}$\!};
\draw[shorten >= -1pt+\pgflinewidth,dashed] (-\dist,{3*\dist}) to [out=-90,in=90] node[right]{$\id_{F(x)}$} (-\dist,\dist) to [out=-90,in=165] (rg);
\draw[shorten >= -1pt+\pgflinewidth] (hx) to [out=180,in=90] (\dist,\dist) to [out=-90,in=15] node[right]{$\rho_{x}$} (rg);
\draw[shorten <= -1pt+\pgflinewidth]  (rg) to [out=-165,in=90] node[left]{$\rho_{x}$} (-\dist,-\dist) to [out=-90,in=0] (ey);
\draw[shorten <= -1pt+\pgflinewidth,dashed] (rg) to [out=-15,in=90] (\dist,-\dist) to [out=-90,in=90] node[left]{$\id_{G(x)}$} (\dist,{-3*\dist});
\draw[shorten >= -1pt+\pgflinewidth] (ey) [out=180,in=-90] to ({-3*\dist},{-\dist}) to [out=90,in=-90] node[right]{$\s_{x}$} ({-3*\dist},{3*\dist});
\draw[shorten >= -1pt+\pgflinewidth] (hx) [out=0,in=90] to ({3*\dist},{\dist}) to [out=-90,in=90] node[left]{$\s_{x}$} ({3*\dist},{-3*\dist});
\end{tikzpicture}
}}
\underset{\text{for $\rho$}}{\overset{\text{(\ref{eq:oplaxidentity})}}{=\joinrel=\joinrel=}}
\vcenter{\hbox{
\begin{tikzpicture}
\newcommand\dist{0.675};
\node[draw,align=center,rounded corners=.3cm,inner sep=5pt] (hx) at ({\dist},{\dist}){$\eta_{x}$};
\node[draw,align=center,rounded corners=.3cm,inner sep=5pt] (ey) at ({-\dist},{-\dist}){$\e_{x}$};
\draw (hx) to [out=180,in=90] (0,0) node[right]{$\rho_{x}$} to [out=-90,in=0] (ey);
\draw (ey) [out=180,in=-90] to ({-2*\dist},{0}) to [out=90,in=-90] node[right]{$\s_{x}$} ({-2*\dist},{2*\dist});
\draw (hx) [out=0,in=90] to ({2*\dist},{0}) to [out=-90,in=90] node[left]{$\s_{x}$} ({2*\dist},{-2*\dist});
\end{tikzpicture}
}}
\overset{\text{(\ref{eq:zigzag1})}}{=\joinrel=\joinrel=}
\vcenter{\hbox{
\begin{tikzpicture}
\newcommand\dist{0.675};
\draw ({0},{2*\dist}) to node[right]{$\s_{x}$} ({0},{-2*\dist});
\end{tikzpicture}
}}
,
\ee
which proves (\ref{eq:oplaxidentity}) for $\s.$ 
Now, consider a 2-morphism $\S:\a\Rightarrow\g$ in $\mC$, as in 
(\ref{eq:2morphismSigma}). Then 
\be
\vcenter{\hbox{
\begin{tikzpicture}[scale=0.75, every node/.style={scale=0.75}]
\newcommand\dist{1.15};
\node[draw,align=center,rounded corners=.3cm,inner sep=8pt] (sg) at (0,0){$\s_{\g}$};
\node[draw,align=center,rounded corners=.3cm,inner sep=6pt] (FS) at ({\dist},{\dist}){$F(\S)$};
\draw[shorten >= -1pt+\pgflinewidth] (-\dist,{2*\dist}) to (-\dist,\dist) node[right]{$\s_{y}$} to [out=-90,in=165] (sg);
\draw[shorten >= -1pt+\pgflinewidth] ({\dist},{2*\dist}) to node[left]{$F(\a)$} (FS);
\draw[shorten >= -1pt+\pgflinewidth] (FS) to [out=-90,in=15] node[left,xshift=0.2cm,yshift=0.2cm]{$F(\g)$} (sg);
\draw[shorten <= -1pt+\pgflinewidth]  (sg) to [out=-165,in=90] (-\dist,-\dist) to node[right]{$G(\g)$} (-\dist,{-1.5*\dist});
\draw[shorten <= -1pt+\pgflinewidth] (sg) to [out=-15,in=90] (\dist,-\dist) to node[left]{$\s_{x}$} (\dist,{-1.5*\dist});
\end{tikzpicture}
}}
\hspace{-4mm}
\overset{\text{(\ref{eq:sigma})}}{=\joinrel=\joinrel=}
\vcenter{\hbox{
\begin{tikzpicture}[scale=0.75, every node/.style={scale=0.75}]
\newcommand\dist{0.675};
\node[draw,align=center,rounded corners=.3cm,inner sep=7pt] (rg) at (0,0){$\overline{\r_{\g}}$};
\node[draw,align=center,rounded corners=.3cm,inner sep=8pt] (hx) at ({2*\dist},{2*\dist}){$\eta_{x}$};
\node[draw,align=center,rounded corners=.3cm,inner sep=7pt] (ey) at ({-2*\dist},{-2*\dist}){$\e_{y}$};
\node[draw,align=center,rounded corners=.3cm,inner sep=6pt] (FS) at ({-\dist},{2*\dist}){$F(\S)$};
\draw[shorten >= -1pt+\pgflinewidth] (-\dist,{3.5*\dist}) node[right]{$F(\a)$} to (FS);
\draw[shorten >= -1pt+\pgflinewidth] (FS) to ({-\dist},{\dist}) node[right]{$F(\g)$} to [out=-90,in=165] (rg);
\draw[shorten >= -1pt+\pgflinewidth] (hx) to [out=180,in=90] (\dist,\dist) to [out=-90,in=15] node[right]{$\rho_{x}$} (rg);
\draw[shorten <= -1pt+\pgflinewidth]  (rg) to [out=-165,in=90] node[left]{$\rho_{y}$} (-\dist,-\dist) to [out=-90,in=0] (ey);
\draw[shorten <= -1pt+\pgflinewidth] (rg) to [out=-15,in=90] (\dist,-\dist) to [out=-90,in=90] node[left]{$G(\g)$} (\dist,{-3*\dist});
\draw[shorten <= -1pt+\pgflinewidth] (ey) [out=180,in=-90] to ({-3*\dist},{-\dist}) to [out=90,in=-90] node[right]{$\s_{y}$} ({-3*\dist},{3.5*\dist});
\draw[shorten <= -1pt+\pgflinewidth] (hx) [out=0,in=90] to ({3*\dist},{\dist}) to [out=-90,in=90] node[left]{$\s_{x}$} ({3*\dist},{-3*\dist});
\end{tikzpicture}
}}
\underset{\text{for $\rho$}}{\overset{\text{(\ref{eq:oplaxnaturalSigma})}}{=\joinrel=\joinrel=}}
\vcenter{\hbox{
\begin{tikzpicture}[scale=0.75, every node/.style={scale=0.75}]
\newcommand\dist{0.675};
\node[draw,align=center,rounded corners=.3cm,inner sep=7pt] (rg) at (0,0){$\overline{\r_{\a}}$};
\node[draw,align=center,rounded corners=.3cm,inner sep=8pt] (hx) at ({2*\dist},{2*\dist}){$\eta_{x}$};
\node[draw,align=center,rounded corners=.3cm,inner sep=7pt] (ey) at ({-2*\dist},{-2*\dist}){$\e_{y}$};
\node[draw,align=center,rounded corners=.3cm,inner sep=6pt] (GS) at ({\dist},{-2*\dist}){$G(\S)$};
\draw[shorten >= -1pt+\pgflinewidth] (-\dist,{3*\dist}) to [out=-90,in=90] node[right]{$F(\a)$} (-\dist,\dist) to [out=-90,in=165] (rg);
\draw[shorten >= -1pt+\pgflinewidth] (hx) to [out=180,in=90] (\dist,\dist) to [out=-90,in=15] node[right]{$\rho_{x}$} (rg);
\draw[shorten <= -1pt+\pgflinewidth]  (rg) to [out=-165,in=90] node[left]{$\rho_{y}$} (-\dist,-\dist) to [out=-90,in=0] (ey);
\draw[shorten <= -1pt+\pgflinewidth] (rg) to [out=-15,in=90] (\dist,-\dist) node[left]{$G(\a)$} to [out=-90,in=90] (GS);
\draw[shorten <= -1pt+\pgflinewidth] (GS) to (\dist,{-3.5*\dist}) node[left]{$G(\g)$};
\draw[shorten <= -1pt+\pgflinewidth] (ey) [out=180,in=-90] to ({-3*\dist},{-\dist}) to [out=90,in=-90] node[right]{$\s_{y}$} ({-3*\dist},{3*\dist});
\draw[shorten <= -1pt+\pgflinewidth] (hx) [out=0,in=90] to ({3*\dist},{\dist}) to [out=-90,in=90] node[left]{$\s_{x}$} ({3*\dist},{-3.5*\dist});
\end{tikzpicture}
}}
\overset{\text{(\ref{eq:sigma})}}{=\joinrel=\joinrel=}
\hspace{-4mm}
\vcenter{\hbox{
\begin{tikzpicture}[scale=0.75, every node/.style={scale=0.75}]
\newcommand\dist{1.15};
\node[draw,align=center,rounded corners=.3cm,inner sep=8pt] (sa) at (0,0){$\s_{\a}$};
\node[draw,align=center,rounded corners=.3cm,inner sep=6pt] (GS) at ({-\dist},{-\dist}){$G(\S)$};
\draw[shorten >= -1pt+\pgflinewidth] (-\dist,{1.5*\dist}) to (-\dist,\dist) node[right]{$\s_{y}$} to [out=-90,in=165] (sa);
\draw[shorten >= -1pt+\pgflinewidth] (GS) to node[right]{$G(\g)$} ({-\dist},{-2*\dist});
\draw[shorten >= -1pt+\pgflinewidth] (\dist,{1.5*\dist}) to (\dist,\dist) node[left]{$F(\a)$} to [out=-90,in=15] (sa);
\draw[shorten <= -1pt+\pgflinewidth]  (sa) to [out=-165,in=90] node[right,xshift=-0.2cm,yshift=-0.2cm]{$G(\a)$} (GS);
\draw[shorten <= -1pt+\pgflinewidth] (sa) to [out=-15,in=90] (\dist,-\dist) to node[left]{$\s_{x}$} (\dist,{-2*\dist});
\end{tikzpicture}
}}
\ee
proves condition~(\ref{eq:oplaxnaturalSigma}) for $\s$ and concludes the proof that $\s$ is an oplax-natural transformation. 
It remains to show that $\e:\id_{F}\Rrightarrow\begin{matrix}[0.6]\s\\\rho\end{matrix}$ and $\eta:\begin{matrix}[0.6]\rho\\\s\end{matrix}\Rrightarrow\id_{G}$ are modifications. 
For this, fix a 1-morphism $y\xleftarrow{\g}x$ in $\mC.$ Then
\be
\hspace{-3mm}
\vcenter{\hbox{
\begin{tikzpicture}[scale=0.9, every node/.style={scale=0.9}]
\newcommand\dist{0.675};
\node[draw,align=center,rounded corners=.35cm,inner sep=7pt] (sg) at (\dist,0){$\s_{\g}$};
\node[draw,align=center,rounded corners=.4cm,inner sep=8pt] (hy) at ({-1*\dist},{2*\dist}){$\eta_{y}$};
\node[draw,align=center,rounded corners=.35cm,inner sep=7pt] (rg) at ({-1*\dist},{-2*\dist}){$\rho_{\g}$};
\draw[shorten >= -1pt+\pgflinewidth] ({2*\dist},{3*\dist}) to [out=-90,in=90] node[left]{$F(\g)$} ({2*\dist},\dist) to [out=-90,in=15] (sg);
\draw[shorten >= -1pt+\pgflinewidth] (hy) to [out=0,in=90] ({0},\dist) node[right]{$\s_{y}$} to [out=-90,in=165] (sg);
\draw[shorten <= -1pt+\pgflinewidth]  (sg) to [out=-165,in=90] (0,-\dist) node[right]{$G(\g)$} to [out=-90,in=0] (rg);
\draw[shorten <= -1pt+\pgflinewidth] (sg) to [out=-15,in=90] ({2*\dist},-\dist) to [out=-90,in=90] node[left]{$\s_{x}$} ({2*\dist},{-3.5*\dist});
\draw[shorten >= -1pt+\pgflinewidth] (hy) to [out=180,in=90] ({-2*\dist},{\dist}) to [out=-90,in=90] node[right]{$\rho_{y}$} ({-2*\dist},{-1*\dist}) to [out=-90,in=165] (rg);
\draw[shorten <= -1pt+\pgflinewidth] (rg) to [out=-15,in=90] ({0},{-3*\dist}) to [out=-90,in=90] node[right]{$\rho_{x}$} ({0},{-3.5*\dist});
\draw[shorten <= -1pt+\pgflinewidth]  (rg) to [out=-165,in=90] ({-2*\dist},{-3*\dist}) to [out=-90,in=90] node[right]{$F(\g)$} ({-2*\dist},{-3.5*\dist});
\end{tikzpicture}
}}
=
\vcenter{\hbox{
\begin{tikzpicture}[scale=0.75, every node/.style={scale=0.8}]
\newcommand\dist{0.675};
\node[draw,align=center,rounded corners=.3cm,inner sep=6pt] (rg) at ({1*\dist},{1*\dist}){$\overline{\r_{\g}}$};
\node[draw,align=center,rounded corners=.3cm,inner sep=7pt] (hx) at ({3*\dist},{3*\dist}){$\eta_{x}$};
\node[draw,align=center,rounded corners=.3cm,inner sep=7pt] (ey) at ({-1*\dist},{-1*\dist}){$\e_{y}$};
\node[draw,align=center,rounded corners=.3cm,inner sep=7pt] (hy) at ({-3*\dist},{3*\dist}){$\h_{y}$};
\node[draw,align=center,rounded corners=.3cm,inner sep=7pt] (rgb) at ({1*\dist},{-3*\dist}){$\rho_{\g}$};
\draw[shorten >= -1pt+\pgflinewidth] (0,{4*\dist}) to [out=-90,in=90] node[right]{$F(\g)$} (0,{2*\dist}) to [out=-90,in=165] (rg);
\draw[shorten >= -1pt+\pgflinewidth] (hx) to [out=180,in=90] ({2*\dist},{2*\dist}) to [out=-90,in=15] node[right]{$\rho_{x}$} (rg);
\draw[shorten <= -1pt+\pgflinewidth]  (rg) to [out=-165,in=90] node[left]{$\rho_{y}$} (0,0) to [out=-90,in=0] (ey);
\draw[shorten <= -1pt+\pgflinewidth] (rg) to [out=-15,in=90] ({2*\dist},0) to [out=-90,in=90] node[left]{$G(\g)$} ({2*\dist},{-2*\dist}) to [out=-90,in=15] (rgb);
\draw[shorten <= -1pt+\pgflinewidth] (ey) [out=180,in=-90] to ({-2*\dist},{0}) to [out=90,in=-90] node[right]{$\s_{y}$} ({-2*\dist},{2*\dist}) to [out=90,in=0] (hy);
\draw[shorten <= -1pt+\pgflinewidth] (hx) [out=0,in=90] to ({4*\dist},{2*\dist}) to [out=-90,in=90] node[left]{$\s_{x}$} ({4*\dist},{-4.5*\dist});
\draw[shorten <= -1pt+\pgflinewidth] (hy) [out=180,in=90] to ({-4*\dist},{2*\dist}) [out=-90,in=90] to node[right]{$\r_{y}$} ({-4*\dist},{-1*\dist}) [out=-90,in=165] to (rgb);
\draw[shorten <= -1pt+\pgflinewidth]  (rgb) to [out=-165,in=90] (0,{-4*\dist}) to node[left]{$F(\g)$} (0,{-4.5*\dist});
\draw[shorten <= -1pt+\pgflinewidth] (rgb) to [out=-15,in=90] ({2*\dist},{-4*\dist}) to node[right]{$\r_{x}$} ({2*\dist},{-4.5*\dist});
\end{tikzpicture}
}}
\overset{\text{(\ref{eq:zigzag2})}}{=\joinrel=\joinrel=}
\hspace{-4mm}
\vcenter{\hbox{
\begin{tikzpicture}[scale=0.75, every node/.style={scale=0.8}]
\newcommand\dist{0.675};
\node[draw,align=center,rounded corners=.3cm,inner sep=7pt] (rg) at ({1*\dist},{1*\dist}){$\overline{\r_{\g}}$};
\node[draw,align=center,rounded corners=.3cm,inner sep=8pt] (hx) at ({3*\dist},{3*\dist}){$\eta_{x}$};
\node[draw,align=center,rounded corners=.3cm,inner sep=8pt] (rgb) at ({1*\dist},{-3*\dist}){$\rho_{\g}$};
\draw[shorten >= -1pt+\pgflinewidth] (0,{4*\dist}) to [out=-90,in=90] node[right]{$F(\g)$} (0,{2*\dist}) to [out=-90,in=165] (rg);
\draw[shorten >= -1pt+\pgflinewidth] (hx) to [out=180,in=90] ({2*\dist},{2*\dist}) to [out=-90,in=15] node[right]{$\rho_{x}$} (rg);
\draw[shorten <= -1pt+\pgflinewidth]  (rg) to [out=-165,in=90] (0,0) to [out=-90,in=90] node[left]{$\rho_{y}$} (0,{-2*\dist}) to [out=-90,in=165] (rgb);
\draw[shorten <= -1pt+\pgflinewidth] (rg) to [out=-15,in=90] ({2*\dist},0) to [out=-90,in=90] node[left]{$G(\g)$} ({2*\dist},{-2*\dist}) to [out=-90,in=15] (rgb);
\draw[shorten <= -1pt+\pgflinewidth] (hx) [out=0,in=90] to ({4*\dist},{2*\dist}) to [out=-90,in=90] node[left]{$\s_{x}$} ({4*\dist},{-4.5*\dist});
\draw[shorten <= -1pt+\pgflinewidth]  (rgb) to [out=-165,in=90] (0,{-4*\dist}) to node[right]{$F(\g)$} (0,{-4.5*\dist});
\draw[shorten <= -1pt+\pgflinewidth] (rgb) to [out=-15,in=90] ({2*\dist},{-4*\dist}) to node[right]{$\r_{x}$} ({2*\dist},{-4.5*\dist});
\end{tikzpicture}
}}
=
\vcenter{\hbox{
\begin{tikzpicture}[scale=0.8, every node/.style={scale=0.8}]
\newcommand\dist{0.675};
\node[draw,align=center,rounded corners=.3cm,inner sep=7pt] (hx) at ({3*\dist},{3*\dist}){$\eta_{x}$};
\draw[shorten >= -1pt+\pgflinewidth] (0,{6*\dist}) to node[right]{$F(\g)$} (0,{-0.5*\dist});
\draw[shorten >= -1pt+\pgflinewidth] (hx) to [out=180,in=90] ({2*\dist},{2*\dist}) to [out=-90,in=90] node[left]{$\rho_{x}$} ({2*\dist},{-0.5*\dist});
\draw[shorten >= -1pt+\pgflinewidth] (hx) [out=0,in=90] to ({4*\dist},{2*\dist}) to [out=-90,in=90] node[left]{$\s_{x}$} ({4*\dist},{-0.5*\dist});
\end{tikzpicture}
}}
,
\ee
where the last equality follows from the fact that $\overline{\rho}$ is the vertical inverse of $\rho.$ 
This proves (\ref{eq:coherencemodification}) for $\eta$, and therefore shows $\eta$ is a modification. 
A similar proof shows $\e$ is a modification. Finally, 
$\s$ is unique up to canonical isomorphism by the uniqueness of adjoints in 2-categories (cf.~Lemma~A.4 in~\cite{PaGNS}). 
\eprf

\br
\label{rmk:simplifiedStinespringproof}
Proposition~\ref{prop:adjunctionsinfunctorcategories} seems to be a useful fact for adjunctions in 2-categories of functors. 
It offers a slightly shorter proof of Theorem~\ref{thm:stinespring}. One merely has to define the functors $\mathbf{Stine}_{\mA}$ and $\rest_{\mA}$ and define the natural transformation $m_{\mA}$. Then one has to show $\begin{matrix}[0.7]\mathbf{Stine}\\\rest\end{matrix}=\id_{\mathbf{OCP}}$ and prove the zig-zag identities for $(\mathbf{Stine}_{\mA},\rest_{\mA},\id,m_{\mA}),$ the last of which were essentially tautologies. Proving oplax-naturality of $\mathbf{Stine}$ and that $m$ is a modification is not necessary thanks to Proposition~\ref{prop:adjunctionsinfunctorcategories}. 
\er
}
\bibliographystyle{plainnat}
\bibliography{Stinespring}

\end{document}